\documentclass{article}
\usepackage{latexsym}
\usepackage[varumlaut]{yfonts}
\usepackage{epsfig}
\usepackage{amsmath}
\usepackage{amssymb}
\usepackage{amsthm}
\usepackage{multirow}
\newtheorem{theorem}{Theorem}[section]
\newtheorem{corollary}[theorem]{Corollary}

\newtheorem{property}[theorem]{Property}
\newtheorem{lemma}[theorem]{Lemma}

\theoremstyle{remark}
\newtheorem{remark}[theorem]{Remark}

\theoremstyle{definition}
\newtheorem{definition}[theorem]{Definition}

\theoremstyle{definition}
\newtheorem{example}[theorem]{Example}

\DeclareMathOperator\Sym{Sym}
\DeclareMathOperator\diag{diag}

\DeclareMathOperator\rk{rk}
\DeclareMathOperator\tr{tr}

\DeclareMathOperator\spa{span}

\DeclareMathOperator\Han{Han}
\DeclareMathOperator\Tri{Tri}
\DeclareMathOperator\Toep{Toep}

\begin{document}

\title{Spectrahedral cones generated by rank 1 matrices}

\author{Roland Hildebrand \thanks{%
WIAS, Mohrenstrasse 39, 10117 Berlin, Germany
({\tt roland.hildebrand@wias-berlin.de}). The author kindly acknowledges support from the ANR G\'eoLMI of the French National Research Agency.}}

\maketitle

\begin{abstract}
Let ${\cal S}_+^n \subset {\cal S}^n$ be the cone of positive semi-definite matrices as a subset of the vector space of real symmetric $n \times n$ matrices. The intersection of ${\cal S}_+^n$ with a linear subspace of ${\cal S}^n$ is called a spectrahedral cone. We consider spectrahedral cones $K$ such that every element of $K$ can be represented as a sum of rank 1 matrices in $K$. We shall call such spectrahedral cones rank one generated (ROG). We show that ROG cones which are linearly isomorphic as convex cones are also isomorphic as linear sections of the positive semi-definite matrix cone, which is not the case for general spectrahedral cones. We give many examples of ROG cones and show how to construct new ROG cones from given ones by different procedures. We provide classifications of some subclasses of ROG cones, in particular, we classify all ROG cones for matrix sizes not exceeding 4. Further we prove some results on the structure of ROG cones. We also briefly consider the case of complex or quaternionic matrices. ROG cones are in close relation with the exactness of semi-definite relaxations of quadratically constrained quadratic optimization problems or of relaxations approximating the cone of nonnegative functions in squared functional systems.
\end{abstract}



\section{Introduction}

Let ${\cal S}^n$ be the real vector space of $n \times n$ real symmetric matrices and ${\cal S}_+^n \subset {\cal S}^n$ the cone of positive semi-definite matrices. The intersection of the cone ${\cal S}_+^n$ with an affine subspace of ${\cal S}^n$ is called a {\it spectrahedron}. Spectrahedra appear as the feasible sets of semi-definite programs and are thus of importance for convex optimization. If the affine subspace happens to be a linear subspace $L \subset {\cal S}^n$, then the intersection $K = L \cap {\cal S}_+^n$ is a {\it spectrahedral cone}. The facial structure of spectrahedra and spectrahedral cones has been studied in \cite{RamanaGoldman95}.

The subject of this contribution are spectrahedral cones $K$ satisfying the following property.

{\property \label{ROGproperty} Every matrix in $K$ can be represented as a sum of rank 1 matrices in $K$. }

\medskip

We shall call such spectrahedral cones {\it rank 1 generated} (ROG). A convex cone in some real vector space will be called a ROG cone if it is linearly isomorphic to a spectrahedral cone possessing Property \ref{ROGproperty}. The corresponding linear isomorphism will define a representation of the ROG cone. Clearly the cone ${\cal S}_+^n$ itself is ROG.

\subsection{Motivation}

In this subsection we elaborate on the role ROG cones play in optimization. The main result is that in some commonly arising situations, semi-definite relaxations of non-convex optimization problems are exact if and only if they lead to conic programs over ROG cones.

The condition of being a ROG spectrahedral cone can equivalently be stated in terms of bounded spectrahedra. Namely, the conic hull $K$ of a bounded spectrahedron $C$ not containing the zero matrix is ROG if and only if $C$ is the convex hull of the rank 1 matrices in $C$. Therefore, if $C$ is a compact section of a ROG spectrahedral cone, then the difficult problem of minimizing a linear function over the nonconvex set of rank 1 matrices in $C$ is equivalent to the easy problem of minimizing this linear function over the bounded spectrahedron $C$.

Property \ref{ROGproperty} is hence in close relation with the exactness of semi-definite relaxations of nonconvex problems in the case when the relaxation is obtained by dropping a rank constraint. Many nonconvex optimization problems which are arising in computational practice fall into this framework, i.e., they can be cast as semi-definite programs with an additional rank constraint. It is this rank constraint which makes the problem nonconvex and difficult to solve. At the same time, dropping the rank constraint provides a convenient way of relaxing the problem to an easily solvable semi-definite program.

A classical example is the MAXCUT problem \cite{GoemansWilliamson}, which can be formulated as the problem of maximizing a linear function over the set of positive semi-definite rank 1 matrices whose diagonal elements all equal 1. By dropping the rank 1 condition, one obtains a semi-definite program which yields an upper bound on the maximum cut. Whether this bound is tight is, however, NP-hard to determine \cite{GareyJohnson79}.

In this context we shall now consider two applications of ROG spectrahedral cones.

\bigskip

{\bf Quadratically constrained quadratic problems.} The most general class of problems which can be formulated as semi-definite programs with an additional rank 1 constraint are the quadratically constrained quadratic problems \cite{RamanaGoldman95},\cite{LMSYZ10}. This class includes also problems with binary decision variables, as the condition $x \in \{a,b\}$ can be cast as the quadratic condition $(x-a)(x-b) = 0$.

A generic quadratically constrained quadratic problem can be written as
\[ \min_{x \in \mathbb R^n} x^TSx\ :\quad x^TA_ix = 0,\ i = 1,\dots,k;\quad x^TBx = 1.
\]
Here $A_1,\dots,A_k;B;S$ are real symmetric $n \times n$ matrices defining the homogeneous quadratic constraints, the inhomogeneous quadratic constraint, and the quadratic cost function, respectively. Introducing the matrix variable $X = xx^T \in {\cal S}_+^n$, we can write the problem as
\begin{equation} \label{nonrelaxed}
\min_{X \in K} \langle S,X \rangle\ :\quad \langle B,X \rangle = 1,\quad \rk X = 1,
\end{equation}
where $K = L \cap {\cal S}_+^n$, and $L \subset {\cal S}^n$ is the linear subspace given by $\{ X \in {\cal S}^n \,|\, \langle A_i,X \rangle = 0 \ \forall\ i = 1,\dots,k \}$. The cone $K$ is hence a linear section of the positive semi-definite matrix cone. Problem \eqref{nonrelaxed} can be relaxed to a semi-definite program by dropping the rank constraint,
\begin{equation} \label{relaxed}
\min_{X \in K} \langle S,X \rangle\ :\quad \langle B,X \rangle = 1.
\end{equation}

Naturally, the question arises when the semi-definite relaxation \eqref{relaxed} obtained from the nonconvex problem \eqref{nonrelaxed} is {\it exact}, i.e., yields the same optimal value as \eqref{nonrelaxed}. In general, this question is NP-hard \cite{RamanaGoldman95}. However, a simple sufficient condition on the spectrahedral cone $K$ is Property \ref{ROGproperty}.

{\lemma \label{QCQPlemma} Let the linear subspace $L \subset {\cal S}^n$ defined above be such that the cone $K = L \cap {\cal S}_+^n$ is rank 1 generated. Then either problems \eqref{nonrelaxed},\eqref{relaxed} are both infeasible, or problem \eqref{relaxed} is unbounded, or problems \eqref{nonrelaxed},\eqref{relaxed} have the same optimal value. }

\begin{proof}
Define the spectrahedron $C = \{ X \in K \,|\, \langle B,X \rangle = 1 \}$. Then the feasible set of problem \eqref{relaxed} is $C$, while that of problem \eqref{nonrelaxed} is $C_1 = \{ X \in C \,|\, \rk X = 1 \}$. If $C = \emptyset$, then both problems are infeasible. Assume that $C \not= \emptyset$. Then $K \not= \{0\}$, and by Property \ref{ROGproperty} every extreme ray of the cone $K$ is generated by a rank 1 matrix. If problem \eqref{relaxed} is bounded, then its optimal value is achieved at an extreme point $X \in C$. Since $X$ generates an extreme ray of $K$, we must have $\rk X = 1$. Thus $X$ is feasible also for problem \eqref{nonrelaxed}, and the optimal value of \eqref{nonrelaxed} is not greater than that of \eqref{relaxed}. But $C_1 \subset C$, and hence the optimal value of \eqref{nonrelaxed} is not smaller than that of \eqref{relaxed}. Therefore both optimal values must coincide.
\end{proof}

In particular, if the spectrahedron $C$ is bounded, then problems \eqref{nonrelaxed} and \eqref{relaxed} are equivalent under the conditions of Lemma \ref{QCQPlemma}.

\bigskip

{\bf Squared functional systems.} Another motivation for the study of ROG spectrahedral cones comes from squared functional systems \cite{NesterovSOS}. Let $\Delta$ be an arbitrary set and $F$ an $n$-dimensional real vector space of real-valued functions on $\Delta$. Choose basis functions $u_1,\dots,u_n \in F$. The {\it squared functional system} generated by these basis functions is the set $\{ u_iu_j \,|\, i,j = 1,\dots,n \}$ of product functions. This system spans another real vector space $V$ of real-valued functions on $\Delta$. Clearly $V$ does not depend on the choice of the basis functions $u_i$, since it is also the linear span of the squares $f^2$, $f \in F$.

Let us define a linear map $\Lambda: V^* \to {\cal S}^n$ and its adjoint $\Lambda^*: {\cal S}^n \to V$ by $\Lambda^*(A) = \sum_{i,j=1}^n A_{ij}u_iu_j$. Here the space ${\cal S}^n$ is identified with its dual by means of the Frobenius scalar product\footnote{The reason for defining the operator $\Lambda$ by virtue of its adjoint $\Lambda^*$ is to stay in line with the notations in \cite{NesterovSOS}. This definition explicitly uses a basis of the space $F$. In a coordinate-free definition, the source space of $\Lambda^*$ should be the space $\Sym^2(F)$ of contravariant symmetric 2-tensors over $F$, and the operator $\Lambda^*$ itself should be defined by linear continuation of the map $f \otimes f \mapsto f^2$, $f \in F$.}. By definition of $V$ the map $\Lambda^*$ is surjective, and hence the map $\Lambda$ is injective.


The {\it sum of squares} (SOS) cone $\Sigma \subset V$, given by the set of all functions of the form $\sum_{k = 1}^N f_k^2$ for $f_1,\dots,f_N \in F$, can be represented as the image $\Lambda^*[{\cal S}_+^n]$ of the positive semi-definite matrix cone and has nonempty interior. The dual $\Sigma^*$ of the SOS cone is given by the set of all dual vectors $w \in V^*$ such that $\Lambda(w) \succeq 0$ \cite[Theorem 17.1]{NesterovSOS}. By injectivity of $\Lambda$ it follows that $\Sigma^*$ is linearly isomorphic to its image $K = \Lambda[\Sigma^*] \subset {\cal S}^n$. This image equals the intersection of ${\cal S}_+^n$ with the linear subspace $L = Im\Lambda$. It follows that $\Sigma^*$ is linearly isomorphic to a spectrahedral cone.

Let $P \subset V$ be the cone of nonnegative functions in $V$. Since every sum of squares of real numbers is nonnegative, we have the inclusion $\Sigma \subset P$. It is then interesting to know when the cones $P$ and $\Sigma$ coincide. The following result shows that the cone $K$ being ROG is a necessary condition.

{\lemma Assume above notations. If $P = \Sigma$, then the spectrahedral cone $K = L \cap {\cal S}_+^n$ is rank 1 generated. }

\begin{proof}
For $x \in \Delta$, define the dual vector $w_x \in V^*$ by $\langle w_x,v \rangle = v(x)$ for all $v \in V$. We first show that for all $x \in \Delta$ the matrix $\Lambda(w_x)$ is contained in the set $K_1 = \{ X \in K \,|\, \rk X \leq 1 \}$.

Fix $x \in \Delta$ and define the vector $s \in \mathbb R^n$ element-wise by $s_i = u_i(x)$, $i = 1,\dots,n$. Then we have for all $A \in {\cal S}^n$ that
\[ \langle \Lambda(w_x),A \rangle = \langle w_x,\Lambda^*(A) \rangle = \sum_{i,j=1}^n A_{ij} \langle w_x,u_iu_j \rangle = \sum_{i,j=1}^n A_{ij} u_i(x)u_j(x) = \langle ss^T,A \rangle.
\]
It follows that $\Lambda(w_x) = ss^T$. Hence the rank of $\Lambda(w_x)$ does not exceed 1. Moreover, we have $\Lambda(w_x) \succeq 0$ and $w_x \in V^*$, and therefore $\Lambda(w_x) \in K$. This proves our claim.

For the sake of contradiction, assume now that $K = \Lambda[\Sigma^*]$ is not ROG. Then there exists a dual vector $y \in \Sigma^*$ such that the matrix $\Lambda(y)$ can be strictly separated from the convex hull of $K_1$. In other words, there exists $A \in {\cal S}^n$ such that $\langle A,\Lambda(y) \rangle < 0$, but $\langle A,X \rangle \geq 0$ for every $X \in K_1$.

Consider the function $q = \Lambda^*(A) \in V$. For every $x \in \Delta$ we have $q(x) = \langle w_x,\Lambda^*(A) \rangle = \langle \Lambda(w_x),A \rangle \geq 0$, because $\Lambda(w_x) \in K_1$. Hence we have $q \in P$. But $\langle q,y \rangle = \langle \Lambda^*(A),y \rangle = \langle A,\Lambda(y) \rangle < 0$, and therefore $y \not\in P^*$.

It follows that $P^* \not= \Sigma^*$ and hence $P \not= \Sigma$. This completes the proof.
\end{proof}

Thus in every squared functional system where the cone of nonnegative functions coincides with the SOS cone $\Sigma$, the dual SOS cone $\Sigma^*$ is linearly isomorphic to a ROG spectrahedral cone. This allows us to construct ROG cones from such squared functional systems. Let us consider two examples.

\begin{itemize}
\item The first example is taken from \cite[Section 3.1]{NesterovSOS}. Here $\Delta = \mathbb R$, and $F$ is the space of all polynomials $q(t)$ of degree not exceeding $n-1$, equipped with the basis of monomials $1,t,\dots,t^{n-1}$. It is well-known that a univariate polynomial $p(t)$ is nonnegative if and only if it is a sum of squares of polynomials $q(t)$ of lower degree. The corresponding ROG cone $K$ is the cone of all Hankel matrices in ${\cal S}_+^n$ and has dimension $2n-1$. We shall denote this cone by $\Han_+^n$. \\
    This result can be generalized to the space of polynomials $q(t,x)$ on $\Delta = \mathbb R \times \mathbb R^m$ which are of degree not exceeding $n-1$ in $t$ and homogeneous of degree 1 in $x = (x_1,\dots,x_m)^T$, equipped with the basis $\{x_1,\dots,x_m,tx_1,\dots,tx_m,\dots,t^{n-1}x_1,\dots,t^{n-1}x_m\}$. The sums of squares representability of the corresponding nonnegative polynomials $p(t,x)$ follows from \cite{Yakubovitch70}. The corresponding ROG cone $K$ is the cone of all block-Hankel matrices in ${\cal S}_+^{nm}$ with $n \times n$ blocks of size $m \times m$ each, and has dimension $\frac{(2n-1)m(m+1)}{2}$. We shall denote this cone by $\Han_+^{n,m}$. Of course, $\Han_+^n = \Han_+^{n,1}$.
\item Let $\Delta = \mathbb R^3$ and let $F$ be the 6-dimensional space of homogeneous quadratic polynomials on $\mathbb R^3$, equipped with the basis $x_1^2,x_2^2,x_3^2,x_2x_3,x_1x_3,x_1x_2$. The space $V$ is then the 15-dimensional space of ternary quartics, and in this space the cone of nonnegative polynomials coincides with the SOS cone \cite{Hilbert}. The corresponding ROG cone $K$ is given by all matrices in ${\cal S}_+^6$ of the form
    \[ A = \begin{pmatrix} a_1 & a_6 & a_5 & a_7 & a_{12} & a_{14} \\ a_6 & a_2 & a_4 & a_{15} & a_8 & a_{10} \\ a_5 & a_4 & a_3 & a_{11} & a_{13} & a_9 \\ a_7 & a_{15} & a_{11} & a_4 & a_9 & a_8 \\ a_{12} & a_8 & a_{13} & a_9 & a_5 & a_7 \\ a_{14} & a_{10} & a_9 & a_8 & a_7 & a_6 \end{pmatrix},\qquad a_1,\dots,a_{15} \in \mathbb R.
    \]
\end{itemize}

\subsection{Outlook}

In this subsection we summarize the contents of the paper.

In Section \ref{sec:preliminaries} we introduce two notions of isomorphisms, a wider one for general convex cones, and the other for spectrahedral cones.

In Section \ref{sec:basic} we study fundamental properties of ROG cones. We establish that the minimal polynomial of a ROG cone, when the latter is viewed as an algebraic interior, is determinantal, and the degree of the cone is given by the maximal rank of the matrices it contains (Subsection \ref{subs:MDP}). In Subsection \ref{subs:FS} we study the facial structure of ROG cones and establish that the rank and the Carath\'eodory number of its elements coincide. In particular, the rank is an invariant of the elements of a ROG cone under linear isomorphisms. In Subsection \ref{subs:II} we prove that the geometry of a ROG cone as a conic convex subset of a real vector space determines its representations as ROG spectrahedral cones uniquely up to isomorphism, which is not true for spectrahedral cones in general.

In Section \ref{sec:construction} we describe different methods to construct ROG cones of higher degree from ROG cones of lower degree. The most simple way is taking direct sums, which is considered in Subsection \ref{subs:direct_sums}. This leads to the notion of simple ROG cones, which are defined as those not representable as a non-trivial direct sum. In Subsections \ref{subs:full_extensions}, \ref{subs:intertwinings} we consider two other ways of constructing ROG cones. The second one can be seen as a generalization of taking direct sums.

In Section \ref{sec:examples} we consider some examples of ROG cones. In Subsection \ref{subs:chordal} we investigate ROG cones defined by conditions of the type that a subset of entries in the representing matrices vanishes. This class of ROG cones is linked to chordal graphs and has been studied in \cite{AHMR88},\cite{PaulsenPowerSmith87}, see also \cite{Laurent01} for a generalization to higher matrix ranks. We show that these cones can be constructed from full matrix cones ${\cal S}_+^k$ by the methods presented in Section \ref{sec:construction}. In Subsection \ref{subs:continuous_family} we construct an example of a continuous family of mutually non-isomorphic ROG cones.

In Section \ref{sec:dim} we consider ROG cones of low codimension (Subsections \ref{subs:codim1}, \ref{subs:codim2}) and simple ROG cones of low dimension (Subsection \ref{subs:dim_bound}).

In Section \ref{sec:discrete} we consider the variety of extreme rays of ROG cones. We show that the discrete part of this variety factors out and does not interfere with the part corresponding to the continuous components.

In Section \ref{sec:small} we give a complete classification of ROG cones for degrees $n \leq 4$ up to isomorphism.

Finally, we briefly consider the case of complex and quaternionic Hermitian matrices in Section \ref{sec:complex}.

We conclude the paper with a summary and an outlook on future work.

\section{Preliminaries} \label{sec:preliminaries}

In this section we formalize some properties of general spectrahedral cones, in particular we rigorously define the notion of isomorphism. Although this notion is widely used implicitly, to our best knowledge is has not yet been explicitly defined in the literature.

\subsection{Notions of isomorphism} \label{subs:notion_iso}

When studying a class of mathematical objects, one has to distinguish between intrinsic properties of the object and those induced by the often necessary coordinate representation. The intrinsic properties are those which are preserved by the isomorphisms of the class. In this paper, the objects under consideration are spectrahedral cones. It is hence necessary to define when spectrahedral cones are considered isomorphic. We shall consider two different notions of isomorphism. The weaker notion forgets about the matricial nature of spectrahedral cones and considers them just as subsets of a real vector space.

{\definition \label{lin_iso} Let $K \subset \mathbb R^n$, $K' \subset \mathbb R^{n'}$ be convex cones. We say that $K$ and $K'$ are {\it linearly isomorphic} if there exists a bijective linear map $l: \spa K \to \spa K'$ such that $l[K] =  K'$. }

\medskip

Here the dimensions $n,n'$ may be different. The dimension of the cone itself is of course invariant under linear isomorphisms.

We shall now take the matricial nature of the spectrahedral cones into account. Naturally, two spectrahedral cones should be considered isomorphic if one cone can be bijectively mapped to the other cone by a coordinate transformation. Such a transformation acts on the matrices in the source space by a map $X \mapsto AXA^T$, where $A$ is a fixed invertible matrix. In order to accommodate different matrix sizes, we allow the matrices in the cone of smaller matrix size to be padded with zeros before the coordinate transformation. Equivalently, we may relax the invertibility condition on the coordinate transformation matrix $A$ and allow it to be rectangular of full column rank. This leads to the following definition.

{\definition \label{def_auto} Let $K \subset {\cal S}_+^n$, $K' \subset {\cal S}_+^{n'}$ be spectrahedral cones, and suppose that $n \leq n'$. We call $K,K'$ {\it isomorphic} if there exists an injective linear map $f: \mathbb R^n \to \mathbb R^{n'}$ with coefficient matrix $A \in \mathbb R^{n' \times n}$ of full column rank such that the induced injective linear map $\tilde f: {\cal S}^n \to {\cal S}^{n'}$ given by $\tilde f: X \mapsto AXA^T$ maps $K$ bijectively onto $K'$. }

\medskip

It is not immediately evident that this notion is well-posed, i.e., defines an equivalence relation. For this we shall need the following reformulation of \cite[Lemma 2.3]{HeltonVinnikov07}.

{\lemma \label{HVlemma} Let $K \subset {\cal S}_+^n$ be a spectrahedral cone, and let $m = \max_{X \in K}\rk X$ be the maximal rank of the matrices in $K$. Then the set $R_{\max} = \{ X \in K \,|\, \rk X = m \}$ equals the (relative) interior of $K$, and the subspace $H(X) = Im X$ is constant over $R_{\max}$. \qed }

In other words, for every spectrahedral cone $K \subset {\cal S}_+^n$ there exists a coordinate system in $\mathbb R^n$ such that the matrices in $K$ are all of the form $X = \begin{pmatrix} M & 0 \\ 0 & 0 \end{pmatrix}$, where $M$ is positive definite if and only if $X$ is in the interior of $K$.

We shall now show that Definition \ref{def_auto} determines an equivalence relation.

\begin{itemize}
\item Reflexivity: By taking $f$ to be the identity map we see that $K$ is isomorphic to itself.
\item Symmetry: If $n \not= n'$, then there is nothing to show. Let now $n = n'$ and let $f: \mathbb R^n \to \mathbb R^{n'}$ with coefficient matrix $A \in \mathbb R^{n' \times n}$ realize an isomorphism between $K$ and $K'$. Then the map $f$ is bijective, and its inverse $f^{-1}: \mathbb R^{n'} \to \mathbb R^n$ with coefficient matrix $A^{-1} \in \mathbb R^{n \times n'}$ realizes an isomorphism between $K'$ and $K$.
\item Transitivity: Let $K \subset {\cal S}_+^n$, $K_1 \subset {\cal S}_+^{n_1}$, $K_2 \subset {\cal S}_+^{n_2}$ be spectrahedral cones and suppose that $K$ is isomorphic to both $K_1$ and $K_2$, the isomorphisms being generated by the injective maps $f_1,f_2$, respectively. Assume that $n_1 \leq n_2$ without loss of generality. Then we have to construct an injective map $f: \mathbb R^{n_1} \to \mathbb R^{n_2}$ which generates an isomorphism between $K_1$ and $K_2$. We have to distinguish several cases:

$n_1 \leq n \leq n_2$: Then $f_1: \mathbb R^{n_1} \to \mathbb R^n$, $f_2: \mathbb R^n \to \mathbb R^{n_2}$, and we may set $f = f_2 \circ f_1$.

$n \leq n_1 \leq n_2$: Then $f_1: \mathbb R^n \to \mathbb R^{n_1}$, $f_2: \mathbb R^n \to \mathbb R^{n_2}$. It is not hard to see that there exists an injective map $f: \mathbb R^{n_1} \to \mathbb R^{n_2}$ such that $f_2 = f \circ f_1$. Any such map generates an isomorphism between $K_1$ and $K_2$.

$n_1 \leq n_2 \leq n$: This is the non-trivial case, since an injective map $f: \mathbb R^{n_1} \to \mathbb R^{n_2}$ satisfying $f_1 = f_2 \circ f$ may not exist. Let $A_i$ be the coefficient matrix of the map $f_i$, $i = 1,2$. Denote the generic elements of $K,K_1,K_2$ by $X,X_1,X_2$, respectively, and set $m = \max_{X \in K} \rk X$. Then we also have $m = \max_{X_i \in K_i} \rk X_i$, $i = 1,2$, because $\rk X_i = \rk A_iX_iA_i^T$. By Lemma \ref{HVlemma} we may assume that all elements of the cones $K,K_1,K_2$ are of the block-diagonal form $\begin{pmatrix} M & 0 \\ 0 & 0 \end{pmatrix}$, where the $m \times m$ matrix $M$ is positive definite if and only if the element is in the interior of the respective cone. Let now $X \in K$ be an element of maximal rank $m$, and let $X_1 \in K_1$, $X_2 \in K_2$ be its preimages under the isomorphisms generated by $f_1,f_2$, respectively. The matrices $X_1,X_2$ are also of rank $m$. Then we have
\[ X_i = \begin{pmatrix} M_i & 0 \\ 0 & 0 \end{pmatrix},\quad X = \begin{pmatrix} M & 0 \\ 0 & 0 \end{pmatrix} = A_i\begin{pmatrix} M_i & 0 \\ 0 & 0 \end{pmatrix}A_i^T,\quad i = 1,2,
\]
where $M,M_1,M_2$ are positive definite $m \times m$ matrices. Partition $A_i = \begin{pmatrix} A_{i,11} & A_{i,12} \\ A_{i,21} & A_{i,22} \end{pmatrix}$, where the block $A_{i,11}$ is of size $m \times m$. Then we obtain $A_{i,11}M_iA_{i,11}^T = M$, $A_{i,21}M_iA_{i,11}^T = 0$, $i = 1,2$. Since $M,M_i$ are invertible, it follows from the first equation that $A_{i,11}$ is invertible. Then the second equation yields $A_{i,21} = 0$. Define the invertible $m \times m$ matrix $A_{11} = A_{2,11}^{-1}A_{1,11}$, let $A_{22}$ be an arbitrary $(n_2 - m) \times (n_1 - m)$ matrix of full column rank, and set $A = \diag(A_{11},A_{22})$. Then $A$ is an $n_2 \times n_1$ matrix of full column rank, and we have by construction $AX_1A^T = X_2$ for all pairs $(X_1,X_2) \in K_1 \times K_2$ such that $A_1X_1A_1^T = A_2X_2A_2^T$. Hence the injective map $f$ defined by the coefficient matrix $A$ generates the sought isomorphism between $K_1,K_2$.
\end{itemize}

Isomorphisms in the sense of Definition \ref{def_auto} are also linear isomorphisms, i.e., in the sense of Definition \ref{lin_iso}. For general spectrahedral cones, the former is, however, a much stronger condition than the latter. For instance, a linear isomorphism between spectrahedral cones in general does not preserve the rank of the matrices in the cone, while the map $X \mapsto AXA^T$ is rank-preserving if $A$ is of full column rank. In particular, isomorphisms in the sense of Definition \ref{def_auto} preserve Property \ref{ROGproperty}. Whether a particular spectrahedral cone is ROG thus depends only on the isomorphism class of this cone.

We may now reformulate Lemma \ref{HVlemma} in terms of Definition \ref{def_auto}.

{\definition We call a spectrahedral cone {\it non-degenerate} if its interior consists of positive definite matrices. }

{\lemma \label{HVlemma2} Let $K$ be a spectrahedral cone and set $m = \max_{X \in K} \rk X$. Then there exists a non-degenerate spectrahedral cone $K' \subset {\cal S}_+^m$ which is isomorphic to $K$. Every non-degenerate spectrahedral cone which is isomorphic to $K$ consists of matrices of size $m \times m$. }

\begin{proof}
Set $H = Im X \subset \mathbb R^n$, where $X \in K$ is a matrix of maximal rank $m$. Identify $H$ with the space $\mathbb R^m$ by introducing an arbitrary coordinate system on $H$. Then the inclusion $i: H \to \mathbb R^n$ defines an injective map $f: \mathbb R^m \to \mathbb R^n$. By Lemma \ref{HVlemma} the interior of $K$ and hence the whole cone $K$ lies in the image of the injective map $\tilde f: {\cal S}^m \to {\cal S}^n$ induced by $f$. The preimage $\tilde f^{-1}[K]$ is then a non-degenerate spectrahedral cone $K' \subset {\cal S}_+^m$. By construction $K'$ is isomorphic to $K$ by virtue of the map $f$. The last assertion of the lemma follows from the fact that the quantity $\max_{X \in K'} \rk X$ is constant over the isomorphism class of $K$, and the matrix size is for a non-degenerate spectrahedral cone always equal to the maximal rank.
\end{proof}

\subsection{Notations} \label{subs:notations}

In this subsection we introduce some notations which simplify the expositions in the next sections.

For $n \in \mathbb N$, define two operators ${\cal L}_n,{\cal F}_n$ from the set of linear subspaces of $\mathbb R^n$ into the set of linear subspaces of ${\cal S}^n$ and the set of faces of the cone ${\cal S}_+^n$, respectively. Let $H \subset \mathbb R^n$ be a linear subspace. Then ${\cal L}_n(H)$, ${\cal F}_n(H)$ will be defined as the linear span and the convex hull of the set $\{ xx^T \in {\cal S}^n \,|\, x \in H \}$, respectively. Note that ${\cal F}_n(H)$ is isomorphic to the cone ${\cal S}_+^{\dim H}$. For a matrix $X \in {\cal S}_+^n$, the smallest face of ${\cal S}_+^n$ containing $X$ is then given by ${\cal F}_n(Im X)$.

In order to indicate the size $n$ of the matrices making up a spectrahedral cone $K$, we shall write $K = L \cap {\cal S}_+^n$ or $K \subset {\cal S}_+^n$, where $L \subset {\cal S}^n$ is a linear subspace. Later in the paper we shall also work with ROG cones as abstract convex conic subsets of a real vector space. They may then have representations in matrix spaces of different sizes.

Let us also define an operator ${\cal H}_n$ from the set of spectrahedral cones $K \subset {\cal S}_+^n$ to the set of linear subspaces of $\mathbb R^n$. For $X \in K$ a matrix of maximal rank, we set ${\cal H}_n(K) = Im X$. By Lemma \ref{HVlemma} the image $Im X$ does not depend on the choice of $X$, and ${\cal H}_n$ is indeed well-defined. The operator ${\cal H}_n$ maps $K$ to the whole space $\mathbb R^n$ if and only if $K$ is non-degenerate.

\section{Basic properties of ROG cones} \label{sec:basic}

In Subsection \ref{subs:MDP} we shall consider ROG cones from the viewpoint of real algebraic geometry. This approach has been quite successful in the study of spectrahedral cones in general, by describing the boundary of these cones as subsets of the zero set of some hyperbolic polynomial. For a certain subclass of spectrahedral cones the minimal such polynomial is precisely the determinant of the matrices making up the cone, and the ROG cones are shown to belong to this subclass. This links the rank of the matrices in the ROG cone to the degree of the cone as an algebraic set.

We then pass on to the facial structure of ROG cones in Subsection \ref{subs:FS}. We show that the facial hierarchy of ROG cones is much more tightly bound to the rank of the matrices in the faces than for general spectrahedral cones. The main result of this subsection is the equality of rank and Carath\'eodory number for the elements of ROG cones, a relation which is familiar and widely used for the full matrix cone ${\cal S}_+^n$.

Subsection \ref{subs:II} contains the main result of the paper. It states that if some convex cone is linearly isomorphic to some ROG cone, then all such ROG cones must be mutually isomorphic as spectrahedral cones. The non-triviality of this assertion comes from the fact that the isomorphism of convex cones as subsets of a real vector space is a much weaker notion than the isomorphism between spectrahedral cones, the latter taking into account also the structure of the matrices making up the cones. While the first two subsections in this section use only elementary tools, in Subsection \ref{subs:II} we will need to consider a certain property of the Pl\"ucker embedding of real Grassmanians. This result does not refer to spectrahedral cones and can be found in the Appendix.

\subsection{Minimal defining polynomial} \label{subs:MDP}

In this subsection we consider the boundary of spectrahedral and ROG cones as a subset of the zero locus of a polynomial. The main goal is to unveil the relation between these polynomials and the determinant of the matrices in the cone. The material in this subsection is basically an application of the theory of algebraic interiors which has been elaborated in \cite{HeltonVinnikov07}.

{\definition \cite[Section 2.2]{HeltonVinnikov07} A closed set $C \subset \mathbb R^m$ is an {\it algebraic interior} if there exists a polynomial $p$ on $\mathbb R^m$ such that $C$ equals the closure of a connected component of the set $\{ x \in \mathbb R^m \,|\, p(x) > 0 \}$. Such a polynomial is called a {\it defining polynomial} of the algebraic interior. }

{\lemma \cite[Lemma 2.1]{HeltonVinnikov07} \label{min_def_poly} Let $C$ be an algebraic interior. Then the defining polynomial $p$ of $C$ with minimal degree is unique up to multiplication by a positive constant. Any other defining polynomial of $C$ is divisible by $p$. }

{\definition The defining polynomial with minimal degree of an algebraic interior $C$ is called {\it minimal defining polynomial}. The {\it degree} of $C$ is defined as the degree of the minimal defining polynomial. }

{\lemma \cite[Theorem 2.2]{HeltonVinnikov07} Every spectrahedron is a convex algebraic interior. }

From Lemma \ref{min_def_poly} it follows that the minimal defining polynomial of a spectrahedral cone is invariant under linear isomorphisms up to a multiplicative positive constant. It also follows that it is homogeneous. Indeed, under a homothety of the cone the minimal defining polynomial transforms to another minimal defining polynomial, which must differ from the original one by a multiplicative positive constant.

For a non-degenerate spectrahedral cone $K \subset {\cal S}_+^n$, a defining polynomial of $K$ is given by the restriction of the determinant in ${\cal S}^n$ to $\spa K$. We shall call this polynomial the {\it determinantal defining polynomial}. Since $\det(AXA^T) = (\det A)^2\cdot\det X$ for square matrices $A,X$, the determinantal defining polynomials of two isomorphic non-degenerate spectrahedral cones differ only by a multiplicative positive constant. We may hence extend the notion of determinantal defining polynomial to the whole isomorphism class of a non-degenerate spectrahedral cone, keeping in mind that for degenerate cones the polynomial is determined only up to a positive constant factor. The degree of the determinantal defining polynomial for general spectrahedral cones $K$ is by virtue of Lemma \ref{HVlemma2} equal to $\max_{X \in K}\rk X$.

In contrast to the minimal defining polynomial, the determinantal defining polynomial is in general not invariant under linear isomorphisms.

{\example \label{exampleL3} Consider the three-dimensional Lorentz cone $L_3 = \{ x = (x_0,x_1,x_2)^T \in \mathbb R^3 \,|\, x_0 \geq \sqrt{x_1^2 + x_2^2} \}$ and two of its spectrahedral representations $K_i = \{ A_i(x) \,|\, x \in L_3 \}$, $i = 1,2$, given by
\[ A_1(x) = \begin{pmatrix} x_0 & x_1 & x_2 \\ x_1 & x_0 & 0 \\ x_2 & 0 & x_0 \end{pmatrix} \in {\cal S}^3,\qquad A_2(x) = \begin{pmatrix} x_0+x_1 & x_2 \\ x_2 & x_0-x_1 \end{pmatrix} \in {\cal S}^2.
\]
It is easily seen that $A_i(x) \succeq 0$ if and only if $x \in L_3$, so $K_1,K_2$ are indeed spectrahedral representations of $L_3$. Clearly they are both non-degenerate. However, the corresponding determinantal defining polynomials are given by $p_{1,\det} = x_0(x_0^2-x_1^2-x_2^2)$, $p_{2,\det} = x_0^2-x_1^2-x_2^2$ and are hence not proportional, while the minimal defining polynomial of $L_3$ is given by $p_{\min} = x_0^2-x_1^2-x_2^2$. This proves that the linearly isomorphic spectrahedral cones $K_1,K_2$ are not isomorphic in the sense of Definition \ref{def_auto}. }

We now come to the main result in this subsection, namely that for ROG cones the determinantal and minimal defining polynomials coincide. Actually, we shall prove this assertion for a somewhat larger subclass of spectrahedral cones\footnote{The extension from ROG cones to the subclass considered in Lemma \ref{Greg_lemma} is due to Gregory Blekherman.}.

{\lemma \label{Greg_lemma} Let $K = L \cap {\cal S}_+^n$ be a non-degenerate spectrahedral cone. Suppose that there exist linearly independent vectors $x_1,\dots,x_n \in \mathbb R^n$ such that $x_ix_i^T \in K$, $i = 1,\dots,n$. Then the determinantal defining polynomial $d$ of $K$ is a minimal defining polynomial. }

\begin{proof}
Denote the linear span of the matrices $x_1x_1^T,\dots,x_nx_n^T \in K$ by $D$, and the intersection $D \cap K$ by $K_D$. We have $D \subset L$, and hence $D \cap {\cal S}_+^n = D \cap L \cap {\cal S}_+^n = K_D$. However, in the coordinates defined by the basis $\{ x_1,\dots,x_n \}$ of $\mathbb R^n$ the subspace $D \subset {\cal S}^n$ is the subspace of diagonal matrices. Hence $K_D = D \cap {\cal S}_+^n$ equals the convex conic hull of $\{ x_1x_1^T,\dots,x_nx_n^T \}$, which in turn is linearly isomorphic to the nonnegative orthant $\mathbb R_+^n$. Moreover, the relative interior of $K_D$ consists of positive definite matrices and is hence contained in the relative interior of $K$. On the other hand, the boundary of $K_D$ is contained in the boundary of $K$ by Lemma \ref{HVlemma}.

Let $p: L \to \mathbb R$ be a minimal defining polynomial of $K$. Since the determinantal defining polynomial $d$ has degree $n$, the degree of $p$ is at most $n$. By Lemma \ref{min_def_poly} $p$ divides $d$. Since $d > 0$ on the relative interior of $K$, we also have $p > 0$ on the relative interior of $K$. Hence $p > 0$ on the relative interior of $K_D$. On the other hand, $p = 0$ on the boundary of $K_D$, because $p = 0$ on the boundary of $K$. Therefore the restriction of $p$ on $D$ is a defining polynomial for the cone $K_D \cong \mathbb R_+^n$.

However, the degree of the algebraic interior $\mathbb R_+^n$ is $n$, and hence $p$ has degree at least $n$. It follows that $\deg p = n$, and $d$ must be a minimal defining polynomial of $K$.
\end{proof}

{\theorem \label{min_eq_det} The determinantal defining polynomial of a ROG spectrahedral cone $K$ is a minimal defining polynomial. }

\begin{proof}
Recall that both the determinantal and the minimal defining polynomial is invariant under isomorphisms up to multiplication by a positive constant. We may then assume without loss of generality that $K \subset {\cal S}_+^n$ is non-degenerate, otherwise we pass to an isomorphic non-degenerate ROG spectrahedral cone by virtue of Lemma \ref{HVlemma2}.

Let $X \in K$ be positive definite. Since $K$ is ROG, there exist vectors $x_1,\dots,x_N \in \mathbb R^n$ such that $X = \sum_{i=1}^N x_ix_i^T$ and $x_ix_i^T \in K$ for all $i = 1,\dots,N$. By virtue of $X \succ 0$ the linear span of $\{ x_1,\dots,x_N \}$ equals $\mathbb R^n$. In particular, among the $x_i$ there are $n$ linearly independent vectors, let these be $x_1,\dots,x_n$. The proof is concluded by application of Lemma \ref{Greg_lemma}.
\end{proof}

We may now link the degree of a ROG spectrahedral cone to the rank of the matrices in the cone.

{\corollary \label{deg_K} The degree of a ROG spectrahedral cone $K$ is given by $\deg K = \max_{X \in K}\rk X$. }

\begin{proof}
The right-hand side of the relation is the degree of the determinantal defining polynomial of $K$, while the left-hand side is the degree of the minimal defining polynomial. The assertion now follows from Theorem \ref{min_eq_det}.
\end{proof}

In this subsection we have shown that for the subclass of ROG spectrahedral cones, two different kinds of associated polynomials coincide. These are on the one hand the {\it determinantal defining polynomial}, which is determined by the matricial structure of the elements of the cone and is invariant under isomorphisms in the sense of Definition \ref{def_auto}, and on the other hand the {\it minimal defining polynomial}, which is a notion from real algebraic geometry and is invariant under linear isomorphisms in the sense of Definition \ref{lin_iso}.

\subsection{Facial structure and rank} \label{subs:FS}

In this subsection we study the facial structure of ROG cones and its connection to the Carath\'eodory number. This allows us to establish a number of representation lemmas which bound the number of rank 1 matrices which enter the sum in Property \ref{ROGproperty}. The results in this subsection follow from properties of the facial structure of general spectrahedral cones and from standard convex analysis arguments.

We shall call an element of a cone $K$ {\it extreme} if it generates an extreme ray of $K$.

{\lemma \label{lem_extreme} Let $K \subset {\cal S}_+^n$ be a ROG spectrahedral cone. Then the set of extreme elements of $K$ is given by $\{ X \in K \,|\, \rk X = 1 \}$. }

\begin{proof}
Since $K$ is ROG, every $X \in K$ with $\rk X > 1$ can be represented as sum of elements $X_i \in K$ of rank 1. Hence such $X$ cannot be extreme. On the other hand, every $X \in K$ with $\rk X = 1$ generates an extreme ray of ${\cal S}_+^n$. Extremality in $K$ for such $X$ follows immediately.
\end{proof}

Let us recall the results of \cite{RamanaGoldman95} on the facial structure of general spectrahedral cones. Let $K = L \cap {\cal S}_+^n$ be a spectrahedral cone. Then the faces of $K$ are given by the intersections of $L$ with the faces of ${\cal S}_+^n$ \cite[Theorem 1]{RamanaGoldman95}, see also \cite[Prop.\ 2.1]{Ycart82}. In particular, the kernel of the matrices $X \in K$ is constant over the relative interior of each face of $K$, and every face of $K$ is exposed \cite[Corollary 1]{RamanaGoldman95}. It follows that the faces of spectrahedral cones are also spectrahedral cones.

The smallest face of $K = L \cap {\cal S}_+^n$ containing a matrix $X \in K$ is given by the intersections $L \cap {\cal F}_n(Im X) = L \cap {\cal S}_+^n \cap {\cal L}_n(Im X) = K \cap {\cal L}_n(Im X)$, because ${\cal F}_n(Im X)$ is the smallest face of ${\cal S}_+^n$ containing $X$. The smallest face of ${\cal S}_+^n$ containing $K$ is given by ${\cal F}_n(Im X)$, where $X$ is an arbitrary matrix in the interior of $K$. Here the operators ${\cal F}_n,{\cal L}_n$ are defined in Subsection \ref{subs:notations}.

{\lemma \label{face_ROG} Every face of a ROG cone is a ROG cone. }

\begin{proof}
Let $K = L \cap {\cal S}_+^n$ be a ROG cone and $K' \subset K$ a face of $K$. Then there exists a face $F$ of ${\cal S}_+^n$ such that $K' = L \cap F$, e.g., $F = ({\cal F}_n \circ {\cal H}_n)(K')$. Let $X \in K'$ be an arbitrary nonzero matrix. Since $X \in K$ and $K$ is ROG, there exist rank 1 matrices $X_1,\dots,X_N \in K$ such that $X = \sum_{i=1}^N X_i$. At the same time, $X \in F$. Since $F$ is a face of ${\cal S}_+^n$, the rank 1 matrices $X_i \in {\cal S}_+^n$ must also be elements of this face. It follows that $X_i \in K'$, and $X$ can be represented as sum of rank 1 matrices in $K'$. Thus $K'$ is ROG.
\end{proof}

{\definition \cite[p.59]{GulerTuncel98} 
Let $K \subset \mathbb R^m$ be a closed pointed convex cone. The {\it Carath\'eodory number} $\kappa(x)$ of a point $x \in K$ is the minimal number $k$ such that there exist extreme elements $x_1,\dots,x_k$ of $K$ satisfying $x = \sum_{i=1}^k x_i$.

The {\it Carath\'eodory number} $\kappa(K)$ of the cone $K$ is the maximum of $\kappa(x)$ over $x \in K$. }

{\lemma \label{Cara_X_upper} Let $K = L \cap {\cal S}_+^n$ be a spectrahedral cone. The Carath\'eodory number of $X \in K$ satisfies $\kappa(X) \leq \rk X$. }

\begin{proof}
We proceed by induction. If $\rk X \leq 1$, then by virtue of Lemma \ref{lem_extreme} we trivially have $\kappa(X) = \rk X$. Suppose the relation $\kappa(X) \leq \rk X$ is proven for $\rk X \leq k-1$, and let $X \in K$ with $\rk X = k \geq 2$.

Without loss of generality we may assume $n = k$, otherwise we replace $K$ by $K_X = L \cap {\cal F}_n(Im X)$, the minimal face of $K$ which contains $X$. Neither the rank nor the Carath\'eodory number of $X$ will change by this substitution of the ambient cone, but now ${\cal F}_n(Im X) \cong {\cal S}_+^k$ and $K_X$ can be seen as a spectrahedral cone defined by $k \times k$ matrices.

Then the boundary of $K$ consists of matrices $Y$ with $\rk Y < k = n$, and hence $\kappa(Y) < k$ by the induction hypothesis. Let $E \in K$ an extreme element of $K$, normalized such that $\tr E = \tr X$. Consider the compact line segment $l$ which is defined by the intersection of $K$ with the affine line passing through $X$ and $E$. Since $X$ is in the interior of $K$, it is also in the interior of the segment $l$. One endpoint of $l$ is given by $E$, while the other one is some matrix $Y \in \partial K$. Then there exists $\lambda \in (0,1)$ such that $X = \lambda E + (1-\lambda)Y$. Hence $\kappa(X) \leq \kappa(E) + \kappa(Y) \leq 1 + (k-1) = k$. This completes the proof.
\end{proof}

{\lemma \label{Cara_X} Let $K \subset {\cal S}_+^n$ be a ROG spectrahedral cone. The Carath\'eodory number of $X \in K$ is given by $\kappa(X) = \rk X$. }

\begin{proof}
We have $\kappa(X) \geq \rk X$, because by virtue of Lemma \ref{lem_extreme} all generators of extreme rays of $K$ have rank 1, and a matrix $X$ cannot be the sum of less that $\rk X$ matrices of rank 1. On the other hand, $\kappa(X) \leq \rk X$ by Lemma \ref{Cara_X_upper}.
\end{proof}

{\corollary \label{Cara_deg} The Carath\'eodory number of a ROG cone equals its degree. }

\begin{proof}
The claim follows immediately from Lemma \ref{Cara_X} and Corollary \ref{deg_K}.
\end{proof}

{\corollary \label{cor_diag1} Let $K \subset {\cal S}_+^n$ be a ROG spectrahedral cone, and let $X \in K$ be an element of rank $k$. Then there exist rank 1 matrices $X_i = x_ix_i^T \in K$, $i = 1,\dots,k$, such that $X = \sum_{i=1}^k X_i$ and the vectors $x_1,\dots,x_k$ are linearly independent. }

\begin{proof}
The Corollary is a consequence of Lemmas \ref{lem_extreme} and \ref{Cara_X}.
\end{proof}

{\corollary \label{cor_diag2} Let $K \subset {\cal S}_+^n$ be a ROG spectrahedral cone of degree $d = \deg K$. Then there exist $d$ linearly independent vectors $r_1,\dots,r_d \in \mathbb R^n$ such that $r_ir_i^T \in K$ for $i = 1,\dots,d$. }

\begin{proof}
The claim follows from Corollaries \ref{deg_K} and \ref{cor_diag1}.
\end{proof}

As a consequence, we have the following stand-alone result on the diagonalization of matrices in a ROG cone.

{\lemma \label{lem_diag} Let $K \subset {\cal S}_+^n$ be a ROG spectrahedral cone, and let $X \in K$ be an element of rank $k$. Then there exists a basis of $\mathbb R^n$ such that in the corresponding coordinates we have $X = \diag(1,\dots,1,0,\dots,0)$, and all diagonal matrices of the form $\diag(d_1,\dots,d_k,0,\dots,0)$, where $d_i \geq 0$, $i = 1,\dots,k$, are in $K$. }

\begin{proof}
By Corollary \ref{cor_diag1} there exist linearly independent vectors $x_1,\dots,x_k \in \mathbb R^n$ such that $X_i = x_ix_i^T \in K$, $i = 1,\dots,k$, and $X = \sum_{i=1}^k X_i$. Extend the set $\{ x_1,\dots,x_k \}$ to a basis of $\mathbb R^n$, then in the coordinates defined by this basis we have $X = \diag(1,\dots,1,0,\dots,0)$.

Moreover, for all $d_1,\dots,d_k \geq 0$ we have $\sum_{i=1}^k d_i x_ix_i^T \in K$, and in the coordinates defined above this matrix has the form $\diag(d_1,\dots,d_k,0,\dots,0)$.
\end{proof}

In this subsection we have shown that the equality between the Carath\'eodory number of a matrix $X \in {\cal S}_+^n$ and its rank which is trivially valid for the cone ${\cal S}_+^n$ extends to ROG spectrahedral cones in general. This allowed us to establish the representation result Corollary \ref{cor_diag1} and the existence result Corollary \ref{cor_diag2}. Lemma \ref{face_ROG} asserts that the subclass of ROG cones is closed under the operation of taking faces, a result which is known to be valid for general spectrahedral cones too.

\subsection{Isomorphisms and linear isomorphisms} \label{subs:II}

This subsection contains the main structural result of the paper, namely that the notion of linear isomorphism of convex cones from Definition \ref{lin_iso} and the notion of isomorphism of spectrahedral cones from Definition \ref{def_auto} coincide on the class of ROG spectrahedral cones. More precisely, we show that if two ROG spectrahedral cones are linearly isomorphic, then they are also isomorphic as spectrahedral cones. The proof requires Lemma \ref{Pluecker_tech}, which follows from an auxiliary result on the image of the Pl\"ucker embedding of real Grassmanians. These are provided in the Appendix.

{\theorem \label{theorem_iso} Let $K \subset {\cal S}_+^n,K' \subset {\cal S}_+^{n'}$ be linearly isomorphic ROG cones. Then they are also isomorphic in the sense of Definition \ref{def_auto}. }

\begin{proof}
By Lemma \ref{HVlemma2} we may assume without loss of generality that the cones $K,K'$ are non-degenerate. By Lemma \ref{min_def_poly} the degrees of $K$ and $K'$ coincide, and hence by Corollary \ref{deg_K} also the maximal ranks $n,n'$ of matrices in the cones $K,K'$, respectively, coincide. We may thus consider both cones as linear sections of the matrix cone ${\cal S}_+^n$. Let $L,L' \subset {\cal S}^n$ be the linear hulls of $K,K'$, respectively, and set $m = \dim L = \dim L'$.

Let $\tilde f: L \to L'$ be a bijective linear map realizing the linear isomorphism between $K$ and $K'$. Let $x_1,\dots,x_m \in \mathbb R^n$ be such that the set $\{ x_ix_i^T \,|\, i = 1,\dots,m \}$ forms a basis of $L$. This is possible because $K$ is a ROG cone. Note also that the vectors $x_1,\dots,x_m$ span the whole space $\mathbb R^n$ because $L$ contains non-singular matrices. For every $i = 1,\dots,m$ we have that $x_ix_i^T$ is an extreme element of $K$ by Lemma \ref{lem_extreme}. Its image $\tilde f(x_ix_i^T)$ must then be an extreme element of $K'$, and again by Lemma \ref{lem_extreme} a positive semi-definite rank 1 matrix. Hence there exist nonzero vectors $y_i \in \mathbb R^n$ such that $\tilde f(x_ix_i^T) = y_iy_i^T$. Moreover, the images $\tilde f(x_ix_i^T)$ form a basis of $L'$, because $\tilde f$ is a bijection, and the vectors $y_1,\dots,y_m$ span $\mathbb R^n$ because $L'$ contains non-singular matrices.

Denote the determinantal polynomial on ${\cal S}^n$ by $d$. Then $p = d|_L$, $p' = d|_{L'}$ are the determinantal defining polynomials of $K,K'$, respectively. By Theorem \ref{min_eq_det} both $p,p'$ are minimal defining polynomials. By Lemma \ref{min_def_poly} there exists a positive constant $c > 0$ such that $p = c\cdot(p' \circ \tilde f)$.

Then the conditions of Lemma \ref{Pluecker_tech} are fulfilled and by this lemma there exists an automorphism of ${\cal S}^n$, given by the map $X \mapsto AXA^T$ for some non-singular matrix $A$, which coincides with $\tilde f$ on $L$ and hence maps $K$ bijectively onto $K'$. This completes the proof.
\end{proof}

Theorem \ref{theorem_iso} states that the geometry of a ROG cone as a subset of real space determines its representations as linear sections of a positive semi-definite matrix cone uniquely up to isomorphisms in the sense of Definition \ref{def_auto}. Of course, this does not preclude the existence of other, nonisomorphic, representations as a spectrahedral cone, but in these the cone will not be ROG. For instance, the spectrahedral cone $K_2$ in Example \ref{exampleL3} is ROG, but the linearly isomorphic spectrahedral cone $K_1$ is not. In the sequel, when we speak of a representation of a ROG cone, we will always mean a spectrahedral representation where the cone is ROG.

%
%
%

\section{Construction of new ROG cones from given ones} \label{sec:construction}

In this section we consider several ways to construct ROG spectrahedral cones of higher degree from given ones. By iterating these procedures, one may construct ROG cones of arbitrarily high complexity.

\subsection{Direct sums} \label{subs:direct_sums}

In this subsection we consider direct sums of ROG cones and introduce the notion of a simple ROG cone\footnote{We propose to reserve the notion {\it irreducible} for ROG cones $K \subset {\cal S}_+^n$ such that the real projective variety defined by the set $\{ x \in \mathbb R^n \,|\, xx^T \in K \}$ is irreducible.}, which is a cone that cannot be represented as a non-trivial direct sum. First we shall consider general spectrahedral cones and pinpoint the difficulties which are associated to the notion of direct sum. Then we show that for ROG cones, the situation looks much more favorable.

Recall that in Subsection \ref{subs:notion_iso} we considered two different notions of isomorphisms of spectrahedral cones. The notion of linear isomorphism, given in Definition \ref{lin_iso}, disregarded the matricial structure of the cone, while the second, stronger notion in Definition \ref{def_auto} took it into account. Similarly we may define the notion of direct sum of spectrahedral cones in different ways, first disregarding the matricial nature of the cones and then taking it into account.

{\definition \label{direct_sum_gen} Let $K \subset \mathbb R^n$, $K' \subset \mathbb R^{n'}$ be convex cones. Their {\it direct sum} $K \oplus K'$ is defined as the set $\{ (x,x') \in \mathbb R^n \oplus \mathbb R^{n'} \,|\, x \in K,\ x' \in K' \}$. }

{\definition \label{direct_sum_spectrahedral} Let $K \subset {\cal S}_+^n$, $K' \subset {\cal S}_+^{n'}$ be spectrahedral cones. Their {\it direct sum} $K \oplus K'$ is defined as the set $\{ \diag(X,X') \in {\cal S}_+^{n+n'} \,|\, X \in K,\ X'\in K'\}$. }

\medskip

Note that for the direct sum $K \oplus K'$ of spectrahedral cones $K \subset {\cal S}_+^n$, $K' \subset {\cal S}_+^{n'}$, the ambient vector space in Definition \ref{direct_sum_gen} is the product ${\cal S}^n \times {\cal S}^{n'}$, while in Definition \ref{direct_sum_spectrahedral} it is the matrix space ${\cal S}^{n+n'}$. However, the former can be naturally regarded as a subspace of the latter, namely the subspace of appropriately partitioned block-diagonal matrices. With this identification, both definitions obviously lead to the same result. The direct sum $K \oplus K'$ is also a spectrahedral cone, because a block-diagonal matrix is positive semi-definite if and only if all blocks are. All these considerations naturally extend to an arbitrary number of factors.

Let now a spectrahedral cone $K$ be isomorphic to a direct sum $K_1 \oplus K_2$ of convex cones in the sense of Definition \ref{direct_sum_gen}. Note that the factors $K_1,K_2$ are faces of $K$, and faces of spectrahedral cones are spectrahedral cones. Therefore the factors $K_1,K_2$ inherit from $K$ the structure of spectrahedral cones, and we may consider their direct sum in the sense of Definition \ref{direct_sum_spectrahedral}. By construction this direct sum is linearly isomorphic to the original cone $K$, but it turns out that it is not necessarily isomorphic to $K$ in the sense of Definition \ref{def_auto}.

{\example \label{dir_sum_example} Consider the spectrahedral cones $K = \{ A(x_1,x_2) \,|\, (x_1,x_2)^T \in \mathbb R_+^2 \}$, $K_i = \{ A_i(x_i) \,|\, x_i \geq 0 \}$, $i = 1,2$, given by
\[ A(x_1,x_2) = \begin{pmatrix} x_1 & 0 & 0 \\ 0 & x_1+x_2 & 0 \\ 0 & 0 & x_2 \end{pmatrix},\quad A_1(x_1) = \begin{pmatrix} x_1 & 0 & 0 \\ 0 & x_1 & 0 \\ 0 & 0 & 0 \end{pmatrix},\quad A_2(x_2) = \begin{pmatrix} 0 & 0 & 0 \\ 0 & x_2 & 0 \\ 0 & 0 & x_2 \end{pmatrix}.
\]
The cone $K$ is linearly isomorphic to a direct sum of two copies of $\mathbb R_+$, and the cones $K_1,K_2$ are its faces corresponding to the factors. However, their direct sum in the sense of Definition \ref{direct_sum_spectrahedral} is given by $K_1 \oplus K_2 = \{ \diag(x_1,x_1,0,0,x_2,x_2) \,|\, (x_1,x_2)^T \in \mathbb R_+^2 \}$. We have $\max_{X \in K} \rk X = 3$, $\max_{X \in K_1 \oplus K_2} \rk X = 4$, and $K$ cannot be isomorphic to $K_1 \oplus K_2$. }

We shall show that ROG cones behave nicely in this respect, and such a situation cannot arrive. However, the main result of this subsection is much stronger. We show that a ROG cone consisting of block-diagonal matrices must be a direct sum whose decomposition into factors is given by the block partition. This is not true for general spectrahedral cones, as the example above demonstrates. The cones $K_1,K_2$ both consist of block-diagonal matrices, but neither of them decomposes into factors.

{\lemma \label{lem_product1} Let $K_1,\dots,K_m$ be ROG cones of degrees $n_1,\dots,n_m$. Then their direct sum $K = \oplus_{k=1}^m K_k$ is also a ROG cone, of degree $n = \sum_{k=1}^m n_k$. The cone $K$ is isomorphic to a block-diagonal non-degenerate ROG cone $K'$ with block sizes $n_1,\dots,n_m$, such that block $k$ defines a non-degenerate ROG cone $K_k'$ which is isomorphic to $K_k$. }

\begin{proof}
Let $X_k \in K_k$ be arbitrary and let $X = \diag(X_1,\dots,X_m) \in K$. Since the factor cones $K_k$ are ROG, every $X_k$ decomposes into a sum of rank 1 matrices $r_{k,j} \in K_k$, $j = 1,\dots,\eta_k$. For every such rank 1 matrix $r_{k,j}$, the matrix $R_{k,j} = \diag(0,\dots,0,r_{k,j},0,\dots,0)$ is a rank 1 matrix in $K$, where the non-zero block is located at position $k$. Then $X = \sum_{k=1}^m \sum_{j=1}^{\eta_k} R_{k,j}$ is a rank 1 decomposition of $X$ as required in Property \ref{ROGproperty}. Since $X \in K$ was an arbitrary element, the cone $K$ is ROG.

By Corollary \ref{deg_K} and Lemma \ref{HVlemma2}, for every $k = 1,\dots,m$ there exists a non-degenerate ROG cone $K_k' \subset {\cal S}_+^{n_k}$ which is isomorphic to $K_k$. Their direct sum $K' = \oplus_{k=1}^m K_k'$ is isomorphic to $K$ and has the required block structure. It is also non-degenerate, because a block-diagonal matrix with all blocks being positive definite is itself positive definite. Hence $K$ has degree $n = \sum_{k=1}^m n_k$ by Corollary \ref{deg_K}.
\end{proof}

{\corollary \label{ROG_product} Let $K_1,\dots,K_m$ be convex cones and $K = \oplus_{k=1}^m K_k$ their direct sum in the sense of Definition \ref{direct_sum_gen}. Then $K$ has a ROG spectrahedral representation if and only if all cones $K_k$ have ROG spectrahedral representations. }

\begin{proof}
If the $K_k$ have ROG representations, then their direct sum in the sense of Definition \ref{direct_sum_spectrahedral} is ROG by Lemma \ref{lem_product1} and defines the required representation of $K$.

Let now $K$ have a ROG representation. Each of the cones $K_k$ is linearly isomorphic to a face of $K$, and this face defines a ROG representation of $K_k$ by Lemma \ref{face_ROG}.
\end{proof}

We shall now apply the powerful Theorem \ref{theorem_iso} to the situation in the preceding corollary. Since all ROG spectrahedral representations of $K$ are isomorphic, they are in particular isomorphic to the simple block-diagonal one. Hence every such representation has itself a particularly simple structure.

{\lemma \label{lem_product2} Let the cone $K = \oplus_{k=1}^m K_k$ be a direct sum of lower-dimensional cones in the sense of Definition \ref{direct_sum_gen} and suppose that $K$ possesses a ROG spectrahedral representation $K \subset {\cal S}_+^n$. Then there exists a direct sum decomposition $\mathbb R^n = \oplus_{k=1}^m H_k$ into subspaces of dimensions $\dim H_k \geq \deg K_k$ such that the intersection $F_k = {\cal L}_n(H_k) \cap K$ is linearly isomorphic to $K_k$ for all $k = 1,\dots,m$, and $K = \sum_{k=1}^m F_k$. If $n = \deg K$, then $\dim H_k = \deg K_k$ for $k = 1,\dots,m$. }

\begin{proof}
By Corollary \ref{ROG_product} each factor cone $K_k$ is ROG, denote its degree by $n_k$.

By Lemma \ref{lem_product1} we have $\deg K = \sum_{k=1}^m n_k$ and $K$ possesses a block-diagonal representation as a linear section of ${\cal S}_+^{\deg K}$ with block sizes $n_k$, such that block $k$ defines a representation of the factor cone $K_k$. By Theorem \ref{theorem_iso} the original representation of $K$ as a linear section of ${\cal S}_+^n$ is isomorphic to this block-diagonal representation. Let $f: \mathbb R^{\deg K} \to \mathbb R^n$ be the injective linear map from Definition \ref{def_auto} which defines the isomorphism, and denote by $H \subset \mathbb R^n$ the image of $f$. The map $f$ then puts the direct sum decomposition of $\mathbb R^{\deg K}$ defined by the block structure of the block-diagonal representation in correspondence to some direct sum decomposition $H = \oplus_{k=1}^m H'_k$, where $\dim H'_k = \deg K_k$. Let $\mathbb R^n = \oplus_{k=1}^m H_k$ be an arbitrary direct sum decomposition such that $H'_k \subset H_k$ for all $k = 1,\dots,m$. By construction this decomposition has the required properties.

If $n = \deg K$, then $f$ is bijective, and $H_k = H'_k$ is the only possible choice for $H_k$. It follows that $\dim H_k = \deg K_k$ in this case.
\end{proof}

We now come to the main result of this subsection. The key idea is that in a block-diagonal rank 1 matrix, only one block is non-zero. Hence a block-diagonal ROG cone must be a direct sum.

{\lemma \label{lem_product} Let $K = L \cap {\cal S}_+^n$ be a ROG cone. Let $\mathbb R^n = H_1 \oplus \dots \oplus H_m$ be a direct sum decomposition of $\mathbb R^n$ and suppose that $L \subset \sum_{k=1}^m {\cal L}_n(H_k)$. Then $K$ is the sum of the ROG cones $K_k = K \cap {\cal L}_n(H_k)$, $k = 1,\dots,m$, and is canonically isomorphic to their direct sum. }

\begin{proof}
First note that the cones $K_k$ are faces of $K$ and hence indeed ROG cones by Lemma \ref{face_ROG}. Moreover, the sum $\sum_{k=1}^m K_k$ is canonically isomorphic to the direct sum $\oplus_{k=1}^m K_k$, because we have $\dim(\sum_{k=1}^m {\cal L}_n(H_k)) = \sum_{k=1}^m \dim {\cal L}_n(H_k)$.

Clearly $\sum_{k=1}^m K_k \subset K$, because $K_k \subset K$ for all $k$ and $K$ is a convex cone.

Let now $X \in K$ be arbitrary. By Property \ref{ROGproperty} there exist rank 1 matrices $X_i \in K$, $i = 1,\dots,N$, such that $X = \sum_{i=1}^N X_i$. Now for every $i$ we have $X_i \in L \subset \sum_{k=1}^m {\cal L}_n(H_k)$. Since $X_i$ is rank 1, there must exist $k_i \in \{1,\dots,m\}$ such that $X_i \in {\cal L}_n(H_{k_i})$. It follows that $X_i \in {\cal L}_n(H_{k_i}) \cap K = K_{k_i}$. Therefore $X \in \sum_{k=1}^m K_k$, and hence $K \subset \sum_{k=1}^m K_k$.

Thus we get $K = \sum_{k=1}^m K_k$, which completes the proof.
\end{proof}

{\definition \label{def_simple} We call a ROG cone $K$ {\it simple} if it is not isomorphic to a nontrivial direct sum of lower-dimensional cones. }

By Corollary \ref{ROG_product} it is irrelevant for this definition whether we suppose a direct sum decomposition in the sense of Definition \ref{direct_sum_gen} or in the sense of Definition \ref{direct_sum_spectrahedral} here. Note that the decomposition of a cone $K$ into simple factor cones in the sense of Definition \ref{direct_sum_gen} is unique up to a permutation of the factors. The next result provides another criterion for simplicity.

{\lemma \label{simple_char} A non-degenerate ROG cone $K \subset {\cal S}_+^n$ is simple if and only if there does not exist a nontrivial decomposition $\mathbb R^n = H_1 \oplus \dots \oplus H_m$ such that $\spa K \subset \sum_{k=1}^m {\cal L}_n(H_k)$. }

\begin{proof}
If $K$ is not simple, then Lemma \ref{lem_product2} applies with a nontrivial direct sum decomposition of $\mathbb R^n$. The assertion of this lemma then implies $\spa K \subset \sum_{k=1}^m {\cal L}_n(H_k)$.

On the other hand, if there exists a nontrivial decomposition $\mathbb R^n = H_1 \oplus \dots \oplus H_m$ such that $\spa K \subset \sum_{k=1}^m {\cal L}_n(H_k)$, then by Lemma \ref{lem_product} $K$ is isomorphic to the direct sum of the ROG cones $K_k = K \cap {\cal L}_n(H_k)$, $k = 1,\dots,m$. Since $K \subset {\cal S}_+^n$ is non-degenerate, we have $\deg K_k = \dim H_k > 0$ for all $k$, and the direct sum is nontrivial.
\end{proof}

{\lemma Let $K \subset {\cal S}_+^n$ be a non-degenerate ROG cone. Then there exists a unique (up to a permutation of factors) direct sum decomposition $\mathbb R^n = H_1 \oplus \dots \oplus H_m$ such that $K$ is the sum of the faces $K_k = {\cal L}_n(H_k) \cap K$, and such that the factor cones $K_k$ are simple. }

\begin{proof}
The claim of the lemma follows from Lemma \ref{lem_product2}, applied to the unique decomposition of $K$ into simple factor cones in the sense of Definition \ref{direct_sum_gen}, and the fact that the subspaces $H_k$ are uniquely determined by the faces $K_k$ representing the factor cones.
\end{proof}

Thus there are two different criteria that allow to check whether a ROG cone $K$ is not simple. On the one hand, one may consider the geometric decomposition of $K$ into factor cones, disregarding the matricial structure. On the other hand, one has the algebraic criterion whether in a non-degenerate representation, $K$ is contained in the sum of two complementary faces of the ambient matrix cone.

\subsection{Full extensions} \label{subs:full_extensions}

In this subsection we consider a method for constructing larger spectrahedral cones $K$ from smaller cones $K'$, such that $K$ is the preimage of $K'$ under a specific linear projection. The interest in this procedure is motivated by the fact that it preserves Property \ref{ROGproperty}, and hence allows to construct larger ROG cones (both in terms of degree and of dimension) from smaller ones. The asserted invariance of Property \ref{ROGproperty} is the main result of this subsection.

Before describing the construction formally, we shall provide a non-formal explanation and an example. Let $n'< n$ and consider the partition $\mathbb R^n = H \oplus E$, where $H$ is spanned by the first $n'$ and $E$ by the last $n - n'$ basis vectors, and the corresponding partition of the matrices $X \in {\cal S}^n$ into four blocks $X_{ij}$, $i,j = 1,2$. Given a spectrahedral cone $K' \subset {\cal S}_+^{n'}$, we define a cone $K \subset {\cal S}_+^n$ as the set of all matrices $X \in {\cal S}_+^n$ such that the upper left $n' \times n'$ submatrix of $X$ is an element of $K'$. The crucial observation is that the cone $K \subset {\cal S}_+^n$ is a spectrahedral cone, obtained by imposing linear conditions on the block $X_{11}$ only.

{\example \label{full_ext_example} Let $n' = 2$, $n = 3$, and let $K'$ be the direct sum ${\cal S}_+^1 \oplus {\cal S}_+^1$, or equivalently, the cone of diagonal positive semi-definite matrices. The cone $K$ is then given by
\[ K = \left\{ \begin{pmatrix} a_1 & 0 & a_3 \\ 0 & a_2 & a_4 \\ a_3 & a_4 & a_5 \end{pmatrix} \succeq 0 \,|\, a_1,\dots,a_5 \in \mathbb R \right\}.
\]
By a permutation of rows and columns one obtains that $K$ is isomorphic to the cone of positive semi-definite tri-diagonal matrices. Both cones $K',K$ are ROG. }

\medskip

The procedure above relies on a very specific decomposition $\mathbb R^n = H \oplus E$, determined by the chosen basis of $\mathbb R^n$. It is not hard to see that the essential objects linking the cones $K'$ and $K$ are the projection $\pi: \mathbb R^n \to \mathbb R^{n'}$ which truncates the last $n-n'$ elements of vectors $x \in \mathbb R^n$ and the induced projection $\tilde\pi = \pi \otimes \pi: {\cal S}^n \to {\cal S}^{n'}$ which assigns to a matrix $X$ its subblock $X_{11}$. In a coordinate-free setting, we thus have to depart from an arbitrary projection. By an appropriate choice of coordinates in $\mathbb R^{n'}$ and $\mathbb R^n$ we may, however, always achieve the block partition described above.

Let $n' < n$ be positive integers and consider a surjective linear map $\pi: \mathbb R^n \to \mathbb R^{n'}$. The map $\pi$ induces a surjective linear map $\tilde\pi = \pi \otimes \pi: {\cal S}^n \to {\cal S}^{n'}$, acting on rank 1 matrices by $\tilde\pi: xx^T \mapsto \pi(x)\pi(x)^T$. Let $E \subset \mathbb R^n$ be the kernel of $\pi$, then the kernel of $\tilde\pi$ is given by the linear subspace $L_E \subset {\cal S}^n$ spanned by all matrices of the form $xy^T+yx^T$, $x \in \mathbb R^n$, $y \in E$. For any subspace $H \subset \mathbb R^n$ which is complementary to $E$, there exists a unique right inverse $\mu_H$ of $\pi$ such that $Im\,\mu_H = H$. The map $\mu_H$ is injective and generates an injective map $\tilde\mu_H = \mu_H \otimes \mu_H: {\cal S}^{n'} \to {\cal S}^n$. The map $\tilde\mu_H$ is the unique right inverse to $\tilde\pi$ with image ${\cal L}_n(H)$. The following result formalizes the assertions made about the cones $K',K$ above.

{\lemma \label{full_wellD} Assume above notations, and let $K' \subset {\cal S}_+^{n'}$ be a spectrahedral cone. Then the intersection $K = \tilde\pi^{-1}[K'] \cap {\cal S}_+^n$ is also a spectrahedral cone, and $K'= \tilde\pi[K]$. Moreover, for any subspace $H \subset \mathbb R^n$ which is complementary to $E$, the map $\tilde\mu_H$ is an isomorphism between $K'$ and the face ${\cal L}_n(H) \cap K$ of $K$. }

\begin{proof}
By definition of $K$ we have $\tilde\pi[K] \subset K'$.

Let now $H$ be a complementary subspace to $E$. By definition $\tilde\mu_H$ is an isomorphism between $K'$ and its image $\tilde\mu_H[K'] = \tilde\pi^{-1}[K'] \cap {\cal L}_n(H)$. However, $\tilde\mu_H[K'] \subset {\cal S}_+^n$, and hence $\tilde\mu_H[K'] = {\cal S}_+^n \cap \tilde\pi^{-1}[K'] \cap {\cal L}_n(H) = K \cap {\cal L}_n(H)$.

Moreover, since $\tilde\mu_H$ is a right inverse of $\tilde\pi$, we have $\tilde\pi[{\cal L}_n(H) \cap K] = K'$ and hence $K' \subset \tilde\pi[K]$, which proves that $K' = \tilde\pi[K]$.

There exists a linear subspace $L' \subset {\cal S}^{n'}$ such that $K' = L' \cap {\cal S}_+^{n'}$. The preimage $L = \tilde\pi^{-1}[L']$ is then a subspace of ${\cal S}^n$. We claim that $K = L \cap {\cal S}_+^n$.

Since $K' \subset L'$, we have $\tilde\pi^{-1}[K'] \subset L$ and hence $K \subset L \cap {\cal S}_+^n$.

On the other hand, $\tilde\pi[L \cap {\cal S}_+^n] \subset \tilde\pi[L] \cap \tilde\pi[{\cal S}_+^n] = L' \cap {\cal S}_+^{n'} = K'$, and hence $L \cap {\cal S}_+^n \subset \tilde\pi^{-1}[K']$. This yields $L \cap {\cal S}_+^n \subset K$. Thus $K = L \cap {\cal S}_+^n$ is a spectrahedral cone.
\end{proof}

Since any two linear surjections $\pi_1,\pi_2: \mathbb R^n \to \mathbb R^{n'}$ are conjugated by an automorphism of the source space $\mathbb R^n$, the isomorphism class of $K$ depends only on $K'$ and $n$, but not on the concrete realization of the surjection $\pi$. This implies that the structure of the cone $K$ is fully determined by the smaller cone $K'$ and motivates the following definition.

{\definition \label{def_full} Let $n' < n$ be integers and let $K' \subset {\cal S}_+^{n'}$, $K \subset {\cal S}_+^n$ be spectrahedral cones. We call $K$ a {\it full extension} of $K'$ if there exists a surjective linear map $\pi: \mathbb R^n \to \mathbb R^{n'}$ such that $K = \tilde\pi^{-1}[K'] \cap {\cal S}_+^n$, where $\tilde\pi = \pi \otimes \pi$. }

We now consider when a given spectrahedral cone $K \subset {\cal S}_+^n$ is a full extension of some smaller cone $K'$. The following result gives a sufficient condition.

{\lemma \label{full_suff} Let $K \subset {\cal S}_+^n$ be a spectrahedral cone and suppose there exist a subspace $E \subset \mathbb R^n$ of dimension $k > 0$ and a subspace $L \subset {\cal S}^n$ such that $K = L \cap {\cal S}_+^n$ and $L_E \subset L$. Then for every linear surjective map $\pi: \mathbb R^n \to \mathbb R^{n-k}$ with kernel $E$, the conditions of Definition \ref{def_full} are satisfied with $n'= n - k$, $K' = \tilde\pi[K]$. }

\begin{proof}
We have $\tilde\pi^{-1}[\tilde\pi[K]] = K + \ker\tilde\pi = (L \cap {\cal S}_+^n) + L_E$.

Since $L_E \subset L$, we get $(L \cap {\cal S}_+^n) + L_E \subset L$ and hence $((L \cap {\cal S}_+^n) + L_E) \cap {\cal S}_+^n \subset K$.

On the other hand, we trivially have $K = L \cap {\cal S}_+^n \subset ((L \cap {\cal S}_+^n) + L_E) \cap {\cal S}_+^n$.

Hence $K = \tilde\pi^{-1}[\tilde\pi[K]] \cap {\cal S}_+^n$, which is what we had to show.
\end{proof}

The conditions in Lemma \ref{full_suff} are also necessary for $K$ to be a full extension. In order to see this, one may choose $E$ equal to the kernel of the projection $\pi$ and $L$ as in the proof of Lemma \ref{full_wellD}.

We now come to the main result of this subsection.

{\theorem \label{basic_full} Let $n' < n$ and let the spectrahedral cone $K \subset {\cal S}_+^n$ be a full extension of the spectrahedral cone $K' \subset {\cal S}_+^{n'}$. Then $K$ is ROG if and only if $K'$ is ROG.  }

\begin{proof}
Let $\pi: \mathbb R^n \to \mathbb R^{n'}$ be a surjective linear map satisfying the conditions of Definition \ref{def_full}. By an appropriate choice of the coordinates in $\mathbb R^n$ and $\mathbb R^{n'}$ we may achieve that the projection $\pi$ truncates the last $n - n'$ elements of vectors $x \in \mathbb R^n$, and the map $\tilde\pi$ takes a matrix $X \in {\cal S}^n$ to its upper left subblock $X_{11}$, as in the explanation at the beginning of this subsection.

Let $K$ be ROG. By Lemma \ref{full_wellD} $K'$ is isomorphic to a face of $K$, and by Lemma \ref{face_ROG} faces of ROG cones are also ROG. This proves that $K'$ is ROG.

Suppose now that $K'$ is a ROG cone. Let $X \in K$ be arbitrary. Then $X_{11} = \tilde\pi(X) \in K'$ by Lemma \ref{full_wellD}, and there exist nonzero vectors $v_1,\dots,v_N \in \mathbb R^{n'}$ such that $X_{11} = \sum_{i=1}^N v_iv_i^T$ and $v_iv_i^T \in K'$ for all $i$. Let $V$ be the $n' \times N$ matrix formed of the column vectors $v_i$. The condition $X \succeq 0$ implies that the columns of the block $X_{12}$ are in the image of $X_{11} = VV^T$. Therefore there exists a $k \times N$ matrix $W$ such that $X_{12} = VW^T$. Let the columns of $W$ be $w_1,\dots,w_N \in \mathbb R^k$. We then have the representation
\[ X = \begin{pmatrix} V \\ W \end{pmatrix} \begin{pmatrix} V \\ W \end{pmatrix}^T + \begin{pmatrix} 0 & 0 \\ 0 & X_{22} - WW^T \end{pmatrix} = \sum_{i=1}^N \begin{pmatrix} v_i \\ w_i \end{pmatrix} \begin{pmatrix} v_i \\ w_i \end{pmatrix}^T + \begin{pmatrix} 0 & 0 \\ 0 & X_{22} - WW^T \end{pmatrix}.
\]
Denote the rank 1 matrix $\begin{pmatrix} v_i \\ w_i \end{pmatrix} \begin{pmatrix} v_i \\ w_i \end{pmatrix}^T$ by $U_i$, $i = 1,\dots,N$. By construction $U_i \succeq 0$ and its upper left $n' \times n'$ submatrix $v_iv_i^T$ is an element of $K'$. Hence $U_i \in \tilde\pi^{-1}[v_iv_i^T] \cap {\cal S}_+^n \subset K$ for all $i$. The $k \times k$ matrix $X_{22} - WW^T$ is the Schur complement of $X_{11}$ in $X$ and is hence positive semi-definite. It can then be written as a sum $\sum_{j=1}^{N'} z_jz_j^T$ with $z_j \in \mathbb R^k$. The rank 1 matrices $Z_j = \begin{pmatrix} 0 & 0 \\ 0 & z_jz_j^T \end{pmatrix}$ are also in $K$, and hence $X = \sum_{i=1}^N U_i + \sum_{j=1}^{N'} Z_j$ is a sum of rank 1 matrices in $K$. This shows that $K$ is also ROG and proves the other direction of the equivalence.
\end{proof}

Given a ROG cone $K' \subset {\cal S}_+^{n'}$, Theorem \ref{basic_full} allows us to construct ROG cones $K$ consisting of matrices of size $n$ for any $n > n'$.

It is not hard to see that under the conditions of Definition \ref{def_full}, $\max_{X \in K} \rk X = n - n' + \max_{X' \in K'} \rk X'$. In particular, the cone $K \subset {\cal S}_+^n$ is non-degenerate if and only if $K' \subset {\cal S}_+^{n'}$ is non-degenerate. If $K',K$ are ROG cones, then by Corollary \ref{deg_K} $\deg K = n - n' + \deg K'$. Note also that the full extension of a ROG cone is always simple.

\subsection{Intertwinings} \label{subs:intertwinings}

In this subsection we present a way to construct new ROG cones from pairs of given ROG cones of smaller degree. Let us start with an example.

{\example \label{ex_intertwining} Consider the cones $K_1 = \Han_+^3$, $K_2 = {\cal S}_+^2$. A composite cone $K \subset {\cal S}_+^4$ can be constructed from these two, consisting of positive semi-definite matrices of the form
\[ \begin{pmatrix} a_1 & a_2 & a_3 & 0 \\ a_2 & a_3 & a_4 & 0 \\ a_3 & a_4 & a_5 & a_6 \\ 0 & 0 & a_6 & a_7 \end{pmatrix},\qquad a_1,\dots,a_7 \in \mathbb R.
\]
One recognizes the cone $K_1$ in the upper left $3 \times 3$ subblock of $K$, and the cone $K_2$ in the lower right $2 \times 2$ subblock. The subblocks intersect in a smaller central subblock of size 1. Both $K_1$ and $K_2$ are canonically isomorphic to faces of $K$ corresponding to the subblocks, and $K$ is equal to the {\it sum} of these faces. }

For generic instances of this construction, the upper left subblock defining $K_1$, the lower right subblock defining $K_2$, and the central subblock representing their intersection can be of any (compatible) sizes. If the central subblock has size zero, then $K$ is the {\it direct} sum $K_1 \oplus K_2$. For non-trivial central subblocks the composite cone $K$ is a projection of the direct sum $K_1 \oplus K_2$. The crucial condition that forces the sum $K_1 + K_2$ to be a spectrahedral cone is that the intersection $K_1 \cap K_2$, which naturally has non-zero elements only in the central subblock, is isomorphic to a full matrix cone and contains {\it every} positive semi-definite matrix which has non-zero elements only in the central subblock. In the example above this condition holds because the variable $a_5$ parameterizes the whole central subblock and does not appear anywhere else. In general, there will be many non-equivalent ways to combine two given cones $K_1,K_2$. Our interest in this procedure is based on the fact that the composite cone is ROG if and only if the smaller cones $K_1,K_2$ are.

We shall now formally define how given spectrahedral cones $K_1,K_2$ can be composed to yield the cone $K$. We shall work in a coordinate-free setting, independent of a specific choice of coordinates, or equivalently, a specific block decomposition of the involved matrices. One should keep in mind, however, that by an appropriate coordinate change one can always achieve the block-structured situation described above. For ease of notation, we introduce the following definition.

{\definition Let $K \subset {\cal S}_+^n$ be a spectrahedral cone. We call a face $F$ of $K$ {\it full} if it is also a face of ${\cal S}_+^n$. The number $k = \max_{X \in F} \rk X$ is called the {\it rank} of the face. }


\medskip

For $i = 1,2$, let $K_i \subset {\cal S}_+^{n_i}$ be spectrahedral cones possessing full faces $F_i \subset K_i$ of rank $k$. Let $H_i = {\cal H}_{n_i}(F_i) \subset \mathbb R^{n_i}$ be the $k$-dimensional linear subspaces corresponding to these faces. Let $\iota_i: \mathbb R^k \to \mathbb R^{n_i}$ be injective linear maps such that $Im\,\iota_i = H_i$. Consider the $k$-dimensional subspace $N = \{ (\iota_1(x),-\iota_2(x)) \,|\, x \in \mathbb R^k \}$ of the direct sum $\mathbb R^{n_1} \oplus \mathbb R^{n_2}$. Set $n = n_1 + n_2 - k$ and identify $\mathbb R^n$ with the quotient space $(\mathbb R^{n_1} \oplus \mathbb R^{n_2})/N$. Let $f_i: \mathbb R^{n_i} \to \mathbb R^n$ be the natural embeddings of the factors $\mathbb R^{n_i}$ into $(\mathbb R^{n_1} \oplus \mathbb R^{n_2})/N$, i.e., $f_1: x_1 \mapsto (x_1,0)+N$, $f_2: x_2 \mapsto (0,x_2)+N$. Then the $f_i$ are injective linear maps such that $f_1 \circ \iota_1 = f_2 \circ \iota_2$. Let $\tilde\iota_i = \iota_i \otimes \iota_i: {\cal S}^k \to {\cal S}^{n_i}$ and $\tilde f_i = f_i \otimes f_i: {\cal S}^{n_i} \to {\cal S}^n$ be the injective maps induced by $\iota_i,f_i$, respectively. Then we have also $\tilde f_1 \circ \tilde\iota_1 = \tilde f_2 \circ \tilde\iota_2$. 

The construction in the preceding paragraph ensures that the images $\tilde f_1[K_1],\tilde f_2[K_2]$ are isomorphic to $K_1,K_2$, respectively, and that they intersect in a full face of rank $k$, namely $\tilde f_1[F_1] = \tilde f_2[F_2]$.

{\definition \label{def_intertwining} Assume above conditions. We call the cone $K = \tilde f_1[K_1] + \tilde f_2[K_2] \subset {\cal S}_+^n$ an {\it intertwining} of the cones $K_1,K_2$ along the full faces $F_i$. }

{\remark The cone $K$ can be seen as the projection of the direct sum $K_1 \oplus K_2$ along the linear subspace generated by the set $\{ (\tilde\iota_1(xx^T),-\tilde\iota_2(xx^T)) \,|\, x \in \mathbb R^k \} \subset {\cal S}^{n_1} \times {\cal S}^{n_2}$. }

We have to show that $K$ is indeed a spectrahedral cone. The following observation is crucial for the proof.

{\lemma \label{3partition} Let $M = \begin{pmatrix} A & B & 0 \\ B^T & C & D \\ 0 & D^T & E \end{pmatrix}$ be a block-partitioned positive semi-definite matrix. Then there exists a decomposition $C = C_1 + C_2$ such that the matrices $\begin{pmatrix} A & B \\ B^T & C_1 \end{pmatrix}$, $\begin{pmatrix} C_2 & D \\ D^T & E \end{pmatrix}$ are positive semi-definite. }

\begin{proof}
The Schur complement of $A$ in $M$ is given by $\begin{pmatrix} C - B^TA^{\dagger}B & D \\ D^T & E \end{pmatrix}$ and is positive semi-definite. Here $A^{\dagger}$ is the pseudo-inverse of $A$, which is also positive semi-definite. Setting $C_1 = B^TA^{\dagger}B$, $C_2 = C - B^TA^{\dagger}B$ yields the desired decomposition.
\end{proof}

{\lemma \label{intertwining_spectra} Assume the conditions of Definition \ref{def_intertwining} and let $L_i \subset {\cal S}^n$ be the linear hull of the image $\tilde f_i[K_i]$, $i = 1,2$. Then $K = L \cap {\cal S}_+^n$, where $L = L_1 + L_2$. Moreover, we have $\tilde f_i[K_i] = \Lambda_i \cap K$, where $\Lambda_i = Im\,\tilde f_i = {\cal L}_n(Im\,f_i)$, and the cones $\tilde f_i[K_i]$ are faces of $K$, $i = 1,2$. }

\begin{proof}
Set $H = Im\,f_1 \cap Im\,f_2$. Then $\Lambda_1 \cap \Lambda_2 = {\cal L}_n(H) = L_1 \cap L_2$. By definition we have $L_i \subset \Lambda_i$ and $K = (L_1 \cap {\cal S}_+^n) + (L_2 \cap {\cal S}_+^n)$.

Introduce a direct sum decomposition $\mathbb R^n = H_1' \oplus H \oplus H_2'$ such that $Im\,f_1 = H_1' \oplus H$, $Im\,f_2 = H \oplus H_2'$. Adopt a coordinate system in $\mathbb R^n$ which is adapted to this decomposition and partition the matrices in ${\cal S}^n$ accordingly. Then every matrix in $\Lambda_1 + \Lambda_2$, and hence also in $L$, has the form $X = \begin{pmatrix} X_{11} & X_{12} & 0 \\ X_{12}^T & X_{22} & X_{23} \\ 0 & X_{23}^T & X_{33} \end{pmatrix}$. Moreover, every matrix whose only nonzero block is $X_{22}$ is in $\Lambda_1 \cap \Lambda_2$ and hence in $L$.

Clearly $K = (L_1 \cap {\cal S}_+^n) + (L_2 \cap {\cal S}_+^n) \subset (L_1 + L_2) \cap {\cal S}_+^n = L \cap {\cal S}_+^n$. Let us show the reverse inclusion.

Let $X \in L \cap {\cal S}_+^n$ be an arbitrary matrix, partitioned as above. By Lemma \ref{3partition} there exists a decomposition $X_{22} = X_{22,1} + X_{22,2}$ such that the matrices
\[ X_1 = \begin{pmatrix} X_{11} & X_{12} & 0 \\ X_{12}^T & X_{22,1} & 0 \\ 0 & 0 & 0 \end{pmatrix} \in \Lambda_1,\quad X_2 = \begin{pmatrix} 0 & 0 & 0 \\ 0 & X_{22,2} & X_{23} \\ 0 & X_{23}^T & X_{33} \end{pmatrix} \in \Lambda_2
\]
are positive semi-definite. On the other hand, by virtue of $X \in L_1 + L_2$ there exists a decomposition $X = X_3 + X_4$ such that
\[ X_3 = \begin{pmatrix} X_{11} & X_{12} & 0 \\ X_{12}^T & X_{22,3} & 0 \\ 0 & 0 & 0 \end{pmatrix} \in L_1,\quad X_4 = \begin{pmatrix} 0 & 0 & 0 \\ 0 & X_{22,4} & X_{23} \\ 0 & X_{23}^T & X_{33} \end{pmatrix} \in L_2.
\]
We have $D_1 = X_1 - X_3 \in \Lambda_1 \cap \Lambda_2 \subset L_1$, $D_2 = X_2 - X_4 \in \Lambda_1 \cap \Lambda_2 \subset L_2$. Hence $X_1 = D_1 + X_3 \in L_1$, $X_2 = D_2 + X_4 \in L_2$. It follows that $X_1 \in \tilde f_1[K_1]$, $X_2 \in \tilde f_2[K_2]$. Therefore $X = X_1 + X_2 \in K$. Thus $L \cap {\cal S}_+^n \subset K$, which proves the first assertion of the lemma.

By construction we have $\tilde f_1[K_1] \subset \Lambda_1 \cap K$. Let us show the reverse inclusion.

Let $X \in \Lambda_1 \cap K$ be an arbitrary element. Since $K = \tilde f_1[K_1] + \tilde f_2[K_2]$, there exists a decomposition $X = X_1 + X_2$ such that $X_i \in \tilde f_i[K_i]$, $i = 1,2$. We have $X,X_1 \in \Lambda_1$, and hence $X_2 \in \Lambda_1 \cap \Lambda_2 \cap {\cal S}_+^n = {\cal F}_n(H) \subset \tilde f_1[K_1]$. Hence $X \in \tilde f_1[K_1]$, which proves $\Lambda_1 \cap K \subset \tilde f_1[K_1]$.

The equality $\Lambda_2 \cap K = \tilde f_2[K_2]$ is shown in a similar way.
\end{proof}

The isomorphism class of $K$ depends not only on the cones $K_i$ and the full faces $F_i \subset K_i$, but also on the maps $\iota_i$, or more precisely, on the linear bijection defined between $H_1$ and $H_2$ by the map $\iota_2 \circ \iota_1^{-1}$, where $\iota_1^{-1}$ is an arbitrary left inverse of $\iota_1$, because it is this bijection which determines the subspace $N$.

{\remark The intertwining operation is not "associative" in the sense that the intertwining $K_{123}$ of a cone $K_3$ with the intertwining $K_{12}$ of cones $K_1,K_2$ can always be represented as an intertwining of $K_2$ with an intertwining $K_{13}$ of $K_1,K_3$. However, since every full face of $K_{12}$ is a subset either of the face of $K_{12}$ isomorphic to $K_1$ or of the face of $K_{12}$ isomorphic to $K_2$, there exists a permutation $\sigma$ of the index set $\{1,2\}$ such that $K_{123}$ is an intertwining of $K_{\sigma(2)}$ with an intertwining $K_{\sigma(1)3}$ of the cones $K_{\sigma(1)},K_3$. }

We now come to the connection with ROG cones.

{\lemma \label{intertwining_ROG} Assume the conditions of Definition \ref{def_intertwining}. Then $K$ is a ROG cone if and only if $K_1,K_2$ are ROG cones. }

\begin{proof}
If $K$ is ROG, then $K_1,K_2$ are also ROG by Lemmas \ref{face_ROG} and \ref{intertwining_spectra}.

Assume that $K_1,K_2$ are ROG, then $\tilde f_1[K_1],\tilde f_2[K_2]$ are also ROG. Let $X \in K$ be arbitrary. Since $K = \tilde f_1[K_1] + \tilde f_2[K_2]$, there exist $X_i \in \tilde f_i[K_i]$, $i = 1,2$, such that $X = X_1 + X_2$. Since $\tilde f_i[K_i]$ are ROG, both $X_1$ and $X_2$ can be represented as a sum of rank 1 matrices in $\tilde f_1[K_1]$ and $\tilde f_2[K_2]$, respectively. Hence $X$ can be represented as a sum of rank 1 matrices in $\tilde f_1[K_1] \cup \tilde f_2[K_2] \subset K$. Thus $K$ is ROG.
\end{proof}

{\remark The preceding proof can in an obvious way be modified to show that if a spectrahedral cone is the sum of a finite number of its faces, then it is ROG if and only if all these faces are ROG. }

Finally we shall consider the special case $k = 1$. This case is simpler than the general case in two respects. Firstly, a face $F$ of a spectrahedral cone satisfying $\max_{X \in F} \rk X = 1$ is always full. For a ROG cone, the set of such faces equals the set of extreme rays. In particular, every ROG cone possesses full faces of rank 1. Therefore every two ROG cones can be intertwined along full faces of rank 1.

The second point is that given spectrahedral cones $K_1,K_2$ with full faces $F_i \subset K_i$ of rank 1, $i = 1,2$, the isomorphism class of the intertwining of $K_1$ and $K_2$ along the faces $F_1$ and $F_2$ is even independent of the maps $\iota_1,\iota_2$. This is because any two bijective maps $\iota,\iota': \mathbb R \to \mathbb R$ can be conjugated by a homothety of the target space, and this homothety can be compensated for by a homothety of one of the cones $K_1,K_2$. The cone constructed in Example \ref{ex_intertwining} is, e.g., the only cone which can be constructed by an intertwining of $\Han_+^3$ and ${\cal S}_+^2$ up to isomorphism, because every extreme ray of each of these cones can be taken to any other by an automorphism of the corresponding cone.

\section{Examples of ROG cones} \label{sec:examples}

In this section we consider two nontrivial families of ROG cones. We show that the class of ROG cones defined by chordal graphs can be constructed from the full matrix cones ${\cal S}_+^k$ by applying the constructive procedures presented in the previous section. We also provide an example of a continuous family of isomorphism classes of ROG cones.

\subsection{Cones defined by chordal graphs} \label{subs:chordal}

In this subsection we consider spectrahedral cones $K_G = L_G \cap {\cal S}_+^n$ defined by linear subspaces of the form $L_G = \{ X \in {\cal S}^n \,|\, X_{ij} = 0\ \forall\ (i,j) \not\in E(G) \}$, where $E(G)$ is the edge set of a graph $G$ on the vertices $1,\dots,n$. Note that the identity matrix is an element of $K_G$. Hence $K_G$ has a nonempty intersection with the interior of ${\cal S}_+^n$, and the linear span of $K_G$ equals $L_G$.

{\lemma \cite[Theorem 2.3]{AHMR88}, \cite[Theorem 2.4]{PaulsenPowerSmith87} \label{chordal} Assume above notations. Then the cone $K_G$ is ROG if and only if the graph $G$ is chordal. }

Chordal graphs are characterized by the condition that they admit a {\it perfect elimination ordering} of the vertices $1,\dots,n$. This is an ordering such that for every $k = 1,\dots,n$, the subset $N_k = \{ l < k \,|\, (l,k) \in E(G) \} \cup \{ k \}$ of vertices forms a {\it clique}, i.e., the subgraph of $G$ defined by $N_k$ is complete.

{\lemma Let $G$ be a chordal graph with vertex set $\{1,\dots,n\}$, and let $K_G$ be the corresponding ROG cone. Then $K_G$ can be constructed out of full matrix cones by iterated intertwinings or taking direct sums. }

\begin{proof}
Assume that the vertices are arranged in a perfect elimination ordering. For a subset $I \subset \{1,\dots,n\}$ of indices, define the linear subspace $H_I = \{x \in \mathbb R^n \,|\, x_i = 0\ \forall\ i \not\in I \}$. For $k = 1,\dots,n$, set $K_k = K_G \cap {\cal F}_n(H_{\{1,\dots,k\}})$.

Note that $K_1$ is isomorphic to the full matrix cone ${\cal S}_+^1$. We shall now show for all $k = 2,\dots,n$ that the cone $K_k$ is either an intertwining of $K_{k-1}$ with a full matrix cone, or a direct sum $K_{k-1} \oplus {\cal S}_+^1$.

Since $G$ is chordal, the set $N_k = \{ l < k \,|\, (l,k) \in E(G) \} \cup \{ k \}$ and its subset $N_k' = \{ l < k \,|\, (l,k) \in E(G) \}$ define cliques of $G$. Therefore the faces ${\cal F}_n(H_{N_k}),{\cal F}_n(H_{N_k'})$ of ${\cal S}_+^n$ are contained in $K$ and are full faces of this cone. In particular, ${\cal F}_n(H_{N_k'})$ is a full face of both ${\cal F}_n(H_{N_k})$ and $K_{k-1}$. On the other hand, $K_k = K_{k-1} + {\cal F}_n(H_{N_k})$ by definition of $N_k$. Hence $K_k$ is an intertwining of $K_{k-1}$ with the full matrix cone ${\cal F}_n(H_{N_k})$ in case that $N_k' \not= \emptyset$, and a direct sum $K_{k-1} \oplus {\cal F}_n(H_{\{k\}})$ in case that $N_k' = \emptyset$.

The proof is completed by the observation that $K_G = K_n$.
\end{proof}

{\lemma \label{chordal_simple} Let $G$ be a chordal graph with vertex set $\{1,\dots,n\}$, and let $K_G$ be the corresponding ROG cone. Then $\deg K_G = n$, and $K_G$ is simple if and only if $G$ is connected. }

\begin{proof}
By construction $K_G \subset {\cal S}_+^n$ contains the identity matrix, and hence $\deg K_G = n$ by Corollary \ref{deg_K}.

Suppose that $K_G$ is not simple. Then there exists a nontrivial direct sum decomposition $\mathbb R^n = H \oplus H'$ such that for every rank 1 matrix $xx^T \in K_G$, either $x \in H$ or $x \in H'$. In particular, if $x = e_i$ is a canonical basis vector, then $e_ie_i^T \in K_G$ by construction of $K_G$ and hence $e_i \in H \cup H'$ for all $i = 1,\dots,n$. Define the index sets $I = \{ i \,|\, e_i \in H \}$ and $I'= \{ i \,|\, e_i \in H' \}$. Then $I \cap I' = \emptyset$ and $I \cup I' = \{ 1,\dots,n \}$, because $\mathbb R^n = H \oplus H'$ is a direct sum decomposition. It follows that $H = \spa\{ e_i \,|\, i \in I \}$ and $H' = \spa\{ e_i \,|\, i \in I' \}$. Let now $x \in \mathbb R^n$ be a nonzero vector such that $X = xx^T \in K_G$. Then for every index pair $(i,j) \in I \times I'$ we have $x_ix_j = 0$ and hence $X_{ij} = 0$. From the fact that $K_G$ is a ROG cone it follows that $X_{ij} = 0$ for all $X \in \spa K_G = L_G$ in general for $(i,j) \in I \times I'$. But then $(i,j) \not\in E(G)$, and there is no edge in $G$ which connects the vertex subsets $I,I'$. Hence $G$ is not connected.

Suppose, on the other hand, that $G$ is not connected. Let $I,I'$ be disjoint nonempty vertex sets such that $I \cup I' = \{ 1,\dots,n \}$ and there is no edge in $G$ which connects $I$ to $I'$. Then by definition for every $X \in L_G$ we have $X_{ij} = x_ix_j = 0$ for every index pair $(i,j) \in I \times I'$. Define subspaces $H = \spa\{ e_i \,|\, i \in I \}$, $H' = \spa\{ e_i \,|\, i \in I' \}$ of $\mathbb R^n$. Then $\mathbb R^n = H \oplus H'$ is by construction a nontrivial direct sum decomposition. It then follows that $L_G \subset {\cal L}_n(H) + {\cal L}_n(H')$, and the cone $K_G$ is not simple.
\end{proof}

\subsection{A continuous family of non-isomorphic cones} \label{subs:continuous_family}

In this subsection we construct a family of ROG cones in ${\cal S}^6$ whose isomorphism class depends on a real parameter. First we shall explain the construction informally. We begin with a full matrix cone ${\cal S}_+^2$ and intertwine consecutively four other copies of ${\cal S}_+^2$ with it along faces of rank 1. This singles out a quadruple of rank 1 faces in the original copy of ${\cal S}_+^2$, or equivalently, a quadruple of points on the projective line $\mathbb RP^1$. Now two composite cones of this form are isomorphic if and only if the corresponding quadruples of faces can be taken to each other by an automorphism of ${\cal S}_+^2$, or equivalently, if the quadruples of points on $\mathbb RP^1$ are projectively equivalent. However, quadruples of points in $\mathbb RP^1$ possess a real invariant, the cross-ratio, which then parameterizes the isomorphism class of the composite cones.

We now define the cones formally. Fix mutually distinct angles $\varphi_1,\dots,\varphi_4 \in [0,\pi)$. For $\varphi \in [0,\pi)$, let $l(\varphi) \subset \mathbb R^2$ be the line through the origin with incidence angle $\varphi$. Then the lines $l(\varphi_1),\dots,l(\varphi_4)$ define a quadruple of points in real projective space $\mathbb RP^1$.

Consider the 11-dimensional subspace $L_{\varphi_1,\varphi_2,\varphi_3,\varphi_4} \subset {\cal S}^6$ of matrices of the form
\begin{equation} \label{cont_cones}
\begin{pmatrix} \alpha_1 & \alpha_2 & \alpha_3\cos\varphi_1 & \alpha_4\cos\varphi_2 & \alpha_5\cos\varphi_3 & \alpha_6\cos\varphi_4 \\
\alpha_2 & \alpha_7 & \alpha_3\sin\varphi_1 & \alpha_4\sin\varphi_2 & \alpha_5\sin\varphi_3 & \alpha_6\sin\varphi_4 \\
\alpha_3\cos\varphi_1 & \alpha_3\sin\varphi_1 & \alpha_8 & 0 & 0 & 0 \\
\alpha_4\cos\varphi_2 & \alpha_4\sin\varphi_2 & 0 & \alpha_9 & 0 & 0 \\
\alpha_5\cos\varphi_3 & \alpha_5\sin\varphi_3 & 0 & 0 & \alpha_{10} & 0 \\
\alpha_6\cos\varphi_4 & \alpha_6\sin\varphi_4 & 0 & 0 & 0 & \alpha_{11} \end{pmatrix},\qquad \alpha_1,\dots,\alpha_{11} \in \mathbb R.
\end{equation}
In the main result of this subsection, Lemma \ref{lem_cont_family} below, we shall show that the 4-dimensional family of spectrahedral cones $K_{\varphi_1,\varphi_2,\varphi_3,\varphi_4} = L_{\varphi_1,\varphi_2,\varphi_3,\varphi_4} \,\cap\, {\cal S}_+^6$ is ROG and under the isomorphism equivalence relation projects to a 1-dimensional family of isomorphism classes. However, first we construct a sequence of intermediate cones, with $K_{\varphi_1,\varphi_2,\varphi_3,\varphi_4}$ as the last element, and show that each one is obtained by an intertwining of the preceding one with the full matrix cone ${\cal S}_+^2$.

Let $H_0,\dots,H_4 \subset \mathbb R^6$ be the two-dimensional subspaces spanned by the columns of the matrices
\begin{equation} \label{H_matrices}
\begin{pmatrix} 1 & 0 \\ 0 & 1 \\ 0 & 0 \\ 0 & 0 \\ 0 & 0 \\ 0 & 0 \end{pmatrix},\ \begin{pmatrix} \cos\varphi_1 & 0 \\ \sin\varphi_1 & 0 \\ 0 & 1 \\ 0 & 0 \\ 0 & 0 \\ 0 & 0 \end{pmatrix},\ \begin{pmatrix} \cos\varphi_2 & 0 \\ \sin\varphi_2 & 0 \\ 0 & 0 \\ 0 & 1 \\ 0 & 0 \\ 0 & 0 \end{pmatrix},\ \begin{pmatrix} \cos\varphi_3 & 0 \\ \sin\varphi_3 & 0 \\ 0 & 0 \\ 0 & 0 \\ 0 & 1 \\ 0 & 0 \end{pmatrix},\ \begin{pmatrix} \cos\varphi_4 & 0 \\ \sin\varphi_4 & 0 \\ 0 & 0 \\ 0 & 0 \\ 0 & 0 \\ 0 & 1 \end{pmatrix},
\end{equation}
respectively. Define subspaces $L_j = \sum_{i=0}^j {\cal L}_6(H_i) \subset {\cal S}^6$ and spectrahedral cones $K_j = L_j \cap {\cal S}_+^6$, $j = 0,\dots,4$. Then $L_{\varphi_1,\varphi_2,\varphi_3,\varphi_4} = L_4$, $K_{\varphi_1,\varphi_2,\varphi_3,\varphi_4} = K_4$.

{\lemma \label{K_j6} The cone $K_j$, $j = 1,\dots,4$, is a ROG cone given by the sum $\sum_{i=0}^j {\cal F}_6(H_i)$. }

\begin{proof}
We prove the lemma by induction over $j$. For $j = 0$ the assertion holds by construction.

Let now $j > 0$ and assume that $K_{j-1} = \sum_{i=0}^{j-1} {\cal F}_6(H_i)$ is a ROG cone. We have to show that $K_j = K_{j-1} + {\cal F}_6(H_j)$ is ROG. To this end we construct $K_j$ as an appropriate intertwining of $K_{j-1}$ and ${\cal S}_+^2 \simeq {\cal F}_6(H_j)$.

Note that the non-zero elements of the matrices in $K_j$ are located in the upper left $(j+1) \times (j+1)$ subblock. For the sake of simplicity, we shall consider $K_{j-1}$ as a subset of ${\cal S}_+^j$ and $K_j$ as a subset of ${\cal S}_+^{j+1}$. Set $k = 1$, $n_1 = j$, $n_2 = 2$, $n = j+1$. Define the injections $\iota_1: \mathbb R \to \mathbb R^j$, $\iota_2: \mathbb R \to \mathbb R^2$, $f_1: \mathbb R^j \to \mathbb R^{j+1}$, $f_2: \mathbb R^2 \to \mathbb R^{j+1}$ by $\iota_1(1) = (\cos\varphi_j,\sin\varphi_j,0,\dots,0)^T$, $\iota_2(1) = (1,0)^T$, $f_1(x_1,\dots,x_j) = (x_1,\dots,x_j,0)^T$, $f_2(x_1,x_2) = (x_1\cos\varphi_j,x_1\sin\varphi_j,0,\dots,0,x_2)$. By construction, the subspace $L_j$ is the linear hull of the intertwining of $K_{j-1}$ and ${\cal S}_+^2$ defined by these maps, and ${\cal F}_6(H_j)$ is the image of ${\cal S}_+^2$ under the induced map $\tilde f_2$. Hence $K_j = K_{j-1} + {\cal F}_6(H_j)$ by Lemma \ref{intertwining_spectra} and $K_j$ is ROG by Lemma \ref{intertwining_ROG}. This completes the proof.
\end{proof}

{\lemma \label{lem_cont_family} The cone $K_{\varphi_1,\varphi_2,\varphi_3,\varphi_4}$ is a ROG cone. Two cones $K_{\varphi_1,\varphi_2,\varphi_3,\varphi_4}$, $K_{\varphi_1',\varphi_2',\varphi_3',\varphi_4'}$ are isomorphic if and only if the corresponding quadruples of lines $l(\varphi_1),\dots,l(\varphi_4) \subset \mathbb R^2$ and $l(\varphi_1'),\dots,l(\varphi_4') \subset \mathbb R^2$ define projectively equivalent quadruples of points in $\mathbb RP^1$. }

\begin{proof}
The first assertion of the lemma follows from Lemma \ref{K_j6} for $j = 4$.

Let us prove the second one. Consider cones $K_{\varphi_1,\varphi_2,\varphi_3,\varphi_4},K_{\varphi_1',\varphi_2',\varphi_3',\varphi_4'}$ for quadruples $(\varphi_1,\dots,\varphi_4)$, $(\varphi_1',\dots,\varphi_4')$ of mutually distinct angles. Let $H_0,\dots,H_4$ and $H_0',\dots,H_4'$, respectively, be the corresponding 2-dimensional subspaces of $\mathbb R^6$ as defined by the column spaces of the matrices \eqref{H_matrices}. Note that $H_0 = H_0'$. By Lemma \ref{K_j6} the set $\{ x \in \mathbb R^6 \,|\, xx^T \in K_{\varphi_1,\varphi_2,\varphi_3,\varphi_4} \}$ is given by the union $\bigcup_{j=0}^4 H_j$, and the set $\{ x \in \mathbb R^6 \,|\, xx^T \in K_{\varphi_1',\varphi_2',\varphi_3',\varphi_4'} \}$ by the union $\bigcup_{j=0}^4 H_j'$. The cones $K_{\varphi_1,\varphi_2,\varphi_3,\varphi_4}$ and $K_{\varphi_1',\varphi_2',\varphi_3',\varphi_4'}$ are then isomorphic if and only if there exists an invertible linear map $f: \mathbb R^6 \to \mathbb R^6$ which takes $\bigcup_{j=0}^4 H_j$ to $\bigcup_{j=0}^4 H_j'$.

Suppose that such a map $f$ exists. The intersections $l(\varphi_i) = H_0 \cap H_i$, $l(\varphi_i') = H_0' \cap H_i'$, $i = 1,2,3,4$, are 1-dimensional, while the intersections $H_i \cap H_j$, $H_i' \cap H_j'$, $i \not= j$, $i,j = 1,\dots,4$, are 0-dimensional. Hence we must have $f[H_0] = H_0'$ and $f[H_i] = H_{\sigma(i)}'$, $i = 1,\dots,4$, where $\sigma \in S_4$ is a permutation of the index set $\{1,\dots,4\}$. Moreover, $f|_{H_0}[l_i] = l_{\sigma(i)}'$, $i = 1,\dots,4$. It follows that $l(\varphi_1),\dots,l(\varphi_4) \subset H_0$ and $l(\varphi_1'),\dots,l(\varphi_4') \subset H_0$ define projectively equivalent quadruples of points in the projectivization of $H_0$.

Suppose now that the lines $l(\varphi_1),\dots,l(\varphi_4) \subset H_0$ and $l(\varphi_1'),\dots,l(\varphi_4') \subset H_0$ define projectively equivalent quadruples of points in the projectivization of $H_0$. Then there exists an invertible linear map $h: H_0 \to H_0$ and a permutation $\sigma \in S_4$ such that $h[l_i] = l_{\sigma(i)}'$, $i = 1,\dots,4$. Let now $x_i \in H_i \setminus l_i$, $x_i' \in H_i' \setminus l_i'$, $i = 1,\dots,4$, be arbitrary points. We then have $H_i = \spa(l_i \cup \{x_i\})$, $H_i' = \spa(l_i' \cup \{x_i'\})$, $i = 1,\dots,4$. Moreover, $\spa(H_0 \cup \{x_1,x_2,x_3,x_4\}) = \spa(H_0 \cup \{x_1',x_2',x_3',x_4'\}) = \mathbb R^6$. We then can extend the map $h$ to a linear map $f: \mathbb R^6 \to \mathbb R^6$ such that $f(x_i) = x_{\sigma(i)}'$, $i = 1,\dots,4$. This map is invertible by construction and $f[H_i] = H_{\sigma(i)}'$, $i = 1,\dots,4$. It follows that $f[\bigcup_{j=0}^4 H_j] = \bigcup_{j=0}^4 H_j'$, which completes the proof.
\end{proof}

It is well-known that there exist infinitely many projectively non-equivalent quadruples of points in $\mathbb RP^1$. The equivalence classes are parameterized by the orbits of the {\it cross-ratio} $\lambda(\varphi_1,\varphi_2,\varphi_3,\varphi_4) = (l_1,l_2;l_3,l_4) = \frac{(\cot\varphi_1-\cot\varphi_3)(\cot\varphi_2-\cot\varphi_4)}{(\cot\varphi_2-\cot\varphi_3)(\cot\varphi_1-\cot\varphi_4)}$ with respect to the action of the symmetric group $S_4$ on the arguments $\varphi_1,\dots,\varphi_4$. Thus there exists a continuum of mutually non-isomorphic ROG cones defined by subspaces $L \subset {\cal S}^6$ of type \eqref{cont_cones}.

The cone $K_{\varphi_1,\varphi_2,\varphi_3,\varphi_4}$ is obtained from the face ${\cal F}_6(H_0) \cong {\cal S}_+^2$ by consecutive intertwining with the faces ${\cal F}_6(H_i) \cong {\cal S}_+^2$, $i = 1,\dots,4$. It is hence an intertwining of 5 full matrix cones ${\cal S}_+^2$. More complicated ROG cones can be obtained by starting with a matrix cone ${\cal S}_+^n$ and consecutively intertwining it with matrix cones ${\cal S}_+^{k_1},\dots,{\cal S}_+^{k_m}$ along full faces of ranks $d_1,\dots,d_m$, where $1 \leq d_i < \min(n,k_i)$, $i = 1,\dots,m$. In this way, families of mutually non-isomorphic ROG cones can be obtained which are parameterized by an arbitrary number of real parameters. Note that all cones obtained in such a way are simple.

\section{Dimension and degree of ROG cones} \label{sec:dim}

In this section we consider the relation between the dimension and the degree of a ROG cone $K$. Evidently we have the inequality chain $\deg K \leq \dim K \leq \frac{\deg K(\deg K + 1)}{2}$, with equality on the left if and only if $K$ is isomorphic to the cone of positive semi-definite diagonal matrices, and equality on the right if and only if $K$ is isomorphic to the full cone of positive semi-definite matrices. We shall say that a ROG cone $K$ has {\it codimension $k$} if $\dim K = \frac{\deg K(\deg K + 1)}{2} - k$. The codimension can be interpreted as the number of linearly independent linear constraints on the matrices $X \in K$ in any non-degenerate representation of $K$.

{\lemma \label{lem:codim_k} Let $K$ be a ROG cone of degree $n$ and dimension $\frac{n(n + 1)}{2} - k$. Then $K$ has a representation $K = \{ X \in {\cal S}_+^n \,|\, \langle X,Q_i \rangle = 0\ \forall\ i = 1,\dots,k \}$, where $Q_1,\dots,Q_k$ are linearly independent quadratic forms on $\mathbb R^n$ such that every nonzero form in the linear span of $\{Q_1,\dots,Q_k\}$ is indefinite. }

\begin{proof}
By Lemma \ref{HVlemma2} and Corollary \ref{deg_K} there exists a non-degenerate representation of $K$ as a linear section of ${\cal S}_+^n$. We have $\dim {\cal S}^n - \dim K = k$, and the orthogonal complement of $\spa K$ in the space of quadratic forms on $\mathbb R^n$ has dimension $k$. Let $\{Q_1,\dots,Q_k\}$ be a basis of this complement. Then by construction we have $K = \spa K \cap {\cal S}_+^n = \{ X \in {\cal S}_+^n \,|\, \langle X,Q_i \rangle = 0\ \forall\ i = 1,\dots,k \}$.

Let $X \in K$ be a positive definite matrix. Suppose for the sake of contradiction that there exists a nonzero linear combination $Q$ of $Q_1,\dots,Q_k$ which is semi-definite. By possibly replacing $Q$ by $-Q$, we may assume that $Q$ is positive semi-definite. Then $\langle Q,X \rangle > 0$, leading to a contradiction.
\end{proof}

In the next subsections we classify ROG cones of codimensions 1 and 2, and give a lower bound on the dimension of simple cones $K$ of fixed degree.

\subsection{ROG cones of codimension 1} \label{subs:codim1}

In this subsection we show that all spectrahedral cones of codimension 1 are ROG. This result is closely linked to Dines' and Brickmans theorems on the convexity of the numerical range of quadratic forms \cite{Dines43},\cite{Brickman61}. All these results are based on the following dimensional argument.

{\lemma \label{lem_codim} Let $L \subset {\cal S}^n$ be a linear subspace of dimension $\frac{n(n + 1)}{2} - d$. Then the spectrahedral cone $K = L \cap {\cal S}_+^n$ has no extreme elements of rank $k > -\frac12+\sqrt{\frac14+2(d+1)}$. }

\begin{proof}
Let $X$ lie on an extreme ray of $K$, and let $k = \rk X$. Then the minimal face of ${\cal S}_+^n$ which contains $X$ has dimension $\frac{k(k+1)}{2}$. Denote this face by $F$. The minimal face of $K$ which contains $X$ is given by the intersection $F \cap L$ and has dimension 1. But since $L$ has codimension $d$, we have $1 = \dim(F \cap L) \geq \dim F - d = \frac{k(k+1)}{2} - d$. This yields $k(k+1) \leq 2(d+1)$, which implies $k \leq -\frac12+\sqrt{\frac14+2(d+1)}$.
\end{proof}

{\corollary \label{cor_codim1} Let $L \subset {\cal S}^n$ be a linear subspace of dimension $\frac{n(n + 1)}{2} - 1$. Then the cone $K = L \cap {\cal S}_+^n$ is ROG. }

\begin{proof}
By Lemma \ref{lem_codim} the cone $K$ has no extreme elements of rank $k \geq 2 > \frac{-1+\sqrt{17}}{2}$. Thus $K$ is ROG.
\end{proof}

{\corollary \label{class_codim1} Every ROG cone of degree $n$ and codimension 1 has a representation of the form $K = \{ X \in {\cal S}_+^n \,|\, \langle X,Q \rangle = 0 \}$ for some indefinite quadratic form $Q$, and every cone of this form is ROG of degree $n$ and codimension 1. Two such cones $K,K'$, defined by indefinite quadratic forms $Q,Q'$, respectively, are isomorphic if and only if either $Q,Q'$ or $Q,-Q'$ have the same signature. }

\begin{proof}
The first claim follows from Lemma \ref{lem:codim_k}.

Let now $Q$ be an indefinite quadratic form. Then the cone $K = \{ X \in {\cal S}_+^n \,|\, \langle X,Q \rangle = 0 \}$ is ROG by Corollary \ref{cor_codim1}. Since $Q$ is indefinite, there exists a positive definite matrix $X$ such that $\langle X,Q \rangle = 0$. Hence $K$ intersects the interior of ${\cal S}_+^n$, and therefore  $\dim K = \dim {\cal S}^n - 1$. Moreover, by Corollary \ref{deg_K} $K$ is of degree $n$.

Let now the cones $K,K'$ be defined by indefinite quadratic forms $Q,Q'$, respectively. The cones $K,K'$ are isomorphic if and only if their linear hulls $L = \{ X \in {\cal S}^n \,|\, \langle X,Q \rangle = 0 \}$, $L' = \{ X \in {\cal S}^n \,|\, \langle X,Q' \rangle = 0 \}$ can be taken to each other by a coordinate transformation of $\mathbb R^n$. This is the case if and only if the orthogonal complements of $L,L'$, namely the 1-dimensional subspaces generated by $Q$ and $Q'$, are related by a coordinate transformation. The last claim now easily follows.
\end{proof}

It is not hard to establish that there are $[\frac{n^2}{4}]$ isomorphism classes of ROG cones of degree $n$ and codimension 1. For $n \geq 3$ all of them are simple.

\subsection{ROG cones of codimension 2} \label{subs:codim2}

In this subsection we classify the ROG cones $K$ of degree $n$ and dimension $\frac{n(n+1)}{2}-2$. If $n = 2$, then the dimension of $K$ is either 2 or 3, and $K$ cannot be of codimension 2. We shall henceforth assume $n \geq 3$. For the classification we shall need the auxiliary Lemmas \ref{synthesis_rk2} and \ref{pencils} which are provided in the Appendix.

{\theorem \label{class_codim2} Let $K$ be a ROG cone of degree $n \geq 3$ and of codimension $d = 2$. Then $K$ is isomorphic to the direct sum ${\cal S}_+^1 \oplus {\cal S}_+^2$ if $n = 3$ and to a full extension of this sum if $n > 3$. }

\begin{proof}
By Lemma \ref{HVlemma2} and Corollary \ref{deg_K} we may assume that $K$ has a non-degenerate representation by matrices of size $n \times n$. By Lemma \ref{lem:codim_k} we have $K = \{ X \in {\cal S}_+^n \,|\, \langle X,Q_1 \rangle = \langle X,Q_2 \rangle = 0 \}$ for some linearly independent quadratic forms $Q_1,Q_2$ on $\mathbb R^n$. Since $K$ is ROG, it has no extremal elements of rank 2 and hence the assumptions of Lemma \ref{synthesis_rk2} in the Appendix are fulfilled.

Suppose that condition (i) of this lemma holds. Then for every $z \in \mathbb R^n$ such that $zz^T$ is an extreme element of $K$, the linear forms $Q_1z,Q_2z$ are linearly dependent. This implies that $z$ is an eigenvector of the pencil $Q_1 + \lambda Q_2$.

Since the degree of $K$ is $n$, by Corollary \ref{cor_diag2} there exist $n$ linearly independent vectors $z_1,\dots,z_n \in \mathbb R^n$ such that the rank 1 matrices $z_kz_k^T$ are in $K$ for $k = 1,\dots,n$. This implies that the pencil $Q_1 + \lambda Q_2$ has $n$ linearly independent real eigenvectors. Therefore the conditions of Lemma \ref{pencils} are satisfied. Let $\mathbb R^n = H_0 \oplus H_1 \oplus \dots \oplus H_m$ be the direct sum decomposition from this lemma. If $m \leq 1$, then the forms $Q_1,Q_2$ are linearly dependent, which contradicts our assumptions. Hence $m \geq 2$.

Let $x_1 \in H_1$, $x_2 \in H_2$ be nonzero vectors. Consider the matrix $X = x_1x_2^T + x_2x_1^T \in {\cal S}^n$. We have $\langle Q_i,X \rangle = 2x_1^TQ_ix_2 = 0$ for $i = 1,2$, and hence $X \in \spa K$. On the other hand, $\spa K$ is generated by all rank 1 matrices in $K$ because $K$ is ROG. However, if $z \in \mathbb R^n$ is such that $zz^T \in K$, then by Lemma \ref{pencils} we have $z \in \bigcup_{k=1}^m (H_0+H_k)$. It follows that $\spa K \subset \sum_{k=1}^m {\cal L}_n(H_0+H_k)$. But $X \not\in \sum_{k=1}^m {\cal L}_n(H_0+H_k)$, leading to a contradiction.

Thus condition (ii) of Lemma \ref{synthesis_rk2} holds. By choosing an appropriate basis of $\mathbb R^n$, we can assume that the linear forms $u,q_1,q_2$ from this lemma are the first elements of the dual basis. Then the cone $K$ is given by the set $\{ X \in {\cal S}_+^n \,|\, X_{12} = X_{13} = 0 \}$. The claim of the theorem now easily follows.
\end{proof}

Hence there is only one isomorphism class of ROG cones of codimension 2 for a given degree $n \geq 3$, in contrast to ROG cones of codimension 1, of which there are many.

\subsection{Lower bound on the dimension of simple ROG cones} \label{subs:dim_bound}

In this section we show that for simple ROG cones $K$ the dimension of $K$ is bounded from below by $2\cdot \deg K-1$. This will be useful later for the classification of simple ROG cones of low degree. We shall need the following auxiliary result.

{\lemma \label{aux_lower} Let $x_1,\dots,x_m \in \mathbb R^n$ be linearly independent vectors, and let $S \subset {\cal S}^n$ be the $m$-dimensional subspace spanned by the rank 1 matrices $x_1x_1^T,\dots,x_mx_m^T$. Let further $H \subset \mathbb R^n$ be a linear subspace. Then the dimension of the intersection $S \cap {\cal L}_n(H)$ is given by the number of indices $i$ such that $x_i \in H$. In particular, $\dim (S \cap {\cal L}_n(H)) \leq \dim H$. }

\begin{proof}
Define the index set $I = \{ i \,|\, x_i \in H \}$. Let $A = \sum_{i=1}^m \alpha_i x_ix_i^T$ be an arbitrary element of $S$, where $\alpha_i$ are scalar coefficients. Suppose there exists an index $j \not\in I$ such that $\alpha_j \not= 0$. Let $y \in \mathbb R^n$ be a vector such that $y^Tx_j = 1$, and $y^Tx_i = 0$ for all $i \not= j$. Such a vector $y$ exists by the linear independence of $x_1,\dots,x_m$. We then get $Ay = \sum_{i=1}^m \alpha_i (y^Tx_i)x_i = \alpha_jx_j \not= H$. Hence $A \not\in {\cal L}_n(H)$.

It follows that every matrix in the intersection $S \cap {\cal L}_n(H)$ is of the form $A = \sum_{i \in I} \alpha_i x_ix_i^T$ for some scalars $\alpha_i$. On the other hand, for every such matrix $A$ and every vector $y \in \mathbb R^n$ we have $Ay = \sum_{i \in I} (y^Tx_i)x_i \in H$, and $A \in {\cal L}_n(H)$. Therefore the intersection $S \cap {\cal L}_n(H)$ equals the linear span of the set $\{ x_ix_i^T \,|\, i \in I \}$. The claims of the lemma now easily follow.
\end{proof}

{\theorem \label{thm_lower} Let $K$ be a simple ROG cone of degree $n$. Then $\dim K \geq 2n-1$. }

\begin{proof}
Represent $K$ as a linear section of ${\cal S}_+^n$. Recall that by Lemma \ref{face_ROG} every face of $K$ is a ROG cone, and that $K$ itself is the face of $K$ of largest degree $n$. Denote by ${\bf F}$ the set of faces $F$ of $K$ such that $\dim F \geq 2\deg F - 1$. The set ${\bf F}$ is not empty, because every extreme ray of $K$ is an element of ${\bf F}$. Set $k = \max_{F \in {\bf F}} \deg F$. Assume for the sake of contradiction that $K \not\in {\bf F}$, and hence $k < n$. Let $F_k \in {\bf F}$ be a face of $K$ which achieves the maximal degree $k$. Denote the linear span of $K$ by $L$, and the linear span of $F_k$ by $L_k$. By construction we have $\dim L_k \geq 2k - 1$.

By Corollary \ref{deg_K} the maximal rank of matrices in $F_k$ equals $k$. Let $Y \in F_k$ be a matrix of maximal rank $k$, and let the $k$-dimensional subspace $H \subset \mathbb R^n$ be its image. Then we have $L_k = L \cap {\cal L}_n(H)$ and $F_k = L \cap {\cal F}_n(H)$. By Corollary \ref{cor_diag1} there exists a basis $\{ r_1,\dots,r_k \}$ of $H$ such that $r_ir_i^T \in K$ for all $i = 1,\dots,k$, and $Y = \sum_{i=1}^k r_ir_i^T$. By virtue of $\deg K = n$ and Corollary \ref{cor_diag2} we may complete this basis of $H$ to a basis $\{ r_1,\dots,r_n \}$ of $\mathbb R^n$ such that $r_ir_i^T \in K$ for all $i = 1,\dots,n$. Adopt the coordinate system defined by this basis. Then all diagonal matrices are in $L$, and the subspace $L_k$ consists of the matrices in $L$ all whose non-zero elements are located in the upper left $k \times k$ block.

Since $K$ is simple, there exists a rank 1 matrix $zz^T \in K$ such that the vector $z = (z_1,\dots,z_n)^T$ is neither in $H$ nor in $\spa\{r_{k+1},\dots,r_n\}$. In other words, the subvector $z_H = (z_1,\dots,z_k)^T$ is not zero, and not all of the elements $z_{k+1},\dots,z_n$ are zero. Without loss of generality, let the nonzero elements in the second group be $z_{k+1},\dots,z_{k+m}$. By scaling the vector $z$, we may also assume that $z^Tz = 1$.

Denote by $F_{k+m}$ the face of $K$ which consists of all matrices in $K$ whose non-zero elements are located in the upper left $(k+m) \times (k+m)$ block. Denote the linear span of $F_{k+m}$ by $L_{k+m}$. Since all diagonal matrices are in $L$, the maximal rank of the matrices in $F_{k+m}$ equals $k+m$. By Corollary \ref{deg_K} we get $\deg F_{k+m} = k+m > k$. By our definition of $k$ we then have $F_{k+m} \not\in {\bf F}$, and hence $\dim L_{k+m} < 2(k+m)-1$. Let $S$ be the $(\dim L_k + m)$-dimensional subspace of $L_{k+m}$ spanned by $L_k$ and the rank 1 matrices $r_{k+1}r_{k+1}^T,\dots,r_{k+m}r_{k+m}^T$.

We have $zz^T \in F_{k+m}$. Consider the matrix $X = \diag(I_{k+m},0,\dots,0) - zz^T \in L_{k+m}$, where $I_{k+m}$ is the $(k+m) \times (k+m)$ identity matrix. By $z^Tz = 1$ the matrix $X$ is positive semi-definite of rank $k+m-1$, with $z$ as kernel vector. It follows that $X \in F_{k+m}$, and by Corollary \ref{cor_diag1} there exist $k+m-1$ linearly independent vectors $x_1,\dots,x_{k+m-1} \in \mathbb R^n$ such that $x_ix_i^T \in F_{k+m}$ for all $i = 1,\dots,k+m-1$, and $X = \sum_{i=1}^{k+m-1} x_ix_i^T$. Since $z^TXz = \sum_{i=1}^{k+m-1} (z^Tx_i)^2 = 0$, it follows that $z^Tx_i = 0$ for all $i = 1,\dots,k+m-1$.

Consider the $(k+m-1)$-dimensional subspace $S' \subset L_{k+m}$ spanned by the rank 1 matrices $x_ix_i^T$, $i = 1,\dots,k+m-1$. Let us bound the dimension of the intersection $S \cap S'$. Let $A \in S \cap S'$ be arbitrary. Since $A \in S$, the matrix $A$ has a block-diagonal structure $A = \diag(A_H,a_{k+1},\dots,a_{k+m},0,\dots,0)$, with $A_H$ a block of size $k \times k$. On the other hand, $A \in S'$ implies $Az = 0$. It follows that $a_{k+1}z_{k+1} = \dots = a_{k+m}z_{k+m} = 0$ and $a_{k+1} = \dots = a_{k+m} = 0$, because the corresponding elements of $z$ are non-zero. The image of $A$ is hence contained in the intersection of the subspace $H$ with the orthogonal complement of $z$. By virtue of $z_H \not= 0$ this intersection has dimension $k-1$. By Lemma \ref{aux_lower} we then get that $\dim (S \cap S') \leq k-1$.

Thus $\dim (S + S') = \dim S + \dim S' - \dim (S \cap S') \geq (\dim L_k + m) + (k+m-1) - (k-1) \geq 2k - 1 + 2m$, leading to a contradiction with the bound $\dim L_{k+m} < 2(k+m)-1$. This completes the proof.
\end{proof}

\section{Isolated extreme rays} \label{sec:discrete}

The extreme rays of a ROG cone are generated by its rank 1 matrices. In this section we study the situation when an extreme ray of a ROG cone $K$ is isolated. We shall show that in this case $K$ is a direct sum of ${\cal S}_+^1$ and a lower-dimensional ROG cone, and the isolated extreme ray is the face of $K$ corresponding to the factor ${\cal S}_+^1$. We deduce a couple of results for simple ROG cones and consider the situation when a simple ROG cone $K$ has a face of codimension 2. We will need the following concept.

{\definition \label{mld_def} The vectors $x_1,\dots,x_{k+1} \in \mathbb R^n$ are called {\it minimally linearly dependent} if they are linearly dependent, but every $k$ of them are linearly independent. }

{\lemma \label{mld_char} A set of vectors $x_1,\dots,x_{k+1} \in \mathbb R^n$ is minimally linearly dependent if and only if their span has dimension $k$ and there exist nonzero real numbers $c_1,\dots,c_{k+1}$ such that $\sum_{i=1}^{k+1} c_ix_i = 0$. }

\begin{proof}
Denote by $L$ the linear span of $\{x_1,\dots,x_{k+1}\}$, and let $X$ be the $n \times (k+1)$ matrix formed of the column vectors $x_i$.

Let $x_1,\dots,x_{k+1} \in \mathbb R^n$ be minimally linearly dependent. Then the dimension of $L$ equals $k$, because there exist $k$ linearly independent vectors in $L$. The matrix $X$ then has rank $k$ and its kernel has dimension 1. Let $(c_1,\dots,c_{k+1})^T \in \mathbb R^{k+1}$ be a generator of $\ker X$. Then $\sum_{i=1}^{k+1} c_ix_i = 0$ and not all $c_i$ are zero. Let $I \subset \{1,\dots,k+1\}$ be the set of indices $i$ such that $c_i \not= 0$. Then the vectors in the set $\{ c_i \,|\, i \in I \}$ are linearly dependent. By assumption, no $k$ vectors are linearly dependent, and therefore $I$ has not less than $k+1$ elements. It follows that $c_i \not= 0$ for all $i$.

Let now $c_1,\dots,c_{k+1}$ be nonzero real numbers such that $\sum_{i=1}^{k+1} c_ix_i = 0$, and suppose $\dim L = k$. Then $x_1,\dots,x_{k+1}$ are linearly dependent. Moreover, $\rk X = k$, and hence the vector $(c_1,\dots,c_{k+1})^T$ generates the kernel of $X$. In particular, there is no nonzero kernel vector with a zero element. It follows that every subset of $k$ vectors is linearly independent. Thus $x_1,\dots,x_{k+1}$ are minimally linearly dependent.
\end{proof}

{\lemma \label{mld_alternative} Let $S \subset \mathbb R^n$ be a subset and $x \in S$ a nonzero vector. Then either

1) there exists a subspace $H \subset \mathbb R^n$ of dimension $n-1$ which does not contain $x$, such that for every $y \in S$ either $y \in H$ or $y$ is a multiple of $x$,

or 2) there exists a minimally linearly dependent subset $T \subset S$ of size at least 3 such that $x \in T$. }

\begin{proof}
Let $L \subset \mathbb R^n$ be the linear span of $S$, and let $k$ be its dimension. Let us complete $x_1 = x$ to a basis $\{x_1,\dots,x_k\} \subset S$ of $L$. Then every vector $y \in S$ can be in a unique way represented as a sum $y = \sum_{i=1}^k c_ix_i$. We have two possibilities.

1) For every vector $y = \sum_{i=1}^k c_ix_i \in S$, either $c_1 = 0$, or $c_2 = \dots = c_k = 0$. Then we can take $H$ as any hyperplane which contains the span of $\{x_2,\dots,x_k\}$ but not $x_1$, and are in the situation 1) of the lemma.

2) There exists $y = \sum_{i=1}^k c_ix_i \in S$ such that $c_1 \not= 0$ and at least one of the coefficients $c_2,\dots,c_k$ is not zero. Let without loss of generality the nonzero coefficients among the $c_2,\dots,c_k$ be the coefficients $c_2,\dots,c_l$, $l \geq 2$. Then we obtain $y - \sum_{i=1}^l c_ix_i = 0$, and the set $\{x_1,\dots,x_l,y\} \subset S$ is minimally linearly dependent by Lemma \ref{mld_char}. Thus we are in the situation 2) of the lemma.
\end{proof}

{\lemma \label{mld_welldefined} Let $K$ be a ROG cone and let $R_1,\dots,R_{k+1} \in K$ be extreme rays of $K$. Let the rank 1 matrices $X_i = x_ix_i^T$ be generators of these extreme rays, respectively, in some representation of $K$ as a linear section of a positive semi-definite matrix cone ${\cal S}_+^n$. Whether the set $\{x_1,\dots,x_{k+1}\} \subset \mathbb R^n$ is minimally linearly dependent then depends only on the extreme rays $R_1,\dots,R_{k+1}$ of $K$, but not on the representation of $K$, its size, or the generators $X_i$. }

\begin{proof}
Let $c_1,\dots,c_{k+1}$ be non-zero real numbers. Then a subset $\{x_1,\dots,x_{k+1}\} \subset \mathbb R^n$ is minimally linearly dependent if and only if the subset $\{c_1x_1,\dots,c_{k+1}x_{k+1}\}$ is minimally linearly dependent. This follows directly from Definition \ref{mld_def}. Hence the property does not depend on the generators $X_i$ of the extreme rays for a given representation of $K$. Let now $X_i = x_ix_i^T$, $Y_i = y_iy_i^T$ be generators of the rays $R_i$ in different representations of sizes $n,m$, respectively. Let $n \leq m$ without loss of generality. By Theorem \ref{theorem_iso} there exists an injective linear map $f: \mathbb R^n \to \mathbb R^m$ such that $f(x_i) = \sigma_iy_i$, where $\sigma_i \in \{-1,+1\}$, for all $i = 1,\dots,k+1$. By the injectivity of $f$, we have for every index subset $I \subset \{1,\dots,k+1\}$ that the set $\{x_i\}_{i \in I}$ is linearly dependent if and only if the set $\{\sigma_iy_i\}_{i \in I}$ is linearly dependent. Hence $\{x_1,\dots,x_{k+1}\}$ is minimally linearly dependent if and only if the set $\{y_1,\dots,y_{k+1}\}$ is. This completes the proof.
\end{proof}

Lemma \ref{mld_welldefined} allows to make the following definition.

{\definition Let $K$ be a ROG cone. We call a subset $\{R_1,\dots,R_{k+1}\}$ of extreme rays of $K$, generated by rank 1 matrices $X_i = x_ix_i^T$, respectively, an {\it MLD set}, if the set $\{x_1,\dots,x_{k+1}\}$ is minimally linearly dependent. }

{\lemma \label{dim_mld} Let $K$ be a ROG cone of degree $n \geq 2$, possessing an MLD set $\{R_1,\dots,R_{n+1}\}$ of extreme rays. Then the following holds:

i) the cone $K$ is simple;

ii) the extreme rays $R_1,\dots,R_{n+1}$ of $K$ are not isolated. }

\begin{proof}
Represent $K$ as a linear section of ${\cal S}_+^n$, and let the rank 1 matrices $X_i = x_ix_i^T$ be generators of the extreme rays $R_i$, $i = 1,\dots,n+1$. Then the set $\{x_1,\dots,x_{n+1}\} \subset \mathbb R^n$ is minimally linearly dependent. Denote the linear span of $K$ by $L$.

Suppose for the sake of contradiction that $K$ is not simple. Then there exists a nontrivial direct sum decomposition $\mathbb R^n = H_1 \oplus H_2$ such that $K \subset {\cal L}_n(H_1) + {\cal L}_n(H_2)$ and hence $x_i \in H_1 \cup H_2$ for all $i = 1,\dots,n+1$. Let $n_1,n_2$ be the dimensions of $H_1,H_2$, respectively, and $n_1',n_2'$ the number of indices $i$ such that $x_i \in H_1$ or $x_i \in H_2$, respectively. Then $n_1',n_2' > 0$, because the vectors $x_1,\dots,x_{n+1}$ span the whole space $\mathbb R^n$ and the decomposition $\mathbb R^n = H_1 \oplus H_2$ is nontrivial. On the other hand, we have $n_1 + n_2 = n$ and $n_1' + n_2' = n + 1$. Hence either $n_1' > n_1$, or $n_2' > n_2$, and there exists a strict subset of the set $\{x_1,\dots,x_{n+1}\}$ which is linearly dependent, leading to a contradiction. This proves i).

We shall now prove ii). For $n = 2$ we have $K = {\cal S}_+^2$, and the assertion is evident. Suppose $n \geq 3$.

By the definition of minimal linear dependence the vectors $x_1,\dots,x_n$ form a basis of $\mathbb R^n$. Choose a coordinate system in which this is the canonical basis. By Lemma \ref{mld_char} there exist nonzero scalars $c_1,\dots,c_{n+1}$ such that $\sum_{i=1}^{n+1}c_ix_i = 0$. We may normalize these scalars by a common factor to achieve $c_{n+1} = -1$. Then we have $x_{n+1} = (c_1,\dots,c_n)^T$.

The subspace $L \subset {\cal S}^n$ contains the $(n+1)$-dimensional linear span $\tilde L$ of the rank 1 matrices $x_ix_i^T$, $i = 1,\dots,n+1$. Let $d_1,\dots,d_n > 0$ be positive scalars, and set $d_{n+1} = -\left(\sum_{i=1}^n d_i^{-1}c_i^2\right)^{-1}$. Then the matrix $M = \sum_{i=1}^{n+1} d_ix_ix_i^T$ is an element of $\tilde L$. Moreover, for every vector $r = (r_1,\dots,r_n)^T$ we have
\[ r^TMr = \sum_{i=1}^{n+1} d_i(r^Tx_i)^2 = \sum_{i=1}^n d_ir_i^2 - \frac{\left( \sum_{i=1}^n c_ir_i \right)^2}{\sum_{i=1}^n d_i^{-1}c_i^2} = \sum_{i=1}^n \left(\sqrt{d_i}r_i - \frac{c_i \sum_{j=1}^n c_jr_j}{\sqrt{d_i}\sum_{j=1}^n d_j^{-1}c_j^2}\right)^2 \geq 0.
\]
It follows that $M \succeq 0$ and hence $M \in K$.

Moreover, we have $r^TMr = 0$ if and only if $\sqrt{d_i}r_i = \frac{c_i \sum_{j=1}^n c_jr_j}{\sqrt{d_i}\sum_{j=1}^n d_j^{-1}c_j^2}$ for all $i = 1,\dots,n$. An equivalent condition is that $r = \alpha s$ for some scalar $\alpha$, where $s = (s_1,\dots,s_n)^T$ is a vector given by $s_i = d_i^{-1}c_i$ for all $i = 1,\dots,n$. Hence $M$ is of rank $n-1$, in particular, it is not rank 1.

Let $H$ be the $(n-1)$-dimensional subspace of vectors $v \in \mathbb R^n$ such that $v^Ts = 0$. Then the minimal face of ${\cal S}_+^n$ which contains $M$ is given by ${\cal F}_n(H)$. It consists of all matrices $X \in {\cal S}_+^n$ such that $Xs = 0$. The linear span ${\cal L}_n(H)$ of this face is given by all $X \in {\cal S}^n$ such that $Xs = 0$. We shall now compute the intersection ${\cal L}_n(H) \cap \tilde L$. Let $X = \sum_{i=1}^{n+1} \alpha_i x_ix_i^T \in {\cal L}_n(H) \cap \tilde L$. Then we have
\[ Xs = \sum_{i=1}^{n+1} \alpha_i (x_i^Ts)x_i = \sum_{i=1}^n \alpha_is_ix_i + \alpha_{n+1} \cdot \sum_{j=1}^n c_js_j \cdot \sum_{i=1}^n c_ix_i = \sum_{i=1}^n \left( \alpha_id_i^{-1} + \alpha_{n+1} \sum_{j=1}^n d_j^{-1}c_j^2 \right)c_ix_i = 0.
\]
It follows that $\alpha_i d_{n+1} = \alpha_{n+1} d_i$ for all $i = 1,\dots,n$. An equivalent condition is that the vectors $\alpha = (\alpha_1,\dots,\alpha_{n+1})^T$ and $d = (d_1,\dots,d_{n+1})^T$ are proportional, and hence $X$ is proportional to $M$. It follows that ${\cal L}_n(H) \cap \tilde L$ is the 1-dimensional subspace generated by $M$.

By Lemma \ref{Cara_X} there exist $n-1$ linearly independent vectors $y_1,\dots,y_{n-1} \in \mathbb R^n$ such that $y_iy_i^T \in K$ for all $i$ and $M = \sum_{i=1}^{n-1} y_iy_i^T$. Note that $y_iy_i^T \in {\cal L}_n(H)$ for all $i$, and hence $y_iy_i^T \in {\cal L}_n(H) \cap L$.


Assume for the sake of contradiction that the extreme ray $R_1$ generated by the rank 1 matrix $x_1x_1^T$ is an isolated extreme ray of $K$. Then there exists $\beta > 0$ such that for every vector $z \in \mathbb R^n$, not proportional to $x_1$ and such that $zz^T \in K$, we have $\sum_{i=2}^n (z^Tx_i)^2 > \beta (z^Tx_1)^2$.

Since the intersection ${\cal L}_n(H) \cap \tilde L$ does not contain a rank 1 matrix, $x_1x_1^T \in \tilde L$, and $y_iy_i^T \in {\cal L}_n(H)$, we have that $y_i$ is not proportional to $x_1$ for every $i = 1,\dots,n-1$. It follows that $\sum_{j=2}^n (y_i^Tx_j)^2 > \beta (y_i^Tx_1)^2$ for all $i = 1,\dots,n-1$. Therefore
\begin{equation} \label{isol_contra}
\sum_{j=2}^n \left( d_j + d_{n+1} c_j^2 \right) = \sum_{j=2}^n x_j^TMx_j = \sum_{i=1}^{n-1} \sum_{j=2}^n (x_j^Ty_i)^2 > \beta \sum_{i=1}^{n-1} (y_i^Tx_1)^2 = \beta x_1^TMx_1 = \beta \left( d_1 + d_{n+1} c_1^2 \right).
\end{equation}
Fix now $d_2,\dots,d_n$ and let $d_1 \to +\infty$. Then $d_{n+1} \to -\left(\sum_{i=2}^n d_i^{-1}c_i^2\right)^{-1}$, and the leftmost term in \eqref{isol_contra} tends to a finite value. On the other hand, the rightmost term in \eqref{isol_contra} tends to $+\infty$, leading to a contradiction.

For the other extreme rays of $K$ the reasoning is similar after an appropriate permutation of the MLD set $\{R_1,\dots,R_{n+1}\}$.
\end{proof}

{\corollary \label{cor_dim_mld} Let $k \geq 2$ and let $K$ be a ROG cone possessing an MLD set $\{R_1,\dots,R_{k+1}\}$ of extreme rays. Then the following holds:

i) the dimension and degree of $K$ satisfy $\dim K \geq 2k - 1$, $\deg K \geq k$;

ii) the extreme rays $R_1,\dots,R_{k+1}$ of $K$ are not isolated. }

\begin{proof}
Represent $K$ as a linear section of ${\cal S}_+^n$ for some $n$, and let the rank 1 matrices $X_i = x_ix_i^T$ be generators of the extreme rays $R_i$, $i = 1,\dots,k+1$. Then the set $\{x_1,\dots,x_{k+1}\} \subset \mathbb R^n$ is minimally linearly dependent. By Lemma \ref{mld_char} the linear span $H$ of the vectors $x_1,\dots,x_{k+1}$ is a subspace of dimension $k$. Then $K_H = {\cal L}_n(H) \cap K$ is a face of $K$ and hence a ROG cone by Lemma \ref{face_ROG}. Moreover, $K_H$ contains the rank 1 matrices $x_1x_1^T,\dots,x_{k+1}x_{k+1}^T$ and is of degree $k$. In particular, the set $\{R_1,\dots,R_{k+1}\}$ is also an MLD set of extreme rays for $K_H$.

Applying Lemma \ref{dim_mld} to the cone $K_H$, we see that $K_H$ is simple and the extreme rays $R_1,\dots,R_{k+1}$ of $K_H$ are not isolated for all $i = 1,\dots,k+1$. By Theorem \ref{thm_lower} we have $\dim K_H \geq 2k-1$. But $\dim K \geq \dim K_H$, $\deg K \geq \deg K_H$, and every extreme ray of $K_H$ is also an extreme ray of $K$. The claim of the corollary now easily follows.
\end{proof}

{\corollary \label{cor_split_isol} Let $K$ be a ROG cone of degree $n$, and let $R$ be an isolated extreme ray of $K$. Then $K$ can be represented as a direct sum $K' \oplus {\cal S}_+^1$, where $K'$ is a ROG cone of degree $n-1$, such that the extreme ray $R$ is given by the set $\{0\} \oplus {\cal S}_+^1$. }

\begin{proof}
Represent $K$ as a linear section of the cone ${\cal S}_+^n$, and let $x \in \mathbb R^n$ be such that $X = xx^T$ generates the isolated extreme ray $R$ of $K$.

Define the set $S = \{ y \in \mathbb R^n \,|\, yy^T \in K \}$ and note that $x \in S$. By virtue of Corollary \ref{cor_dim_mld} the vector $x$ cannot be contained in a minimally linearly dependent subset of $S$ of cardinality at least 3. By Lemma \ref{mld_alternative} there exists a subspace $H \subset \mathbb R^n$ of dimension $n-1$ such that $x \not\in H$ and $S \subset H \cup \spa\{x\}$.

Hence $\spa K = \spa\{ yy^T \,|\, y \in S \} \subset {\cal L}_n(H) + \spa R$, and by Lemma \ref{lem_product} we have $K = K' + R$, where $K' = K \cap {\cal L}_n(H)$ is the face of $K$ generated by $H$, and the sum is isomorphic to the direct sum of the summands. By Lemma \ref{lem_product1} the cone $K'$ has degree $n-1$. This completes the proof.
\end{proof}

{\theorem \label{thm_isol} Let $K$ be a ROG cone of degree $n$. Then the number of its isolated extreme rays does not exceed $n$. Let $R_1,\dots,R_k$ be the isolated extreme rays of $K$. Then $K$ is isomorphic to a direct sum $K' \oplus \mathbb R_+^k$, where $K'$ is a ROG cone of degree $n-k$ without isolated extreme rays, and the extreme rays $R_1,\dots,R_k$ correspond to the extreme rays of the summand $\mathbb R_+^k$. }

\begin{proof}
We prove the theorem by induction over $n$. If $n = 1$, then $K = \mathbb R_+$, and the assertion is evident. Suppose now that $n \geq 2$ and the assertion is proven for cones of degrees not exceeding $n-1$.

If $K$ has no isolated extreme ray, then the assertion of the theorem holds with $K' = K$.

Assume now that $R$ is an isolated extreme ray of $K$. By Corollary \ref{cor_split_isol} $K$ can be represented as a direct sum $K_1 \oplus \mathbb R_+$, where $K_1$ is a ROG cone of degree $n-1$. By the assumption of the induction, the number of isolated extreme rays of $K_1$ is finite and does not exceed $n-1$, let these be $\rho_2,\dots,\rho_{k'}$, $1 \leq k' \leq n$. Moreover, $K_1$ is isomorphic to a direct sum $K' \oplus \mathbb R_+^{k'-1}$, where $K'$ is a ROG cone of degree $n-k'$ without isolated extreme rays. It follows that $K \cong K' \oplus \mathbb R_+^{k'}$.

Now every extreme ray of the direct sum $K' \oplus \mathbb R_+^{k'}$ is either an extreme ray of the factor $K'$ or an extreme ray of the factor $\mathbb R_+^{k'}$, and it is isolated in the direct sum if and only if it is isolated in the factor. The extreme rays of $K'$ are not isolated in $K'$, and hence they are not isolated in $K$. The factor $\mathbb R_+^{k'}$ has $k'$ extreme rays, and all of them are isolated. These $k'$ rays hence exhaust the isolated extreme rays of $K' \oplus \mathbb R_+^{k'}$. It follows that $k'= k$ and the assertion of the theorem readily follows.
\end{proof}

The discrete and the continuous part of the set of extreme rays of $K$ thus generate separate factors of the cone $K$. The factor generated by the discrete part is isomorphic to the nonnegative orthant, with the discrete extreme rays of $K$ being its generators.

{\corollary \label{cor_simple_isol} Let $K$ be a simple ROG cone of degree $\deg K \geq 2$. Then $K$ has no isolated extreme rays. }

\begin{proof}
The corollary is an immediate consequence of Theorem \ref{thm_isol}.
\end{proof}

{\lemma \label{lem_tangent} Let $K \subset {\cal S}_+^n$ be a ROG cone, and let $x \in \mathbb R^n$ be such that the rank 1 matrix $xx^T$ generates an extreme ray of $K$ which is not isolated. Then there exists a vector $y \in \mathbb R^n$, linearly independent of $x$, such that $xy^T + yx^T \in \spa K$. }

\begin{proof}
Assume the conditions of the lemma. Then there exists a sequence $v_1,v_2,\dots$ of nonzero vectors in $\mathbb R^n$ such that $x^Tv_k = 0$, $(x+v_k)(x+v_k)^T \in K$ for all $k$, and $\lim_{k \to \infty} v_k = 0$. Set $y_k = \frac{v_k}{||v_k||}$. Then we have $\frac{(x+v_k)(x+v_k)^T - xx^T}{||v_k||} = xy_k^T + y_kx^T + \frac{v_kv_k^T}{||v_k||} \in \spa K$. Since $\lim_{k \to \infty} \frac{v_kv_k^T}{||v_k||} = 0$ and $\spa K$ is closed, we have $xy^T + yx^T \in \spa K$ for every accumulation point of the sequence $y_1,y_2,\dots$. But such accumulation points exist due to the compactness of the unit sphere, and every such accumulation point is orthogonal to $x$. This completes the proof.
\end{proof}

{\corollary \label{cor_tangent} Let $K \subset {\cal S}_+^n$ be a simple ROG cone of degree $\deg K \geq 2$. Then for every nonzero vector $x \in \mathbb R^n$ such that $xx^T \in K$ there exists a vector $y \in \mathbb R^n$, linearly independent of $x$, such that $xy^T + yx^T \in \spa K$. }

\begin{proof}
The corollary is an immediate consequence of Lemma \ref{lem_tangent} and Corollary \ref{cor_simple_isol}.
\end{proof}

{\lemma \label{face_codim2} Let $K$ be a simple ROG cone of degree $\deg K \geq 2$. If $K$ has a face $F \subset K$ such that $\dim K - \dim F = 2$, then $K$ is isomorphic to an intertwining of $F$ and ${\cal S}_+^2$. }

\begin{proof}
Assume the conditions of the lemma, and set $n = \deg K$, $k = \deg F$. Represent $K$ as a linear section of the cone ${\cal S}_+^n$, and let $X \in F$ be a matrix of maximal rank $k$. Denote the image of $X$ by $H$. Then $F = {\cal L}_n(H) \cap K$. By Corollary \ref{cor_diag2} there exist linearly independent vectors $r_{k+1},\dots,r_n \in \mathbb R^n$ such that $\mathbb R^n = \spa(H \cup \{r_{k+1},\dots,r_n\})$ and $r_jr_j^T \in K$, $j = k+1,\dots,n$. We obtain $\dim K \geq \dim F + \dim\spa\{r_{k+1}r_{k+1}^T,\dots,r_nr_n^T\} = (\dim K - 2) + (n-k)$. It follows that $k \geq n-2$. If $k = n-2$, then $\spa K = \spa F + \spa\{r_{n-1}r_{n-1}^T,r_nr_n^T\}$ and $K$ is isomorphic to the direct sum $F \oplus {\cal S}_+^1 \oplus {\cal S}_+^1$, contradicting the simplicity of $K$.

Hence $k = n - 1$. Set $x = r_n$ for simplicity of notation. By Corollary \ref{cor_tangent} there exists a nonzero vector $y \in H$ such that $xy^T + yx^T \in \spa K$.

Since the codimension of $F$ in $K$ is two, we have $\spa K = \spa F \oplus \spa xx^T \oplus \spa (xy^T+yx^T)$. Since $K$ is simple, there exists a vector $z \in \mathbb R^n$ such that $z \not\in H \cup \spa\{x\}$ and $zz^T \in K$. Let $z = z_H + \beta x$ be the decomposition of $z$ corresponding to the direct sum decomposition $\mathbb R^n = H \oplus \spa\{x\}$. Then $z_H \not= 0$, $\beta \not= 0$, $zz^T = z_Hz_H^T + \beta(z_Hx^T + xz_H^T) + \beta^2 xx^T$. On the other hand, we have the decomposition $zz^T = Z_F + \alpha_1 xx^T + \alpha_2 (xy^T+yx^T)$, where $Z_F \in \spa F$.

Let $l$ be a linear form which is zero on $H$, but $l(x) = 1$. Contracting both decompositions of the rank 1 matrix $zz^T$ with $l$, we obtain $\beta z_H + \beta^2 x = \alpha_1 x + \alpha_2 y$, and hence $\alpha_1 = \beta^2$, $\beta z_H = \alpha_2 y$, $\alpha_2 \not= 0$, $Z_F = z_Hz_H^T = (\beta^{-1}\alpha_2)^2 yy^T \in \spa F$.

Hence $yy^T \in F \subset K$. Thus ${\cal L}_n(\spa\{x,y\}) \subset \spa K$ and ${\cal F}_n(\spa\{x,y\}) = {\cal L}_n(\spa\{x,y\}) \cap K$ is a face of $K$ which is isomorphic to ${\cal S}_+^2$. By construction $K$ is an intertwining of the faces $F$ and ${\cal F}_n(\spa\{x,y\})$, with the intersection $F \cap {\cal F}_n(\spa\{x,y\})$ generated by $yy^T$. This yields the assertion of the lemma.
\end{proof}

\section{Classification for small degrees} \label{sec:small}

In this section we classify all simple ROG cones $K$ of degree $n = \deg K \leq 4$ up to isomorphism. As we already noted in Subsection \ref{subs:continuous_family}, for degree 6 there exist infinitely many isomorphism classes of simple ROG cones. Whether the classification for degree 5 is finite remains an open question. Denote by $\Tri_+^n$ the cone of all tri-diagonal matrices in ${\cal S}_+^n$.

\subsection{Cones of degree $n \leq 3$}

For $n = 1$ the only ROG cone is ${\cal S}_+^1$.

For $n = 2$ we have the ROG cones ${\cal S}_+^1 \oplus {\cal S}_+^1$ and ${\cal S}_+^2$, of which only the latter is simple.

For $n = 3$ the only ROG cone of dimension 6 is ${\cal S}_+^3$, which is simple. By Theorem \ref{thm_lower} any other simple ROG cone must have dimension 5, i.e., is given by $K = \{ X \in {\cal S}_+^3 \,|\, \langle X,Q \rangle = 0 \}$ for some indefinite quadratic form $Q$. The isomorphism class of $K$ depends only on the signature of $Q$, and the forms $\pm Q$ define the same cone $K$. Moreover, every cone $K$ of this form is ROG by Corollary \ref{cor_codim1}. The possible isomorphism classes are hence given by the signatures $(++-)$ and $(+-0)$ of $Q$. It is easily seen that the corresponding ROG cones are isomorphic to $\Han_+^3$ and the full extension of ${\cal S}_+^1 \oplus {\cal S}_+^1$, respectively. The latter cone is isomorphic to $\Tri_+^3$.

\subsection{Cones of degree 4 and codimension $d \leq 2$}

Let $K$ be a ROG cone of degree $\deg K = 4$.

If $\dim K = 10$, then $K \simeq {\cal S}_+^4$.

If $\dim K = 9$, then $K$ is of the form $\{ X \in {\cal S}_+^4 \,|\, \langle X,Q \rangle = 0 \}$ for some indefinite quadratic form $Q$. As in the case $n = 3$, the isomorphism class of $K$ is defined by the signature of $Q$, where $\pm Q$ yield the same cone $K$. The possible isomorphism classes of $K$ are then defined by the signatures $(+-00)$, $(++-0)$, $(++--)$, and $(+++-)$ of $Q$. In the first two cases $K$ is a full extension of ${\cal S}_+^1 \oplus {\cal S}_+^1$ and $\Han_+^3$, respectively. In the third case $K$ is isomorphic to the cone of positive semi-definite $4 \times 4$ block-Hankel matrices $\Han_+^{2,2}$. It can be interpreted as the moment cone of the homogeneous biquadratic forms on $\mathbb R^2 \times \mathbb R^2$. It is not hard to see that all four isomorphism classes consist of simple cones.

Let $\dim K = 8$. If $K$ is simple, then by Theorem \ref{class_codim2} it is isomorphic to a full extension of ${\cal S}_+^1 \oplus {\cal S}_+^2$.

By Theorem \ref{thm_lower} any other simple ROG cone of degree 4 must have dimension 7.

\subsection{Cones of degree 4 and dimension 7}

The simple ROG cones of degree 4 and dimension 7 are somewhat more difficult to classify. We shall first consider a number of special cases and then show that the general case can be reduced to one of these special cases. The most complicated case is that of cones isomorphic to the $4 \times 4$ positive semi-definite Hankel matrices, its consideration can be found in the Appendix.

{\lemma \label{doubleS2} Let $K \subset {\cal S}_+^4$ be a simple ROG cone of dimension 7 and degree 4. Suppose that the subspace of block-diagonal matrices consisting of two blocks of size $2 \times 2$ each is contained in $\spa K$. Then $K$ is isomorphic to the cone $\Tri_+^4$ of positive semi-definite tri-diagonal matrices. }

\begin{proof}
The subspace of block-diagonal matrices as defined in the formulation of the lemma is 6-dimensional. Hence there exist scalars $a_{13},a_{14},a_{23},a_{24}$, not all equal zero, such that the linear span of $K$ is given by all matrices of the form
\[ A = \sum_{i=1}^7 \alpha_i A_i = \begin{pmatrix} \alpha_1 & \alpha_2 & \alpha_7 a_{13} & \alpha_7 a_{14} \\ \alpha_2 & \alpha_3 & \alpha_7 a_{23} & \alpha_7 a_{24} \\ \alpha_7 a_{13} & \alpha_7 a_{23} & \alpha_4 & \alpha_5 \\ \alpha_7 a_{14} & \alpha_7 a_{24} & \alpha_5 & \alpha_6 \end{pmatrix},\quad \alpha_1,\dots,\alpha_7 \in \mathbb R,
\]
where the matrices $A_1,\dots,A_7$ are defined by the above identity. Since $K$ is ROG, there exists a rank 1 matrix $Z = zz^T = \sum_{i=1}^7 \zeta_i A_i$ with $\zeta_7 \not= 0$. Its upper right $2 \times 2$ block is also rank 1. Hence $a_{13}a_{24} = a_{14}a_{23}$ and there exist angles $\varphi_1,\varphi_2$ and a positive scalar $r$ such that $a_{13} = r\cos\varphi_1\cos\varphi_2$, $a_{14} = r\cos\varphi_1\sin\varphi_2$, $a_{23} = r\sin\varphi_1\cos\varphi_2$, $a_{24} = r\sin\varphi_1\sin\varphi_2$. Define a basis of $\mathbb R^4$ by the vectors $x_1 = (-\sin\varphi_1,\cos\varphi_1,0,0)^T$, $x_2 = (\cos\varphi_1,\sin\varphi_1,0,0)^T$, $x_3 = (0,0,\cos\varphi_2,\sin\varphi_2)^T$, $x_4 = (0,0,-\sin\varphi_2,\cos\varphi_2)^T$. In the coordinates given by this basis $K$ equals the cone $\Tri_+^4$, which proves our claim.
\end{proof}

{\lemma \label{lem_row_simple} Let $K \subset {\cal S}_+^n$ be a ROG cone, let $e_1,\dots,e_n$ be the canonical basis vectors of $\mathbb R^n$, and let $y = (0,y_2,\dots,y_n)^T \in \mathbb R^n$ be a vector such that $y_2,\dots,y_n \not= 0$. If $e_1e_1^T,\dots,e_ne_n^T,e_1y^T+ye_1^T \in \spa K$, then $K$ is simple and $\dim K \geq 2n-1$. }

\begin{proof}
Suppose for the sake of contradiction that $K$ is not simple. Then there exists a nontrivial direct sum decomposition $\mathbb R^n = H_1 \oplus H_2$ such that for every rank 1 matrix $xx^T \in K$ we have either $x \in H_1$ or $x \in H_2$. Hence $e_i \in H_1 \cup H_2$ for $i = 1,\dots,n$. It follows that $H_1,H_2$ are spanned by complementary subsets of the canonical basis of $\mathbb R^n$. Hence there exists a permutation of the basis vectors such that in the corresponding coordinate system every matrix in $K$, and hence also in $\spa K$, becomes block-diagonal with a nontrivial block structure. But this is in contradiction with the assumption $e_1y^T+ye_1^T \in \spa K$. Hence $K$ must be simple.

Since the identity matrix is an element of $K$, we have $\deg K = n$. The bound on the dimension now follows from Theorem \ref{thm_lower}.
\end{proof}

{\corollary \label{cor_5_intertw} Let $K \subset {\cal S}_+^4$ be a simple ROG cone of dimension 7 and degree 4. Suppose there exist linearly independent vectors $z_1,z_2,z_3 \in \mathbb R^4$ and nonzero scalars $\alpha,\beta$ such that $z_1z_1^T,z_2z_2^T,z_3z_3^T,\alpha(z_1z_2^T+z_2z_1^T)+\beta(z_1z_3^T+z_3z_1^T) \in \spa K$. Then $K$ is isomorphic to either $\Tri_+^4$, or the full extension of ${\cal S}_+^1 \oplus {\cal S}_+^1 \oplus {\cal S}_+^1$, or an intertwining of $\Han_+^3$ and ${\cal S}_+^2$. }

\begin{proof}
Denote by $H \subset \mathbb R^4$ the hyperplane spanned by $z_1,z_2,z_3$. The face $F = {\cal L}_4(H) \cap K$ of $K$ is a ROG cone by Lemma \ref{face_ROG}. Applying Lemma \ref{lem_row_simple} to $F$, we obtain $\dim F \geq 5$.

However, $\dim F \not= 6$, because a ROG cone $K$ with a face $F$ of codimension 1 in $K$ is isomorphic to $F \oplus {\cal S}_+^1$ and hence not simple. Therefore $F$ has codimension 2 in $K$. By Lemma \ref{face_codim2} the cone $K$ is isomorphic to an intertwining of $F$ with ${\cal S}_+^2$.

In the previous subsection we established that a simple ROG cone of dimension 5 and degree 3 is isomorphic to either $\Han_+^3$ or $\Tri_+^3$. If $F \cong \Han_+^3$, then $K$ is isomorphic to an intertwining of $\Han_+^3$ and ${\cal S}_+^2$. If $F \cong \Tri_+^3$, then there exist two possibilities for $K$, because $\Tri_+^3$ has two non-isomorphic types of extreme rays. It is not hard to see that an intertwining of $\Tri_+^3$ with ${\cal S}_+^2$ along these two types of extreme rays leads to cones which are isomorphic to $\Tri_+^4$ or the full extension of ${\cal S}_+^1 \oplus {\cal S}_+^1 \oplus {\cal S}_+^1$, respectively.
\end{proof}

{\lemma \label{simpleS2} Let $K \subset {\cal S}_+^4$ be a simple ROG cone of dimension 7 and degree 4. Suppose that $K$ has a face which is isomorphic to ${\cal S}_+^2$. Then $K$ fulfills the conditions of Corollary \ref{cor_5_intertw}. }

\begin{proof}
By assumption there exist linearly independent vectors $x_1,x_2 \in \mathbb R^4$ such that $x_1x_1^T,x_2x_2^T,x_1x_2^T+x_2x_1^T \in \spa K$. By Corollary \ref{cor_diag2} we may complete $x_1,x_2$ with vectors $x_3,x_4$ to a basis of $\mathbb R^4$ such that $x_3x_3^T,x_4x_4^T \in K$. Pass to the coordinate system defined by this basis. By Corollary \ref{cor_tangent} there exists a nonzero vector $y = (y_1,y_2,0,y_4)^T$ such that $x_3y^T + yx_3^T \in \spa K$.

If $y_1 = y_2 = 0$, then $x_3x_4^T+x_4x_3^T \in \spa K$ and $K$ fulfills the conditions of Lemma \ref{doubleS2}. Hence $K$ is isomorphic to $\Tri_+^4$. The claim of the lemma then immediately follows in this case.

Suppose now that $y_1,y_2$ are not simultaneously zero.

Let us first consider the case $y_4 = 0$. Let $F = {\cal L}_4(\spa\{x_1,x_2,x_3\}) \cap K$ be the face of $K$ which consists of matrices $X \in K$ whose last column vanishes. Then $x_1x_1^T,x_2x_2^T,x_1x_2^T+x_2x_1^T,x_3x_3^T,x_3y^T+yx_3^T \in \spa F$, and $\dim F \geq 5$. Since $\dim F = 6$ is not possible by the simplicity of $K$, we then must have $\spa F = \spa \{x_1x_1^T,x_2x_2^T,x_1x_2^T+x_2x_1^T,x_3x_3^T,x_3y^T+yx_3^T\}$. It follows that $F$ is isomorphic to $\Tri_+^3$. From Lemma \ref{face_codim2} it follows that $K$ is isomorphic to an intertwining of $\Tri_+^3$ and ${\cal S}_+^2$, which proves the claim of the lemma in this case.

Suppose now that $y_4 \not= 0$. Define the nonzero vector $z_3 = (y_1,y_2,0,0)^T$. Then $z_3z_3^T \in K$ and $y = \alpha x_4 + \beta z_3$ with $\alpha = y_4$ and $\beta = 1$. The linearly independent vectors $z_1 = x_3$, $z_2 = x_4$, $z_3$, and scalars $\alpha,\beta$ then satisfy the conditions of Corollary \ref{cor_5_intertw}.
\end{proof}

We are now in a position to consider the general case.

{\theorem Let $K$ be a simple ROG cone of degree $\deg K = 4$ and dimension $\dim K = 7$. Then $K$ is isomorphic to either $\Tri_+^4$, or the full extension of ${\cal S}_+^1 \oplus {\cal S}_+^1 \oplus {\cal S}_+^1$, or an intertwining of $\Han_+^3$ and ${\cal S}_+^2$, or $\Han_+^4$. }

\begin{proof}
Let $K \subset {\cal S}_+^4$ be a simple ROG cone of degree 4 and dimension 7. By Corollary \ref{cor_diag2} there exist linearly independent vectors $x_1,x_2,x_3,x_4$ such that $x_ix_i^T \in K$, $i = 1,\dots,4$. Pass to the coordinate system defined by the basis $\{x_1,x_2,x_3,x_4\}$. Then all diagonal matrices are in the linear span of $K$. Moreover, by Corollary \ref{cor_tangent} there exist nonzero vectors $y_i = (y_{i1},y_{i2},y_{i3},y_{i4})^T$ such that $x_i^Ty_i = y_{ii} = 0$ and $x_iy_i^T + y_ix_i^T \in \spa K$, $i = 1,2,3,4$. Therefore $\spa K$ contains all matrices of the form
\begin{equation} \label{matrix_alpha}
\begin{pmatrix} \alpha_1 & \alpha_5 y_{12} + \alpha_6 y_{21} & \alpha_5 y_{13} + \alpha_7 y_{31} & \alpha_5 y_{14} + \alpha_8 y_{41} \\ \alpha_5 y_{12} + \alpha_6 y_{21} & \alpha_2 & \alpha_6 y_{23} + \alpha_7 y_{32} & \alpha_6 y_{24} + \alpha_8 y_{42} \\ \alpha_5 y_{13} + \alpha_7 y_{31} & \alpha_6 y_{23} + \alpha_7 y_{32} & \alpha_3 & \alpha_7 y_{34} + \alpha_8 y_{43} \\ \alpha_5 y_{14} + \alpha_8 y_{41} & \alpha_6 y_{24} + \alpha_8 y_{42} & \alpha_7 y_{34} + \alpha_8 y_{43} & \alpha_4 \end{pmatrix},\qquad \alpha_1,\dots,\alpha_8 \in \mathbb R.
\end{equation}
Since the dimension of $\spa K$ is 7, the matrices at the coefficients $\alpha_1,\dots,\alpha_8$ must be linearly dependent. This is equivalent to the condition that the matrix
\begin{equation} \label{matrixY}
Y = \begin{pmatrix} y_{12} & y_{21} & 0 & 0 \\ y_{13} & 0 & y_{31} & 0 \\ y_{14} & 0 & 0 & y_{41} \\ 0 & y_{23} & y_{32} & 0 \\ 0 & y_{24} & 0 & y_{42} \\ 0 & 0 & y_{34} & y_{43} \end{pmatrix}
\end{equation}
is rank-deficient, $\rk Y \leq 3$. Here the rows of $Y$ correspond to the elements $(1,2)$, $(1,3)$, $(1,4)$, $(2,3)$, $(2,4)$, $(3,4)$ of \eqref{matrix_alpha}, respectively, and the columns to the expressions at the coefficients $\alpha_5,\dots,\alpha_8$, respectively. By construction every column of $Y$ is nonzero.

If there exists a column of $Y$ with exactly one nonzero element, let it be $y_{ij}$, then $x_ix_i^T,x_jx_j^T,x_ix_j^T+x_jx_i^T \in \spa K$ and $K$ has a face which is isomorphic to ${\cal S}_+^2$. By Lemma \ref{simpleS2} the cone $K$ is then isomorphic to either $\Tri_+^4$, or the full extension of ${\cal S}_+^1 \oplus {\cal S}_+^1 \oplus {\cal S}_+^1$, or an intertwining of $\Han_+^3$ and ${\cal S}_+^2$.

If there exists a column of $Y$ with exactly two nonzero elements, let them be $y_{ij},y_{ik}$, then the linearly independent vectors $z_1 = x_i$, $z_2 = x_j$, $z_3 = x_k$ and scalars $\alpha = y_{ij}$, $\beta = y_{ik}$ satisfy the conditions of Corollary \ref{cor_5_intertw}, and $K$ is again isomorphic to one of the aforementioned cones.

Let us now assume that all elements $y_{ij}$ for $i \not= j$ are nonzero. Then $\rk Y = 3$, and the subspace spanned by the set $\{ x_ix_i^T, x_iy_i^T+y_ix_i^T \,|\, i = 1,2,3,4 \}$ has dimension 7. Since this subspace is contained in $\spa K$, it must actually equal $\spa K$. There exists a nonzero vector $\beta = (\beta_1,\beta_2,\beta_3,\beta_4)^T$ such that $Y\beta = 0$. It is easy to see that no three columns of $Y$ can be linearly dependent, and hence all elements $\beta_i$ are nonzero. By possibly multiplying $y_i$ by the nonzero constant $\beta_i$, we may assume without loss of generality that $\beta = (1,1,1,1)^T$. Then $y_{ij} = -y_{ji}$ for all $i,j = 1,\dots,4$, $i \not= j$.

It is not hard to check that $\spa K$ can then alternatively be written as the set $\{ X \in {\cal S}^4 \,|\, \langle X,Q_i \rangle = 0,\ i = 1,2,3 \}$, where the linearly independent quadratic forms $Q_1,Q_2,Q_3$ are given by
\[ {\scriptsize \begin{pmatrix} 0 & y_{13}y_{23} & -y_{12}y_{23} & 0 \\
  y_{13}y_{23} & 0 & y_{12}y_{13} & 0 \\
 -y_{12}y_{23} & y_{12}y_{13} & 0 & 0 \\
        0 & 0 & 0 & 0 \end{pmatrix},\
\begin{pmatrix}  0 & y_{14}y_{24} & 0 & -y_{12}y_{24} \\
  y_{14}y_{24} & 0 & 0 & y_{12}y_{14} \\
        0 &       0 & 0 &        0 \\
 -y_{12}y_{24} & y_{12}y_{14} & 0 &        0 \end{pmatrix},\
\begin{pmatrix} 0 & 0 & y_{14}y_{34} & -y_{13}y_{34} \\
        0 & 0 &       0 &        0 \\
  y_{14}y_{34} & 0 &       0 &  y_{13}y_{14} \\
 -y_{13}y_{34} & 0 & y_{13}y_{14} &        0 \end{pmatrix} },
\]
respectively.

The rank 1 matrices in $K$ are then given by $zz^T$ such that $z \not= 0$ and $z^TQ_iz = 0$ for $i = 1,2,3$. Let us determine the set of vectors $z = (z_1,z_2,z_3,z_4)^T$ which satisfy this quadratic system of equations. It is not hard to see that if a solution $z$ is not equal to a canonical basis vector, then all elements of $z$ are nonzero. For such $z$ the quadratic system can be written as
\begin{eqnarray}
y_{23}^{-1}z_1^{-1} - y_{13}^{-1}z_2^{-1} + y_{12}^{-1}z_3^{-1} &=& 0, \nonumber\\
y_{24}^{-1}z_1^{-1} - y_{14}^{-1}z_2^{-1} + y_{12}^{-1}z_4^{-1} &=& 0, \label{z_system} \\
y_{34}^{-1}z_1^{-1} - y_{14}^{-1}z_3^{-1} + y_{13}^{-1}z_4^{-1} &=& 0. \nonumber
\end{eqnarray}
This is a linear system in the unknowns $z_i^{-1}$. If the coefficient matrix of this system is full rank, then the solution $(z_1^{-1},\dots,z_4^{-1})$ is proportional to $(0,y_{12}^{-1},y_{13}^{-1},y_{14}^{-1})$ and does not correspond to a real vector $z$. In this case the only rank 1 matrices in the subspace $\spa K$ are the matrices $x_ix_i^T$, $i = 1,\dots,4$, and $K$ is not ROG.

Thus the coefficient matrix of system \eqref{z_system} is rank deficient. This implies that all $3 \times 3$ minors of this matrix vanish, which leads to the condition $y_{14}^{-1}y_{23}^{-1} - y_{13}^{-1}y_{24}^{-1} + y_{12}^{-1}y_{34}^{-1} = 0$. 
The general solution of system \eqref{z_system} is then given by
\[ \begin{pmatrix} z_1^{-1} \\ z_2^{-1} \\ z_3^{-1} \\ z_4^{-1} \end{pmatrix} = \gamma_1 \begin{pmatrix} y_{12}^{-2}+y_{13}^{-2}+y_{14}^{-2} \\ y_{13}^{-1}y_{23}^{-1}+y_{14}^{-1}y_{24}^{-1} \\ y_{14}^{-1}y_{34}^{-1}-y_{12}^{-1}y_{23}^{-1} \\ -y_{12}^{-1}y_{24}^{-1}-y_{13}^{-1}y_{34}^{-1} \end{pmatrix} + \gamma_2 \begin{pmatrix} 0 \\ y_{12}^{-1} \\ y_{13}^{-1} \\ y_{14}^{-1} \end{pmatrix},\qquad \gamma_1,\gamma_2 \in \mathbb R.
\]
It can be checked by direct calculation that none of the $2 \times 2$ minors of the $4 \times 2$ matrix composed of the two vectors at $\gamma_1,\gamma_2$, respectively, vanishes. Hence the 2-dimensional subspace of solutions of system \eqref{z_system} is transversal to all coordinate planes spanned by pairs of canonical basis vectors. By Lemma \ref{Hankel4} the cone $K$ is then isomorphic to $\Han_+^4$.
\end{proof}

\section{Complex and quaternionic Hermitian matrices} \label{sec:complex}

So far we considered spectrahedral cones defined as linear sections of the cone of positive semi-definite real symmetric matrices. The definition of ROG cones can be applied also to spectrahedral cones defined as linear sections of cones of complex Hermitian or quaternionic Hermitian matrices, or even more general, as linear sections of general symmetric cones, because the rank is well-defined for the elements of these cones. We shall now consider to which extent the results developed in the preceding sections carry over to the complex and quaternionic Hermitian case, and introduce a family of complex and quaternionic Hermitian ROG cones which does not exist in the real case.

The extension of Theorem \ref{theorem_iso} to the case of complex or quaternionic Hermitian matrices is not straightforward and remains open. Recall that the proof of Theorem \ref{theorem_iso} is based on Lemma \ref{lem_Pluecker}, which makes an assertion about the Pl\"ucker embedding of real Grassmanians. In the case of complex or quaternionic Grassmanians, the coefficients $\sigma_i$ in the formulation of this lemma have to be chosen not from the finite set $\{-1,+1\}$, but from the unit circle in the complex plane or from the quaternionic unit sphere $S^3$. But then the argument at the end of the proof of Lemma \ref{lem_Pluecker} is no more valid. For quaternionic matrices, the proof fails even earlier, because determinants of general quaternionic matrices and hence the Pl\"ucker embedding itself are not well-defined.

The results of Subsections \ref{subs:MDP} and \ref{subs:FS} carry over to the case of complex or quaternionic Hermitian matrices without changes. The same holds for Lemmas \ref{lem_product1}, \ref{lem_product}, and Corollary \ref{ROG_product}. Lemma \ref{simple_char} holds if we assume the direct sums in Definition \ref{def_simple} in the sense of Definition \ref{direct_sum_spectrahedral}. Note, however, that the space of quaternionic vectors of length $n$ is not a vector space. The subspaces $H_i$ in the decomposition in Lemmas \ref{lem_product} and \ref{simple_char} have to be assumed being invariant with respect to {\it right} multiplication by quaternionic coefficients. The results of Subsections \ref{subs:full_extensions}, \ref{subs:intertwinings}, \ref{subs:chordal} also carry over. In the complex analog of the construction in Subsection \ref{subs:continuous_family} we have to consider equivalence classes of quadruples of points in the complex projective plane. These are parameterized by the complex cross-ratio which leads to a family of isomorphism classes with a complex parameter. A generalization to the quaternionic case is also straightforward due to the recent development of a theory of the quaternionic cross-ratio in \cite{GwynneLibine12}. The complex and quaternionic analogs of the results in Subsection \ref{subs:codim1} are even stronger than in the real case due to the larger dimension of full faces of rank 2. For complex Hermitian matrices, every spectrahedral cone up to codimension 2 is ROG, for quaternionic Hermitian matrices up to codimension 4. Theorem \ref{thm_lower} holds without changes also for the complex and quaternionic cases. The results of Section \ref{sec:discrete} can be generalized to the complex and quaternionic case with the exception of Lemma \ref{face_codim2}. In the quaternionic case, minimally linearly dependent sets have to be defined with respect to the multiplication by quaternionic coefficients from the {\it right}. The classification of complex and quaternionic simple ROG spectrahedral cones is trivial up to degree 2, for degree 3 the situation is already more complicated than in the real case.

In the complex and the quaternionic case there exists one important class of ROG spectrahedral cones which is missing in the real case, namely the positive semi-definite block-Toeplitz matrices. For the complex case this result is widely known (see, e.g., \cite[Theorem 3.2]{Tismenetsky93}) and is equivalent to the matrix version of the Fej\'er-Riesz theorem \cite[p.118]{RosenblumRovnyak}. Below we shall provide a proof for the quaternionic case, which is valid with appropriate modifications also for the complex case. It is based on the following spectral factorization result for hyperunitary matrices, i.e., quaternionic invertible square matrices $U$ satisfying $U^* = U^{-1}$, where the asterisk denotes the conjugate transpose.

{\lemma \label{spectral} Let $U$ be a hyperunitary matrix. Then there exists another hyperunitary matrix $V$ of the same size and a diagonal matrix with diagonal entries on the quaternionic unit sphere $S^3$ such that $UV = VD$. In particular, for the columns $v_1,\dots,v_n$ of $V$ we have $Uv_k = v_kd_k$, where $d_1,\dots,d_n$ are the diagonal elements of $D$. }

\begin{proof}
There exist a hyperunitary matrix $V$ and an upper triangular matrix $T$ such that $U = VTV^*$ \cite{Brenner51}. Since $U,V$ are hyperunitary, $T$ must also be hyperunitary, $TT^* = I$. It follows that $T$ is diagonal with unit norm diagonal entries. We may hence set $D = T$ and obtain $U = VDV^*$. The claim now easily follows.
\end{proof}

Denote the cone of positive semi-definite Hermitian block-Toeplitz matrices consisting of $n \times n$ blocks of size $m \times m$ each by $\Toep_+^{n,m}$. Recall that Hermitian block-Toeplitz matrices have the form
\begin{equation} \label{toep}
T = \begin{pmatrix} m_0 & m_1^* & \ddots & m_{n-1}^* \\ m_1 & m_0 & \ddots & m_{n-2}^* \\ \ddots & \ddots & \ddots & \ddots \\ m_{n-1} & m_{n-2} & \ddots & m_0 \end{pmatrix},
\end{equation}
where $m_0$ is a Hermitian block and $m_1,\dots,m_{n-1}$ are general blocks of size $m \times m$. First we need a characterization of rank 1 matrices of this form.

{\lemma \label{lem:toep_rk1} A matrix of the form \eqref{toep} is positive semi-definite of rank 1 if and only if there exists a non-zero quaternionic vector $v \in \mathbb H^m$ and a quaternion $q$ with $|q| = 1$ such that
\begin{equation} \label{rk1_Toep}
T = \begin{pmatrix} v \\ vq \\ vq^2 \\ \vdots \\ vq^{n-1} \end{pmatrix} \begin{pmatrix} v \\ vq \\ vq^2 \\ \vdots \\ vq^{n-1} \end{pmatrix}^*.
\end{equation}
}

\begin{proof}
Let $T$ be as in \eqref{rk1_Toep}. Since $q^* = q^{-1}$, we get $(vq^k)(vq^l)^* = vq^{k-l}v^*$, and $T$ is of the form \eqref{toep} with $m_j = vq^jv^*$.

Let now $T \in \Toep_+^{n,m}$ be of rank 1. Then $T = uu^*$, where $u = (u_0^*,\dots,u_{n-1}^*)^*$ is a non-zero quaternionic vector partitioned in $n$ subvectors of length $m$ each. Since $m_0 = u_ku_k^*$ for all $k = 0,\dots,n-1$, every subvector has actually to be non-zero. Set $v = u_0$. Then $u_ku_k^* = vv^*$ yields $u_k = vq_k$ for every $k = 1,\dots,n-1$, where $q_k$ are unit norm quaternions. Set $q = q_1$. Then we have $vq_kv^* = u_ku_0^* = m_k = u_{k+1}u_1^* = vq_{k+1}q^{-1}v^*$ for all $k = 1,\dots,n-2$. This yields $q_{k+1} = q_kq$ and by induction $q_k = q^k$.
\end{proof}

{\lemma Let $T \in \Toep_+^{n,m}$ be of rank $N$. Then $T$ can be represented as a sum of $N$ rank 1 matrices in $\Toep_+^{n,m}$, and $\Toep_+^{n,m}$ is ROG. }

\begin{proof}
There exists a $nm \times N$ matrix $\tilde W$ such that $T = \tilde W\tilde W^*$. Partition $\tilde W$ into blocks $\tilde W_0,\dots,\tilde W_{n-1}$ of size $m \times N$. Define $(n-1)m \times N$ matrices $\tilde W_u,\tilde W_l$, such that $\tilde W_u$ is obtained from $\tilde W$ by removal of the block $\tilde W_{n-1}$, and $\tilde W_l$ is obtained by removal of the block $\tilde W_0$. By virtue of the block-Toeplitz structure of $T$ we have $\tilde W_u\tilde W_u^* = \tilde W_l\tilde W_l^*$. Hence there exists a $N \times N$ hyperunitary matrix $U$ such that $\tilde W_l = \tilde W_uU$. It follows that $\tilde W_k = \tilde W_{k-1}U$ for all $k = 1,\dots,n-1$, and by iterating $\tilde W_k = \tilde W_0U^k$.

By Lemma \ref{spectral} there exist hyperunitary matrices $D,V$, where $D$ is diagonal, such that $UV = VD$. The relation $\tilde W_k = \tilde W_0U^k$ can then be rewritten as $W_k = W_0D^k$, where we have defined $W_k = \tilde W_kV$, $k = 0,\dots,n-1$. Define also $W = \tilde WV$, then we have $T = WW^*$ and $W_0,\dots,W_{n-1}$ are the subblocks of $W$. Let $w_1,\dots,w_N$ be the columns of $W$, then we get $T = \sum_{j=1}^N w_jw_j^*$. By virtue of the relation $W_k = W_0D^k$, each of the rank 1 matrices $w_jw_j^*$ has structure \eqref{rk1_Toep} with $q$ being the $j$-th diagonal element of $D$ and $v$ the $j$-th column of $W_0$. The application of Lemma \ref{lem:toep_rk1} completes the proof.
\end{proof}

\section{Conclusions and open questions}

In this contribution we have defined and considered a special class of spectrahedral cones, the rank 1 generated cones. These cones are characterized by Property \ref{ROGproperty}. They have applications in optimization, namely for the approximation of difficult optimization problems by semi-definite programs, in the common case where the semi-definite program is obtained by dropping a rank 1 constraint on the matrix-valued decision variable. They are closely linked to the property of such a semi-definite relaxation being exact.

We provided many examples of ROG cones and several structural results. One of the main results has been that the geometry of a ROG cone as a convex conic subset of a real vector space uniquely determines its representation as a linear section of the positive semi-definite matrix cone, if this representation is required to satisfy Property \ref{ROGproperty}, up to isomorphism (Theorem \ref{theorem_iso}). In particular, every point of the cone has the same rank in every such representation. The rank also equals its Carath\'eodory number (Lemma \ref{Cara_X}). The Carath\'eodory number of the cone itself equals its degree as an algebraic interior (Corollary \ref{Cara_deg}).

There exist surprisingly many ROG cones. This is due to the fact that there are several non-trivial ways to construct ROG cones of higher degree out of ROG cones of lower degree, which we have called full extensions (Subsection \ref{subs:full_extensions}) and intertwinings (Subsection \ref{subs:intertwinings}). Besides, there is the obvious way of taking direct sums (Subsection \ref{subs:direct_sums}). Iterating these procedures, one may obtain families of mutually non-isomorphic ROG cones with arbitrarily many real parameters. One may call ROG cones that are neither direct sums nor intertwinings nor full extensions of other ROG cones {\it elementary}. Examples of elementary ROG cones are the cones of positive semi-definite block-Hankel matrices and the cones $K = \{ X \in {\cal S}_+^n \,|\, \langle X,Q \rangle = 0 \}$ of codimension 1 (Subsection \ref{subs:codim1}), where $Q$ is an indefinite non-degenerate quadratic form. Besides these infinite series of elementary ROG cones, there exists the exceptional moment cone of the ternary quartics of dimension 15 and degree 6. It is unknown whether for the real symmetric case there exist other elementary cones.

We classified the simple ROG cones, i.e., those not representable as non-trivial direct sums, up to degree 4. There are 1,1,3,10 equivalence classes of such cones for degrees 1,2,3,4, respectively, with respect to isomorphisms.

The set of extreme rays of a ROG cone is an intersection of quadrics and hence defines a real projective variety. The varieties defined by direct sums or intertwinings are finite unions of smaller projective varieties. The classification of the irreducible varieties defined by ROG cones is an open question. It would follow from a classification of the elementary ROG cones.

\appendix

\section{Pl\"ucker embeddings of real Grassmanians}

The purpose of this section is to provide Lemma \ref{Pluecker_tech}, which is needed for the proof of Theorem \ref{theorem_iso}. It turns out that Lemma \ref{Pluecker_tech} is essentially equivalent to a property of the Pl\"ucker embedding of real Grassmanians, which is stated below as Lemma \ref{lem_Pluecker}. However, we start with results on the rank 1 completion of partially specified matrices, which will be needed to prove Lemma \ref{lem_Pluecker}.

{\definition A real partially specified $n \times m$ matrix is defined by an index subset ${\cal P} \subset \{1,\dots,n\} \times \{1,\dots,m\}$, called a {\it pattern}, together with a collection of real numbers $(A_{ij})_{(i,j) \in {\cal P}}$. A {\it completion} of a partially specified matrix $({\cal P},(A_{ij})_{(i,j) \in {\cal P}})$ is a real $n \times m$ matrix $C$ such that $C_{ij} = A_{ij}$ for all $(i,j) \in {\cal P}$. }

We shall be concerned with the question when a partially specified matrix possesses a completion of rank 1. This problem has been solved in \cite{CJRW89}, see also \cite{HHW06}. In order to formulate the result, we need to define a weighted bipartite graph $G$ associated to the partially specified matrix. The two groups of vertices will be the row indices $1,\dots,n$ and the column indices $1,\dots,m$. The edges will be the elements of ${\cal P}$, with the weight of $(i,j)$ equal to $A_{ij}$.

{\lemma \cite[Theorem 5]{HHW06} \label{lem:rk1compl} A partially specified matrix $({\cal P},(A_{ij})_{(i,j) \in {\cal P}})$ has a rank 1 completion if and only if the following conditions are satisfied. If for some $(i,j) \in {\cal P}$ we have $A_{ij} = 0$, then either $A_{ij'} = 0$ for all $(i,j') \in {\cal P}$, or $A_{i'j} = 0$ for all $(i',j) \in {\cal P}$. Further, for every cycle $i_1$-$j_1$-$i_2$-$\cdots$-$i_k$-$j_k$-$i_1$, $k \geq 2$, of the bipartite graph $G$ corresponding to the partially specified matrix, where $i_1$ in the representation of the cycle is a row index, we have $\prod_{l=1}^k A_{i_lj_l} = A_{i_kj_1} \cdot \prod_{l=1}^{k-1} A_{i_{l+1}j_l}$. }

Note that the relation in the second condition of the lemma depends only on the cycle itself, but not on its starting point or on the direction in which the edges are traversed. Since the products in the lemma are multiplicative under the concatenation of paths \cite[p.2171]{HHW06}, we may also restrict the condition to prime cycles (i.e., chordless cycles where each vertex appears at most once). 

{\corollary \label{partially_spec} Let $A = ({\cal P},(A_{ij})_{(i,j) \in {\cal P}})$ be a partially specified matrix such that $A_{ij} = \pm1$ for all $(i,j) \in {\cal P}$, and $G$ the corresponding bipartite graph. Assume further that for every prime cycle $i_1$-$j_1$-$\cdots$-$j_k$-$i_1$ of $G$ with $k \geq 2$, where the representation of the cycle begins with a row index, we have $\prod_{l=1}^k A_{i_lj_l} = A_{i_kj_1} \cdot \prod_{l=1}^{k-1} A_{i_{l+1}j_l}$. Then there exists a rank 1 completion $C = ef^T$ of $A$ such that $e \in \{-1,+1\}^n$, $f \in \{-1,+1\}^m$. }

\begin{proof}
By Lemma \ref{lem:rk1compl} there exists a rank 1 completion $\tilde C = \tilde e\tilde f^T$ of $A$, where $\tilde e \in \mathbb R^n$, $\tilde f \in \mathbb R^m$. Suppose there exists an index $i$ such that $\tilde e_i = 0$. Then all elements of the $i$-th row of $\tilde C$ vanish, and all elements of this row are unspecified in $A$. We may then set $\tilde e_i = 1$ and $\tilde e\tilde f^T$ would still be a completion of $A$. Hence assume without loss of generality that all elements of $\tilde e$ are nonzero. In a similar manner, we may assume that the elements of $\tilde f$ are nonzero.

We then define the vectors $e \in \mathbb R^n$, $f \in \mathbb R^m$ element-wise by the signs of the elements of $\tilde e,\tilde f$, respectively. For every $(i,j) \in {\cal P}$ we then have $e_if_j = \frac{\tilde e_i\tilde f_j}{|\tilde e_i\tilde f_j|} = \frac{A_{ij}}{|A_{ij}|} = A_{ij}$, because $A_{ij} = \pm1$. It follows that $C = ef^T$ is also a completion of $A$.
\end{proof}

We now come to the Grassmanian $Gr(n,\mathbb R^m)$, i.e., the space of linear $n$-planes in $\mathbb R^m$. Fix a basis in $\mathbb R^m$. Then an $n$-plane $\Lambda$ can be represented by an $n$-tuple of linear independent vectors in $\mathbb R^m$, namely those spanning $\Lambda$. Let us treat these vectors as row vectors and stack them into an $n \times m$ matrix $M$. The matrix $M$ is determined only up to left multiplication by a nonsingular $n \times n$ matrix, reflecting the ambiguity in the choice of vectors spanning $\Lambda$. The {\it Pl\"ucker coordinate} $\Delta_{i_1\dots i_n}$ of $\Lambda$, where $1 \leq i_1 < \dots < i_n \leq m$, is defined as the determinant of the $n \times n$ submatrix formed of the columns $i_1,\dots,i_n$ of $M$. The vector $\Delta$ of all Pl\"ucker coordinates is determined by the $n$-plane $\Lambda$ up to multiplication by a nonzero constant and corresponds to a point in the projectivization $\mathbb P(\wedge^n \mathbb R^m)$ of the $n$-th exterior power of $\mathbb R^m$. The map $\Lambda \mapsto \Delta$ from $Gr(n,\mathbb R^m)$ to $\mathbb P(\wedge^n \mathbb R^m)$ is called the {\it Pl\"ucker embedding}. For a more detailed introduction into the subject see \cite[Chapter 7]{HodgePedoeI}.

{\lemma \label{lem_Pluecker} Let $\Lambda,\Lambda' \subset \mathbb R^m$ be two $n$-planes with Pl\"ucker coordinate vectors $\Delta,\Delta'$, respectively. Suppose there exists a positive constant $c$ such that $|\Delta_{i_1\dots i_n}| = c|\Delta'_{i_1\dots i_n}|$ for all $n$-tuples $(i_1,\dots,i_n)$. Then there exists a linear automorphism of $\mathbb R^m$, given by a diagonal coefficient matrix $\Sigma = \diag(\sigma_1,\dots,\sigma_m)$, where $\sigma_i \in \{-1,+1\}$ for all $i = 1,\dots,m$, which takes the $n$-plane $\Lambda$ to $\Lambda'$. }

\begin{proof}
Assume the conditions of the lemma. Let without restriction of generality $\Delta_{1\dots n} \not= 0$, then also $\Delta'_{1\dots n} \not= 0$. Otherwise we may permute the basis vectors of $\mathbb R^m$ to obtain these inequalities. Then we may choose the $n \times m$ matrix $M$ representing $\Lambda$ such that the first $n$ columns of $M$ form the identity matrix. Make a similar choice for the $n \times m$ matrix $M'$ representing $\Lambda'$. Then we have $\Delta_{1\dots n} = \Delta'_{1\dots n} = 1$ and hence $c = 1$ for this choice of $M,M'$. If $m = n$, then we may take $\Sigma$ as the identity matrix. Let $m > n$.

Let $k,l$ be indices such that $1 \leq k \leq n$, $n < l \leq m$. The determinant $\Delta_{1,\dots,k-1,k+1,\dots,n,l}$ is then given by $(-1)^{n-k}M_{kl}$. Likewise, $\Delta'_{1,\dots,k-1,k+1,\dots,n,l} = (-1)^{n-k}M'_{kl}$, and hence $|M_{kl}| = |M_{kl}'|$ by the assumption on $\Delta,\Delta'$. We then get $|M_{kl}| = |M_{kl}'|$ also for all $k = 1,\dots,n$, $l = 1,\dots,m$.

Let now ${\cal P}$ be the set of index pairs $(k,l)$ such that $M_{kl} \not= 0$, and set $A_{kl} = \frac{M'_{kl}}{M_{kl}} \in \{-1,+1\}$ for $(k,l) \in {\cal P}$. Then for every completion $C$ of the partially specified matrix $A = ({\cal P},(A_{kl})_{(k,l) \in {\cal P}})$ we have $M' = M \bullet C$, where $\bullet$ denotes the Hadamard matrix product.

We shall now show that the partially specified matrix $A$ satisfies the condition of Corollary \ref{partially_spec}. Let $i_1$-$j_1$-$\cdots$-$j_k$-$i_1$ be a prime cycle of the bipartite graph $G$ corresponding to $A$, where $k \geq 2$, $i_1,\dots,i_k$ are row indices, and $j_1,\dots,j_k$ are column indices. Since the cycle is prime, the row and column indices are mutually distinct. The $k \times k$ submatrix $\hat M$ of $M$ consisting of elements with row indices $i_1,\dots,i_k$ and column indices $j_1,\dots,j_k$ does not have any nonzero elements except those specified by the edges of the cycle, because any such element would render the cycle non-prime. In particular, every row and every column of $\hat M$ contains exactly two nonzero elements. The index set $\{j_1,\dots,j_k\}$ then has an empty intersection with $\{1,\dots,n\}$, because the first $n$ columns of $M$ contain strictly less than two nonzero elements each. Moreover, in the Leibniz formula for the determinant $\det\hat M$ only two products are nonzero, and the corresponding permutations are related by a cyclic permutation, which has sign $(-1)^{k-1}$. Therefore we have $|\det\hat M| = \left| \prod_{l=1}^k M_{i_lj_l} - (-1)^kM_{i_kj_1} \cdot \prod_{l=1}^{k-1} M_{i_{l+1}j_l} \right|$.

Consider the $n \times n$ submatrix of $M$ consisting of columns with indices in $(\{1,\dots,n\} \setminus \{i_1,\dots,i_k\}) \cup \{ j_1,\dots,j_k\}$. The determinant of this submatrix has absolute value $|\det\hat M|$ by construction. A similar formula holds for the absolute value of the determinant of the corresponding $n \times n$ submatrix of $M'$. By the assumption on $\Delta,\Delta'$ we then have
\[ \left| \prod_{l=1}^k M_{i_lj_l} - (-1)^kM_{i_kj_1} \cdot \prod_{l=1}^{k-1} M_{i_{l+1}j_l} \right| = \left| \prod_{l=1}^k M'_{i_lj_l} - (-1)^kM'_{i_kj_1} \cdot \prod_{l=1}^{k-1} M'_{i_{l+1}j_l} \right|.
\]
It follows that either
\[ \left(1 - \prod_{l=1}^k A_{i_lj_l}\right)\prod_{l=1}^k M_{i_lj_l} = \left(1 - A_{i_kj_1} \cdot \prod_{l=1}^{k-1} A_{i_{l+1}j_l}\right)(-1)^kM_{i_kj_1} \cdot \prod_{l=1}^{k-1} M_{i_{l+1}j_l}
\]
or
\[ \left(1 + \prod_{l=1}^k A_{i_lj_l}\right)\prod_{l=1}^k M_{i_lj_l} = \left(1 + A_{i_kj_1} \cdot \prod_{l=1}^{k-1} A_{i_{l+1}j_l}\right)(-1)^kM_{i_kj_1} \cdot \prod_{l=1}^{k-1} M_{i_{l+1}j_l}.
\]
Note that all the involved elements of $M$ are nonzero, while those of $A$ equal $\pm1$. The relation $\prod_{l=1}^k A_{i_lj_l} = -A_{i_kj_1} \cdot \prod_{l=1}^{k-1} A_{i_{l+1}j_l}$ would then imply that in each of the two equations above, one side is zero while the other is not. Therefore we must have $\prod_{l=1}^k A_{i_lj_l} = A_{i_kj_1} \cdot \prod_{l=1}^{k-1} A_{i_{l+1}j_l}$, and the condition in Corollary \ref{partially_spec} is fulfilled.

By this corollary there exists a rank 1 completion $C = ef^T$ of $A$ such that $e \in \{-1,+1\}^n$, $f \in \{-1,+1\}^m$. We then have $M' = M \bullet (ef^T) = \diag(e)\cdot M\cdot\diag(f)$. Setting $\Sigma = \diag(f)$ completes the proof.
\end{proof}

We now provide the technical result which is necessary for the proof of Theorem \ref{theorem_iso}.

{\lemma \label{Pluecker_tech} Let $x_1,\dots,x_m,y_1,\dots,y_m \in \mathbb R^n$ be such that $\spa\{x_1,\dots,x_m\} = \spa\{y_1,\dots,y_m\} = \mathbb R^n$. Denote by $L \subset {\cal S}^n$ the linear span of the set $\{x_1x_1^T,\dots,x_mx_m^T\}$ and assume there exists a linear map $\tilde f: L \to {\cal S}^n$ such that $\tilde f(x_ix_i^T) = y_iy_i^T$ for all $i = 1,\dots,m$. Assume further that there exists a positive constant $c > 0$ such that $\det Z = c\det\tilde f(Z)$ for all $Z \in L$. Then there exists a non-singular $n \times n$ matrix $S$ such that $\tilde f(Z) = SZS^T$ for all matrices $Z \in L$. }

\begin{proof}
Assemble the column vectors $x_i$ into an $n \times m$ matrix $X$ and the column vectors $y_i$ into an $n \times m$ matrix $Y$. By assumption these matrices have full row rank $n$. For mutually distinct indices $i_1,\dots,i_n \in \{1,\dots,m\}$, let $X_{i_1\dots i_n},Y_{i_1\dots i_n}$ be the $n \times n$ submatrices formed of the columns $i_1,\dots,i_n$ of $X,Y$, respectively. We have $\det(X_{i_1\dots i_n}X_{i_1\dots i_n}^T) = \det(\sum_{k=1}^n x_{i_k}x_{i_k}^T) = c\det(\sum_{k=1}^n y_{i_k}y_{i_k}^T) = c\det(Y_{i_1\dots i_n}Y_{i_1\dots i_n}^T)$, which implies $|\det X_{i_1\dots i_n}| = \sqrt{c}|\det Y_{i_1\dots i_n}|$.

Since the $n$-tuple $(i_1,\dots,i_n)$ was chosen arbitrarily, the $n$-planes spanned in $\mathbb R^m$ by the row vectors of $X,Y$, respectively, fulfill the conditions of Lemma \ref{lem_Pluecker}. By this lemma there exist a nonsingular $n \times n$ matrix $S$ and a diagonal matrix $\Sigma = \diag(\sigma_1,\dots,\sigma_m)$ with $\sigma_i \in \{-1,+1\}$ such that $Y = SX\Sigma$, or equivalently $y_i = \sigma_iSx_i$ for all $i = 1,\dots,m$. For every $i = 1,\dots,m$ we then have $\tilde f(x_ix_i^T) = y_iy_i^T = Sx_ix_i^TS^T$, and by linear extension we get the claim of the lemma.
\end{proof}

\section{Extreme elements of rank 2}

In this section we provide auxiliary results which are needed for the classification of ROG spectrahedral cones of codimension 2 in Subsection \ref{subs:codim2}. By virtue of Lemma \ref{lem_codim} these results allow in principle also a classification of ROG cones of codimensions 3 and 4, but the number and complexity of cases to be considered becomes prohibitive in the framework of this paper.

We first provide the following structural result on real symmetric matrix pencils. Recall that $y$ is called an eigenvector of the pencil $Q_1 + \lambda Q_2$ if the linear forms $Q_1y,Q_2y$ are linearly dependent.

{\lemma \label{pencils} Let $Q_1,Q_2$ be quadratic forms on $\mathbb R^n$ such that the pencil $Q_1 + \lambda Q_2$ possesses $n$ linearly independent real eigenvectors. Then there exists a direct sum decomposition $\mathbb R^n = H_0 \oplus H_1 \oplus \dots \oplus H_m$, non-degenerate quadratic forms $\Phi_k$ on $H_k$, $k = 1,\dots,m$, and mutually distinct angles $\varphi_1,\dots,\varphi_m \in [0,\pi)$ with the following properties. For every vector $x = \sum_{k=0}^m x_k$, where $x_k \in H_k$, we have $Q_1(x) = \sum_{k=1}^m \cos\varphi_k \Phi_k(x_k)$, $Q_2(x) = \sum_{k=1}^m \sin\varphi_k \Phi_k(x_k)$. Moreover, the set of real eigenvectors of the pencil $Q_1 + \lambda Q_2$ is given by the union $\bigcup_{k=1}^m (H_0 + H_k)$. }

\begin{proof}
We define the subspace $H_0$ as the intersection $\ker Q_1 \cap \ker Q_2$. For every real eigenvector $y \not\in H_0$ of the pencil $Q_1 + \lambda Q_2$, the linear span of the set $\{ Q_1y,Q_2y \}$ of linear forms has then dimension 1. Hence there exists a unique angle $\varphi(y) \in [0,\pi)$ such that $\sin\varphi(y) Q_1y - \cos\varphi(y) Q_2y = 0$.

By assumption we find linearly independent real eigenvectors $y_1,\dots,y_{n-\dim H_0}$ of the pencil $Q_1 + \lambda Q_2$ such that $\spa (H_0 \cup \{ y_1,\dots,y_{n-\dim H_0} \}) = \mathbb R^n$. Regroup these vectors into subsets $\{y_{11},\dots,y_{1d_1}\}$, $\dots$, $\{y_{m1},\dots,y_{md_m}\}$ such that $\varphi(y_{kl}) = \varphi_k$, $k = 1,\dots,m$, $l = 1,\dots,d_k$, where $\varphi_1,\dots,\varphi_m \in [0,\pi)$ are mutually distinct angles, and $d_k$ is the number of eigenvectors corresponding to angle $\varphi_k$. Define the subspace $H_k$ as the linear span of $y_{k1},\dots,y_{kd_k}$, $k = 1,\dots,m$. Then by construction we have that $H_0 \oplus H_1 \oplus \dots \oplus H_m$ is a direct sum decomposition of $\mathbb R^n$. Moreover, every vector $y \in H_k$ is an eigenvector and we have $\sin\varphi_k Q_1y - \cos\varphi_k Q_2y = 0$ for all $y \in H_k$, $k = 1,\dots,m$. It follows that there exist quadratic forms $\Phi_k$ on $H_k$, $k = 1,\dots,m$, such that $Q_1|_{H_k} = \cos\varphi_k \Phi_k$, $Q_2|_{H_k} = \sin\varphi_k \Phi_k$.

Let now $k,k' \in \{1,\dots,m\}$ be distinct indices and $y \in H_k$, $y' \in H_{k'}$ be arbitrary vectors. By construction we have $\sin\varphi_k Q_1y - \cos\varphi_k Q_2y = \sin\varphi_{k'} Q_1y' - \cos\varphi_{k'} Q_2y' = 0$. Therefore $\sin\varphi_k y^TQ_1y' - \cos\varphi_k y^TQ_2y' = \sin\varphi_{k'} y^TQ_1y' - \cos\varphi_{k'} y^TQ_2y' = 0$. But $\varphi_k,\varphi_{k'}$ are distinct, and thus this linear system on $y^TQ_iy'$ has only the trivial solution $y^TQ_1y' = y^TQ_2y' = 0$. The decomposition formulas $Q_1(x) = \sum_{k=1}^m \cos\varphi_k \Phi_k(x_k)$, $Q_2(x) = \sum_{k=1}^m \sin\varphi_k \Phi_k(x_k)$ now readily follow.

Let $1 \leq k \leq m$. Suppose there exists a vector $y \in H_k$ such that $\Phi_ky = 0$. Then we have $Q_1y = Q_2y = 0$, and $y \in H_0$. Thus $y = 0$, and the form $\Phi_k$ must be non-degenerate.

Let now $x = \sum_{k=0}^m x_k$ be a real eigenvector of the pencil $Q_1 + \lambda Q_2$, where $x_k \in H_k$. Then we have $\sin\varphi Q_1x - \cos\varphi Q_2x = 0$ for some angle $\varphi \in [0,\pi)$. Let $z_k \in H_k$, $k = 0,\dots,m$ be arbitrary vectors, and set $z = \sum_{k=0}^m z_k$. Then we have $x^TQ_1z = \sum_{k=1}^m \cos\varphi_k \Phi_k(x_k,z_k)$, $x^TQ_2z = \sum_{k=1}^m \sin\varphi_k \Phi_k(x_k,z_k)$. It follows that
\begin{eqnarray*}
0 &=& \sin\varphi x^TQ_1z - \cos\varphi x^TQ_2z = \sum_{k=1}^m (\sin\varphi\cos\varphi_k \Phi_k(x_k,z_k) - \cos\varphi\sin\varphi_k \Phi_k(x_k,z_k)) \\
&=& \sum_{k=1}^m \sin(\varphi-\varphi_k) \Phi_k(x_k,z_k).
\end{eqnarray*}
Here $\Phi_k(x_k,z_k) = \frac14(\Phi_k(x_k+z_k)-\Phi_k(x_k-z_k))$ is as usual the bilinear form defined by the quadratic form $\Phi_k$. Since this holds identically for all $z_k \in H_k$ and the forms $\Phi_k$ are non-degenerate, we must have either $\varphi = \varphi_k$ or $x_k = 0$ for each $k = 1,\dots,m$. Therefore $x \in H_0 + H_k$ for some $k$. On the other hand, every vector $x \in H_0+H_k$ is an eigenvector of the pencil $Q_1 + \lambda Q_2$, since it satisfies $\sin\varphi_k Q_1x - \cos\varphi_k Q_2x = 0$.
\end{proof}

We now come to spectrahedral cones $K \subset {\cal S}_+^n$ of codimension $d$. We represent these as in Lemma \ref{lem:codim_k} by linearly independent quadratic forms $Q_1,\dots,Q_d$ on $\mathbb R^n$, $K = \{ X \in {\cal S}_+^n \,|\, \langle X,Q_i \rangle = 0,\ i = 1,\dots,d \}$. We study the intersections of the spectrahedral cone $K$ with faces ${\cal F}_n(H)$ of ${\cal S}_+^n$ of rank not exceeding 2, i.e., where $H = \spa\{x,y\}$ for some vectors $x,y \in \mathbb R^n$.

{\lemma \label{rk2_lem1} Assume above notations. The face ${\cal F}_n(H) \cap K$ of $K$ is generated by an extreme element of rank 2 if and only if the $d \times 3$ matrix
\[ M(x,y) = \begin{pmatrix} x^TQ_1x & 2x^TQ_1y & y^TQ_1y \\ \vdots & \vdots & \vdots \\ x^TQ_dx & 2x^TQ_dy & y^TQ_dy \end{pmatrix}
\]
has rank $2$ and its kernel is generated by a vector $(a,b,c)^T \in \mathbb R^3$ such that the matrix $A = \begin{pmatrix} a & b \\ b & c \end{pmatrix}$ is definite. }

\begin{proof}
If $x,y$ are linearly dependent, then both the face ${\cal F}_n(H)$ and the matrix $M(x,y)$ have rank at most 1. Hence we may assume that $x,y$ are linearly independent.

The intersection ${\cal F}_n(H) \cap K$ is given by all matrices $X = axx^T + b(xy^T+yx^T) + cyy^T$ such that $\begin{pmatrix} a & b \\ b & c \end{pmatrix} \succeq 0$ and
\[ \begin{pmatrix} \langle X,Q_1 \rangle \\ \vdots \\ \langle X,Q_d \rangle \end{pmatrix} = M(x,y) \cdot \begin{pmatrix} a \\ b \\ c \end{pmatrix} = 0.
\]
The signature of $X = axx^T + b(xy^T+yx^T) + cyy^T$ equals the signature of $A$.

Suppose that $\rk M(x,y) = 2$ and the kernel of $M(x,y)$ is generated by a vector $(a,b,c)^T$ such that $A \succ 0$. Then $X = axx^T + b(xy^T+yx^T) + cyy^T$ is positive semi-definite of rank 2, and the face ${\cal F}_n(H) \cap K$ is generated by $X$, which proves the "if" direction.

Suppose now that $X = axx^T + b(xy^T+yx^T) + cyy^T$ is positive semi-definite of rank 2 and generates the face ${\cal F}_n(H) \cap K$. Then $(a,b,c)^T \in \ker M(x,y)$ and $A \succ 0$. Moreover, the dimension of $\ker M(x,y)$ is 1 and it must be generated by the vector $(a,b,c)^T$, otherwise we would have $\dim({\cal F}_n(H) \cap K) > 1$. It follows that $\rk M(x,y) = 2$, which proves the "only if" direction.
\end{proof}

{\lemma \label{rk2_lem2} Assume the notations of the previous lemma and set $d = 2$. The matrix $M(x,y)$ satisfies the conditions of the previous lemma if and only if the bi-quartic polynomial $p(x,y)$ given by
\[
(y^TQ_1y\cdot x^TQ_2x - x^TQ_1x\cdot y^TQ_2y)^2 - 4(x^TQ_1y\cdot y^TQ_2y - x^TQ_2y\cdot y^TQ_1y)(x^TQ_1x\cdot x^TQ_2y - x^TQ_1y\cdot x^TQ_2x)
\]
is negative on $x,y$. }

\begin{proof}
The matrix $A$ is definite if and only if $b^2 - ac < 0$, and $M(x,y)$ has full rank 2 if and only if the cross product
\[ \begin{pmatrix} x^TQ_1x \\ 2x^TQ_1y \\ y^TQ_1y \end{pmatrix} \times \begin{pmatrix} x^TQ_2x \\ 2x^TQ_2y \\ y^TQ_2y \end{pmatrix} = \begin{pmatrix} 2(x^TQ_1y\cdot y^TQ_2y - x^TQ_2y\cdot y^TQ_1y) \\ y^TQ_1y\cdot x^TQ_2x - x^TQ_1x\cdot y^TQ_2y \\ 2(x^TQ_1x\cdot x^TQ_2y - x^TQ_1y\cdot x^TQ_2x) \end{pmatrix}
\]
is nonzero. In this case the kernel of $M(x,y)$ is generated by this cross product, and hence $b^2 - ac < 0$ if and only if $p(x,y) < 0$.
\end{proof}

{\lemma \label{lem_S23} Assume the notations of the previous lemma and suppose that the polynomial $p(x,y)$ is nonnegative for all $x,y \in \mathbb R^n$. Suppose there exists $z \in \mathbb R^n$ such that $z^TQ_1z = z^TQ_2z = 0$ and the linear forms $q_1 = Q_1z$, $q_2 = Q_2z$ are linearly independent. Then there exists a linear form $u$ which is linearly independent from $q_1,q_2$ and such that $Q_1 = u \otimes q_1 + q_1 \otimes u$, $Q_2 = u \otimes q_2 + q_2 \otimes u$. }

\begin{proof}
By virtue of the condition $z^TQ_1z = z^TQ_2z = 0$ the nonnegative polynomial $p(x,y)$ vanishes for $x = z$ and all $y \in \mathbb R^n$. Therefore $\left.\frac{\partial p(x,y)}{\partial x}\right|_{x = z} = 0$ for all $y \in \mathbb R^n$. By virtue of $z^TQ_1z = z^TQ_2z = 0$, at $x = z$ this gradient is given by
\[ \left.\frac{\partial p(x,y)}{\partial x}\right|_{x = z} = -8(q_1^Ty\cdot y^TQ_2y - q_2^Ty\cdot y^TQ_1y)(q_2^Ty \cdot q_1 - q_1^Ty\cdot q_2) = 0.
\]
Since $q_1,q_2$ are linearly independent, the linear form $q_2^Ty \cdot q_1 - q_1^Ty\cdot q_2$ is nonzero if $q_1^Ty \not= 0$ or $q_2^Ty \not = 0$. Therefore $q_1^Ty\cdot y^TQ_2y = q_2^Ty\cdot y^TQ_1y$ for all such $y$, i.e., for a dense subset of $\mathbb R^n$. It follows that $q_1^Ty\cdot y^TQ_2y = q_2^Ty\cdot y^TQ_1y$ identically for all $y \in \mathbb R^n$.

In particular, for every $y \in \mathbb R^n$ such that $q_1^Ty = 0$, $q_2^Ty \not= 0$ we have $y^TQ_1y = 0$. The subset of such vectors $y$ is dense in the kernel of $q_1$, and hence $Q_1$ is zero on this kernel. It follows that there exists a linear form $u_1$ such that $Q_1 = q_1 \otimes u_1 + u_1 \otimes q_1$. In a similar manner, there exists a linear form $u_2$ such that $Q_2 = q_2 \otimes u_2 + u_2 \otimes q_2$. It follows that $q_1^Ty\cdot q_2^Ty \cdot u_2^Ty = q_2^Ty\cdot q_1^Ty \cdot u_1^Ty$ identically for all $y \in \mathbb R^n$. For all $y \in \mathbb R^n$ such that $q_1^Ty \not= 0$ and $q_2^Ty \not= 0$ it follows that $u_2^Ty = u_1^Ty$. Since the set of such vectors $y$ is dense in $\mathbb R^n$, we get that $u_1,u_2$ are equal to the same linear form $u$.

Note that $q_1^Tz = q_2^Tz = 0$ by assumption. We obtain $q_1 = Q_1z = q_1^Tz \cdot u + u^Tz \cdot q_1 = u^Tz \cdot q_1$, and hence $u^Tz = 1$. Therefore $u$ must be linearly independent of $q_1,q_2$.
\end{proof}

{\lemma \label{synthesis_rk2} Let $Q_1,Q_2$ be linearly independent quadratic forms on $\mathbb R^n$. Consider the spectrahedral cone $K = \{ X \in {\cal S}_+^n \,|\, \langle X,Q_1 \rangle = \langle X,Q_2 \rangle = 0 \}$. Suppose that $K$ does not have extreme elements of rank 2. Then one of the following conditions holds.

(i) For every $z \in \mathbb R^n$ such that $z^TQ_1z = z^TQ_2z = 0$ the linear forms $q_1 = Q_1z$, $q_2 = Q_2z$ are linearly dependent.

(ii) There exist linearly independent linear forms $q_1,q_2,u$ such that $Q_1 = u \otimes q_1 + q_1 \otimes u$, $Q_2 = u \otimes q_2 + q_2 \otimes u$. }

\begin{proof}

Suppose that condition (i) does not hold, then there exists $z \in \mathbb R^n$ such that $z^TQ_1z = z^TQ_2z = 0$ and the linear forms $q_1 = Q_1z$, $q_2 = Q_2z$ are linearly independent. Further, by Lemmas \ref{rk2_lem1}, \ref{rk2_lem2} the polynomial $p(x,y)$ defined in Lemma \ref{rk2_lem2} is nonnegative for all $x,y \in \mathbb R^n$. Hence the conditions of Lemma \ref{lem_S23} are fulfilled and condition (ii) holds.
\end{proof}

\section{Coordinate-free characterization of $\Han_+^4$}

In this section we provide an auxiliary result which is necessary for the classification of simple ROG cones of degree 4 and dimension 7.

{\lemma \label{Hankel4} Let $e_1,\dots,e_4$ be the canonical basis vectors of $\mathbb R^4$. Let $P \subset \mathbb R^4$ be a 2-dimensional subspace which is transversal to all coordinate planes spanned by pairs of basis vectors. Define the set of vectors
\begin{eqnarray*}
{\cal R} &=& \{ \alpha e_i \,|\, \alpha \in \mathbb R,\ i = 1,2,3,4 \} \quad \cup \\
&& \{ z = (z_1,z_2,z_3,z_4)^T \in \mathbb R^4 \,|\, z_i \not= 0\ \forall\ i = 1,\dots,4;\ (z_1^{-1},z_2^{-1},z_3^{-1},z_4^{-1})^T \in P \}.
\end{eqnarray*}
Let $L \subset {\cal S}^4$ be the linear span of the set $\{ zz^T \,|\, z \in {\cal R} \}$. Then $\dim L = 7$, and the spectrahedral cone $K = L \cap {\cal S}_+^4$ is isomorphic to the cone $\Han_+^4$ of positive semi-definite $4 \times 4$ Hankel matrices. }

\begin{proof}
Let $(r_1\cos\varphi_1,\dots,r_4\cos\varphi_4)^T,(r_1\sin\varphi_1,\dots,r_4\sin\varphi_4)^T \in \mathbb R^4$ be two linearly independent vectors spanning $P$. By the transversality property of $P$ all $2 \times 2$ minors of the $4 \times 2$ matrix composed of these vectors are nonzero. Hence the angles $\varphi_1,\dots,\varphi_4$ are mutually distinct modulo $\pi$, and the scalars $r_1,\dots,r_4$ are nonzero. We may also assume without loss of generality that none of the angles $\varphi_i$ is a multiple of $\pi$, otherwise we choose slightly different basis vectors in $P$.

For all $\xi \in [0,\pi)$ we then have that the vector $(r_1\sin(\varphi_1+\xi),\dots,r_4\sin(\varphi_4+\xi))^T$ is an element of $P$. More precisely, we get
\[ {\cal R} = \{ \alpha e_i \,|\, \alpha \in \mathbb R,\ i = 1,2,3,4 \} \cup \left\{ \alpha\left( \frac{1}{r_1\sin(\varphi_1+\xi)},\dots,\frac{1}{r_4\sin(\varphi_4+\xi)}\right)^T \,|\, \alpha \in \mathbb R,\ \xi \not= \varphi_i \mod\pi \right\}.
\]
Now set $\cos\xi = \frac{1-t^2}{1+t^2}$, $\sin\xi = \frac{2t}{1+t^2}$, and $s = t - \frac{1}{t}$. Then $\frac{1}{r_i\sin(\varphi_i+\xi)} = \frac{1+t^2}{r_it(2\cos\varphi_i-s\sin\varphi_i)}$. Define the vector $\mu(s) = \left( \frac{1}{r_1(2\cos\varphi_1-s\sin\varphi_1)},\dots,\frac{1}{r_4(2\cos\varphi_4-s\sin\varphi_4)}\right)^T$ for all $s \in \mathbb R$ except the values $s = 2\cot\varphi_i$, $i = 1,\dots,4$. We then get
\begin{eqnarray*}
{\cal R} &=& \{ \alpha e_i \,|\, \alpha \in \mathbb R,\ i = 1,2,3,4 \} \cup \{ \alpha\mu(s) \,|\, \alpha \in \mathbb R,\ s \in \mathbb R,\ s \not= 2\cot\varphi_i \} \cup \\
&& \cup \left\{ \alpha\left( \frac{1}{r_1\sin\varphi_1},\dots,\frac{1}{r_4\sin\varphi_4}\right)^T \,|\, \alpha \in \mathbb R \right\}.
\end{eqnarray*}
Multiplying the vector $\mu(s)$ by the common denominator of its elements, we obtain the vector
\begin{equation} \label{nus}
\nu(s) = \mu(s) \cdot \prod_{i=1}^4 (2\cos\varphi_i-s\sin\varphi_i) = \diag(r_1^{-1},r_2^{-1},r_3^{-1},r_4^{-1}) \cdot M \cdot \diag(8,-4,2,-1) \cdot \eta(s),
\end{equation}
where $\eta(s) = (1,s,s^2,s^3)^T$ and the matrix $M$ is given by
\[ {\scriptsize \begin{pmatrix}
\cos\varphi_2\cos\varphi_3\cos\varphi_4 & \sin\varphi_2\sin\varphi_3\sin\varphi_4 + \sin(\varphi_2+\varphi_3+\varphi_4) & \cos\varphi_2\cos\varphi_3\cos\varphi_4 - \cos(\varphi_2+\varphi_3+\varphi_4) & \sin\varphi_2\sin\varphi_3\sin\varphi_4 \\
\cos\varphi_1\cos\varphi_3\cos\varphi_4 & \sin\varphi_1\sin\varphi_3\sin\varphi_4 + \sin(\varphi_1+\varphi_3+\varphi_4) & \cos\varphi_1\cos\varphi_3\cos\varphi_4 - \cos(\varphi_1+\varphi_3+\varphi_4) & \sin\varphi_1\sin\varphi_3\sin\varphi_4 \\
\cos\varphi_1\cos\varphi_2\cos\varphi_4 & \sin\varphi_1\sin\varphi_2\sin\varphi_4 + \sin(\varphi_1+\varphi_2+\varphi_4) & \cos\varphi_1\cos\varphi_2\cos\varphi_4 - \cos(\varphi_1+\varphi_2+\varphi_4) & \sin\varphi_1\sin\varphi_2\sin\varphi_4 \\
\cos\varphi_1\cos\varphi_2\cos\varphi_3 & \sin\varphi_1\sin\varphi_2\sin\varphi_3 + \sin(\varphi_1+\varphi_2+\varphi_3) & \cos\varphi_1\cos\varphi_2\cos\varphi_3 - \cos(\varphi_1+\varphi_2+\varphi_3) & \sin\varphi_1\sin\varphi_2\sin\varphi_3
\end{pmatrix} }.
\]
Here for the calculus of $M$ we used the formulas
\begin{eqnarray*}
\sin\varphi_i\cos\varphi_j\cos\varphi_k+\sin\varphi_j\cos\varphi_i\cos\varphi_k+\sin\varphi_k\cos\varphi_i\cos\varphi_j &=& \sin\varphi_i\sin\varphi_j\sin\varphi_k + \sin(\varphi_i+\varphi_j+\varphi_k), \\ \sin\varphi_2\sin\varphi_3\cos\varphi_4+\sin\varphi_2\sin\varphi_4\cos\varphi_3+\sin\varphi_3\sin\varphi_4\cos\varphi_2 &=& \cos\varphi_i\cos\varphi_j\cos\varphi_k - \cos(\varphi_i+\varphi_j+\varphi_k).
\end{eqnarray*}
Note that the vector $\nu(s)$ can also be defined by the right-hand side of \eqref{nus} for $s = 2\cot\varphi_i$, and for this value of $s$ it is proportional to $e_i$. Defining $\eta(\infty) = e_4$ and $\nu(\infty) = \diag(r_1^{-1},r_2^{-1},r_3^{-1},r_4^{-1}) \cdot M \cdot \diag(8,-4,2,-1) \cdot \eta(\infty)$, we finally get
\[ {\cal R} = \{ \alpha\nu(s) \,|\, \alpha \in \mathbb R,\ s \in \mathbb R \cup \{\infty\} \}.
\]

A symbolic computation with a computer algebra system yields
\[ \det M = \sin(\varphi_1-\varphi_2)\sin(\varphi_1-\varphi_3)\sin(\varphi_1-\varphi_4)\sin(\varphi_2-\varphi_3)\sin(\varphi_2-\varphi_4)\sin(\varphi_3-\varphi_4) \not= 0.
\]
Hence the matrix product $\diag(r_1^{-1},r_2^{-1},r_3^{-1},r_4^{-1}) \cdot M \cdot \diag(8,-4,2,-1)$ is non-degenerate, and the subspace $L = \spa\{ zz^T \,|\, z = \nu(s),\ s \in \mathbb R \cup \{\infty\} \}$ is isomorphic to the subspace $L' = \spa\{ zz^T \,|\, z = \eta(s),\ s \in \mathbb R \cup \{\infty\} \}$.

The subspace $L'$, however, is the subspace of Hankel matrices in ${\cal S}^4$. The claim of the lemma now easily follows.
\end{proof}


\end{document}